\documentclass[11pt]{amsart}
\usepackage{amssymb,mathrsfs}
\title{Theory of the Siegel Modular Variety}

\begin{document}

\author{Jae-Hyun Yang}

\address{Department of Mathematics, Inha University,
Incheon 402-751, Korea}
\email{jhyang@inha.ac.kr }


\newtheorem{theorem}{Theorem}[section]
\newtheorem{lemma}{Lemma}[section]
\newtheorem{proposition}{Proposition}[section]
\newtheorem{remark}{Remark}[section]
\newtheorem{definition}{Definition}[section]

\renewcommand{\theequation}{\thesection.\arabic{equation}}
\renewcommand{\thetheorem}{\thesection.\arabic{theorem}}
\renewcommand{\thelemma}{\thesection.\arabic{lemma}}
\newcommand{\BR}{\mathbb R}
\newcommand{\BQ}{\mathbb Q}
\newcommand{\BT}{\mathbb T}
\newcommand{\BM}{\mathbb M}
\newcommand{\bn}{\bf n}
\def\charf {\mbox{{\text 1}\kern-.24em {\text l}}}
\newcommand{\BC}{\mathbb C}
\newcommand{\BZ}{\mathbb Z}

\thanks{\noindent{Subject Classification:} Primary 14K10\\
\indent Keywords and phrases: Siegel modular variety, Siegel
modular forms, abelian varieties, Satake parameters, lifting,
invariant holomorphic differential forms, proportionality theorem,
motives, cohomology}


\begin{abstract}
{In this paper, we discuss the theory of the Siegel modular
variety in the aspects of arithmetic and geometry. This article
covers the theory of Siegel modular forms, the Hecke theory, a
lifting of elliptic cusp forms, geometric properties of the Siegel
modular variety, (hypothetical) motives attached to Siegel modular
forms and a cohomology of the Siegel modular variety. }
\end{abstract}

\maketitle

\newcommand\tr{\triangleright}
\newcommand\al{\alpha}
\newcommand\be{\beta}
\newcommand\g{\gamma}
\newcommand\gh{\Cal G^J}
\newcommand\G{\Gamma}
\newcommand\de{\delta}
\newcommand\e{\epsilon}
\newcommand\z{\zeta}
\newcommand\vth{\vartheta}
\newcommand\vp{\varphi}
\newcommand\om{\omega}
\newcommand\p{\pi}
\newcommand\la{\lambda}
\newcommand\lb{\lbrace}
\newcommand\lk{\lbrack}
\newcommand\rb{\rbrace}
\newcommand\rk{\rbrack}
\newcommand\s{\sigma}
\newcommand\w{\wedge}
\newcommand\fgj{{\frak g}^J}
\newcommand\lrt{\longrightarrow}
\newcommand\lmt{\longmapsto}
\newcommand\lmk{(\lambda,\mu,\kappa)}
\newcommand\Om{\Omega}
\newcommand\ka{\kappa}
\newcommand\ba{\backslash}
\newcommand\ph{\phi}
\newcommand\M{{\Cal M}}
\newcommand\bA{\bold A}
\newcommand\bH{\bold H}
\newcommand\D{\Delta}

\newcommand\Hom{\text{Hom}}
\newcommand\cP{\Cal P}

\newcommand\cH{\Cal H}

\newcommand\pa{\partial}

\newcommand\pis{\pi i \sigma}
\newcommand\sd{\,\,{\vartriangleright}\kern -1.0ex{<}\,}
\newcommand\wt{\widetilde}
\newcommand\fg{\frak g}
\newcommand\fk{\frak k}
\newcommand\fp{\frak p}
\newcommand\fs{\frak s}
\newcommand\fh{\frak h}
\newcommand\Cal{\mathcal}

\newcommand\fn{{\frak n}}
\newcommand\fa{{\frak a}}
\newcommand\fm{{\frak m}}
\newcommand\fq{{\frak q}}
\newcommand\CP{{\mathcal P}_g}
\newcommand\Hgh{{\mathbb H}_g \times {\mathbb C}^{(h,g)}}
\newcommand\BD{\mathbb D}
\newcommand\BH{\mathbb H}
\newcommand\CCF{{\mathcal F}_g}
\newcommand\CM{{\mathcal M}}
\newcommand\Ggh{\Gamma_{g,h}}
\newcommand\Chg{{\mathbb C}^{(h,g)}}
\newcommand\Yd{{{\partial}\over {\partial Y}}}
\newcommand\Vd{{{\partial}\over {\partial V}}}

\newcommand\Ys{Y^{\ast}}
\newcommand\Vs{V^{\ast}}
\newcommand\LO{L_{\Omega}}
\newcommand\fac{{\frak a}_{\mathbb C}^{\ast}}

\begin{center}
\textit{To the memory of my mother}
\end{center}
\vskip 1cm

\centerline{\large \bf Table of Contents}

\vskip 0.5cm $ \qquad\qquad\qquad\qquad\textsf{\large \ 1.
Introduction}$\vskip 0.01cm

$\qquad\qquad\qquad\qquad \textsf{\large\ 2. Invariant Metrics and
Laplacians on Siegel Space }$\vskip 0.01cm

$ \qquad\qquad\qquad\qquad  \textsf{\large\ 3. Invariant
Differential Operators on Siegel Space }$\vskip 0.01cm

$ \qquad\qquad\qquad\qquad \textsf{\large\ 4. Siegel's Fundamental
Domain}$\vskip 0.01cm

$ \qquad\qquad\qquad\qquad  \textsf{\large\ 5. Siegel Modular
Forms}$\vskip 0.01cm

$ \qquad\qquad\qquad\qquad\qquad  \textsf{\ 5.1. Basic Properties
of Siegel Modular Forms}$\vskip 0.01cm

$ \qquad\qquad\qquad\qquad\qquad  \textsf{\ 5.2. The Siegel
Operator}$\vskip 0.01cm

$ \qquad\qquad\qquad\qquad\qquad  \textsf{\ 5.3. Construction of
Siegel Modular Forms }$\vskip 0.01cm

$ \qquad\qquad\qquad\qquad\qquad  \textsf{\ 5.4. Singular Modular
Forms}$\vskip 0.01cm

$ \qquad\qquad\qquad\qquad \textsf{\large\ 6. The Hecke Algebra
}$\vskip 0.01cm

$ \qquad\qquad\qquad\qquad\qquad  \textsf{\ 6.1. The Structure of
the Hecke Algebra }$\vskip 0.01cm
\par

$ \qquad\qquad\qquad\qquad\qquad  \textsf{\ 6.2. Action of the
Hecke Algebra on Siegel Modular Forms}$\vskip 0.01cm

$ \qquad\qquad\qquad\qquad \textsf{\large\ 7. Jacobi Forms
}$\vskip 0.01cm

$ \qquad\qquad\qquad\qquad  \textsf{\large\ 8. Lifting of Elliptic
Cusp Forms to Siegel Modular Forms }$\vskip 0.01cm

$ \qquad\qquad\qquad\qquad  \textsf{\large\ 9. Holomorphic
Differential Forms on Siegel Space }$\vskip 0.01cm

$ \qquad\qquad\qquad\qquad  \textsf{\large 10. Subvarieties of the
Siegel Modular Variety  }$\vskip 0.01cm

$ \qquad\qquad\qquad\qquad \textsf{\large 11. Proportionality
Theorem }$\vskip 0.01cm

$ \qquad\qquad\qquad\qquad  \textsf{\large 12. Motives and Siegel
Modular Forms }$\vskip 0.01cm

$ \qquad\qquad\qquad\qquad\textsf{\large 13. Remark on Cohomology of a Shimura Variety }$
\vskip 0.01cm

$ \qquad\qquad\qquad\qquad\textsf{\large References }$

\newpage

%
%
\begin{section}{{\bf Introduction}}
\setcounter{equation}{0} For a given fixed positive integer $g$,
we let
$${\mathbb H}_g=\,\{\,\Omega\in \BC^{(g,g)}\,|\ \Om=\,^t\Om,\ \ \ \text{Im}\,\Om>0\,\}$$
be the Siegel upper half plane of degree $g$ and let
$$Sp(g,\BR)=\{ M\in \BR^{(2g,2g)}\ \vert \ ^t\!MJ_gM= J_g\ \}$$
be the symplectic group of degree $g$, where $F^{(k,l)}$ denotes
the set of all $k\times l$ matrices with entries in a commutative
ring $F$ for two positive integers $k$ and $l$, $^t\!M$ denotes
the transposed matrix of a matrix $M$ and
$$J_g=\begin{pmatrix} 0&I_g\\
                   -I_g&0\end{pmatrix}.$$
Then $Sp(g,\BR)$ acts on $\BH_g$ transitively by
\begin{equation}
M\cdot\Om=(A\Om+B)(C\Om+D)^{-1},
\end{equation} where $M=\begin{pmatrix} A&B\\
C&D\end{pmatrix}\in Sp(g,\BR)$ and $\Om\in \BH_g.$ Let
$$\G_g=Sp(g,\BZ)=\left\{ \begin{pmatrix} A&B\\
C&D\end{pmatrix}\in Sp(g,\BR) \,\big| \ A,B,C,D\
\textrm{integral}\ \right\}$$ be the Siegel modular group of
degree $g$. This group acts on $\BH_g$ properly discontinuously.
C. L. Siegel investigated the geometry of $\BH_g$ and automorphic
forms on $\BH_g$ systematically. Siegel\,\cite{Si1} found a
fundamental domain ${\mathcal F}_g$ for $\G_g\ba\BH_g$ and
described it explicitly. Moreover he calculated the volume of
$\CCF.$ We also refer to \cite{Ig},\,\cite{M2},\,\cite{Si1} for
some details on $\CCF.$ Siegel's fundamental domain is now called
the Siegel modular variety and is usually denoted by ${\Cal A}_g$.
In fact, ${\Cal A}_g$ is one of the important arithmetic varieties
in the sense that it is regarded as the moduli of principally
polarized abelian varieties of dimension $g$. Suggested by Siegel,
I. Satake \cite{Sa1} found a canonical compactification, now
called the Satake compactification of ${\Cal A}_g$. Thereafter W.
Baily \cite{B1} proved that the Satake compactification of ${\Cal
A}_g$ is a normal projective variety. This work was generalized to
bounded symmetric domains by W. Baily and A. Borel \cite{BB}
around the 1960s. Some years later a theory of smooth
compactification of bounded symmetric domains was develpoed by
Mumford school \cite{AMRT}. G. Faltings and C.-L. Chai \cite{FC}
investigated the moduli of abelian varieties over the integers and
could give the analogue of the Eichler-Shimura theorem that
expresses Siegel modular forms in terms of the cohomology of local
systems on ${\Cal A}_g$. I want to emphasize that Siegel modular
forms play an important role in the theory of the arithmetic and
the geometry of the Siegel modular variety ${\Cal A}_g$.

The aim of this paper is to discuss a theory of the Siegel modular
variety in the aspects of arithmetic and geometry. Unfortunately
two important subjects, which are the theory of harmonic analysis
on the Siegel modular variety, and the Galois representations
associated to Siegel modular forms are not covered in this
article. These two topics shall be discussed in the near future in
the separate papers. This article is organized as follows. In
Section 2, we review the results of Siegel and Maass on invariant
metrics and their Laplacians on $\BH_g$. In Section 3, we
investigate differential operators on $\BH_g$ invariant under the
action (1.1). In Section 4, we review Siegel's fundamental domain
${\mathcal F}_g$ and expound the spectral theory of the abelian
variety $A_{\Omega}$ associated to an element $\Omega$ of
${\mathcal F}_g$. In Section 5, we review some properties of
vector valued Siegel modular forms, and also discuss construction
of Siegel modular forms and singular modular forms. In Section 6,
we review the structure of the Hecke algebra of the group
$GSp(g,{\mathbb Q})$ of symplectic similitudes and investigate the
action of the Hecke algebra on Siegel modular forms. In Section 7,
we briefly illustrate the basic notion of Jacobi forms which are
needed in the next section. We also give a short historical survey
on the theory of Jacobi forms. In Section 8, we deal with a
lifting of elliptic cusp forms to Siegel modular forms and give
some recent results on the lifts obtained by some people. A
lifting of modular forms plays an important role arithmetically
and geometrically. One of the interesting lifts is the so-called
Duke-Imamo${\check g}$lu-Ikeda lift. We discuss this lift in some
detail. In Section 9, we give a short survey of toroidal
compactifications of the Siegel modular variety ${\Cal A}_g$ and
illustrate a relationship between Siegel modular forms and
holomorphic differential forms on ${\mathcal A}_g.$ Siegel modular
forms related to holomorphic differential forms on ${\mathcal
A}_g$ play an important role in studying the geometry of
${\mathcal A}_g.$ In Section 10, We investigate the geometry of
subvarieties of the Siegel modular variety. Recently Grushevsky
and Lehavi \cite{G-L} announced that they proved that the Siegel
modular variety ${\mathcal A}_6$ of genus $6$ is of general type
after constructing a series of new effective geometric divisors on
${\mathcal A}_g.$ Before 2005 it had been known that ${\mathcal
A}_g$ is of general type for $g\geq 7$. In fact, in 1983
Mumford\,\cite{Mf4} proved that ${\mathcal A}_g$ is of general
type for $g\geq 7$. Nearly past twenty years nobody had known
whether ${\mathcal A}_6$ is of general type or not. In Section 11,
we formulate the proportionality theorem for an automorphic vector
bundle on the Siegel modular variety following the work of Mumford
(cf.\,\cite{Mf3}). In Section 12, we explain roughly Yoshida's
interesting results about the fundamental periods of a motive
attached to a Siegel modular form. These results are closely
related to Deligne's conjecture about critical values of an
$L$-function of a motive and the (pure or mixed) Hodge theory. In
the final section, we recall the definition of a Shimura variety
and give some remarks on the cohomology of Shimura varieties.\par
In person I am indebted to C. L. Siegel, one of the great
mathematicians of the 20th century for introducing me to the
beautiful and deep area even though I have never met him before.
Finally I would like to give my hearty thanks to Hiroyuki Yoshida
for explaining his important work kindly and sending two
references \cite{Yo3, Yo4} to me.

\vskip 0.1cm \noindent {\bf Notations:} \ \ We denote by
$\BQ,\,\BR$ and $\BC$ the field of rational numbers, the field of
real numbers and the field of complex numbers respectively. We
denote by $\BZ$ and $\BZ^+$ the ring of integers and the set of
all positive integers respectively. The symbol ``:='' means that
the expression on the right is the definition of that on the left.
For two positive integers $k$ and $l$, $F^{(k,l)}$ denotes the set
of all $k\times l$ matrices with entries in a commutative ring
$F$. For a square matrix $A\in F^{(k,k)}$ of degree $k$,
$\sigma(A)$ denotes the trace of $A$. For any $M\in F^{(k,l)},\
^t\!M$ denotes the transposed matrix of $M$. $I_n$ denotes the
identity matrix of degree $n$. For $A\in F^{(k,l)}$ and $B\in
F^{(k,k)}$, we set $B[A]=\,^tABA.$ For a complex matrix $A$,
${\overline A}$ denotes the complex {\it conjugate} of $A$. For
$A\in \BC^{(k,l)}$ and $B\in \BC^{(k,k)}$, we use the abbreviation
$B\{ A\}=\,^t{\overline A}BA.$ For a number field $F$, we denote
by ${\mathbb A}_F$ the ring of adeles of $F$. If $F=\BQ$, the
subscript will be omitted. We denote by ${\mathbb A}_{F,f}$ and
${\mathbb A}_f$ the finite part of ${\mathbb A}_F$ and ${\mathbb
A}$ respectively. By ${\overline\BQ}$ we mean the algebraic
closure of $\BQ$ in $\BC$.

\vskip 0.1cm

\end{section}
%
%
\begin{section}{{\bf Invariant Metrics and Laplacians on Siegel Space}}
\setcounter{equation}{0}
\newcommand\POB{ {{\partial}\over {\partial{\overline \Omega}}} }
\newcommand\PZB{ {{\partial}\over {\partial{\overline Z}}} }
\newcommand\PX{ {{\partial}\over{\partial X}} }
\newcommand\PY{ {{\partial}\over {\partial Y}} }
\newcommand\PU{ {{\partial}\over{\partial U}} }
\newcommand\PV{ {{\partial}\over{\partial V}} }
\newcommand\PO{ {{\partial}\over{\partial \Omega}} }
\newcommand\PZ{ {{\partial}\over{\partial Z}} }

\vskip 0.21cm For $\Om=(\omega_{ij})\in\BH_g,$ we write $\Om=X+iY$
with $X=(x_{ij}),\ Y=(y_{ij})$ real and $d\Om=(d\om_{ij})$. We
also put
$$\PO=\,\left(\,
{ {1+\delta_{ij}}\over 2}\, { {\partial}\over {\partial \om_{ij} }
} \,\right) \qquad\text{and}\qquad \POB=\,\left(\, {
{1+\delta_{ij}}\over 2}\, { {\partial}\over {\partial {\overline
{\om}}_{ij} } } \,\right).$$ C. L. Siegel \cite{Si1} introduced
the symplectic metric $ds^2$ on $\BH_g$ invariant under the action
(1.1) of $Sp(g,\BR)$ given by
\begin{equation}
ds^2=\s (Y^{-1}d\Om\, Y^{-1}d{\overline\Om})\end{equation} and H.
Maass \cite{M1} proved that its Laplacian is given by
\begin{equation}
\Delta=\,4\,\s \left( Y\,\,{
}^t\!\left(Y\POB\right)\PO\right).\end{equation} And
\begin{equation}
dv_g(\Om)=(\det Y)^{-(g+1)}\prod_{1\leq i\leq j\leq g}dx_{ij}\,
\prod_{1\leq i\leq j\leq g}dy_{ij}\end{equation} is a
$Sp(g,\BR)$-invariant volume element on
$\BH_g$\,(cf.\,\cite{Si3},\,p.\,130).

\vskip 0.2cm\noindent \begin{theorem}\,({\bf Siegel\,\cite{Si1}}).
(1) There exists exactly one geodesic joining two arbitrary points
$\Om_0,\,\Om_1$ in $\BH_g$. The length $\rho(\Om_0,\Om_1)$ of this
geodesic is given by

\begin{equation}
\rho(\Om_0,\Om_1)^2=\s \left( \left( \log {
{1+R(\Om_0,\Om_1)^{\frac 12} }\over {1-R(\Om_0,\Om_1)^{\frac 12} }
}\right)^2\right),\end{equation} where $R(\Om_1,\Om_2)$ is the
cross-ratio defined by
\begin{equation}
R(\Om_1,\Om_0)=(\Om_1-\Om_0)(\Om_1-{\overline
\Om}_0)^{-1}(\overline{\Om}_1-\overline{\Om}_0)(\overline{\Om}_1-\Om_0)^{-1}.
\end{equation}
(2) For $M\in Sp(g,\BR)$, we set
$${\tilde \Om}_0=M\cdot \Om_0\quad \textrm{and}\quad {\tilde \Om}_1=M\cdot
\Om_1.$$ Then $R(\Om_1,\Om_0)$ and
$R({\tilde\Om}_1,{\tilde\Om}_0)$ have the same eigenvalues.

\noindent (3) All geodesics are symplectic images of the special
geodesics
\begin{equation}
\alpha(t)=i\,\textrm{diag}(a_1^t,a_2^t,\cdots,a_g^t),
\end{equation}
where $a_1,a_2,\cdots,a_g$ are arbitrary positive real numbers
satisfying the condition
$$\sum_{k=1}^g \left( \log a_k\right)^2=1.$$
\end{theorem}
\noindent The proof of the above theorem can be found in
\cite{Si1}, pp.\,289-293.

\newcommand\OW{\overline{W}}
\newcommand\OP{\overline{P}}
\newcommand\OQ{\overline{Q}}
\newcommand\Dg{{\mathbb D}_g}
\newcommand\Hg{{\mathbb H}_g}

\vskip 0.1cm Let $$\BD_g=\left\{\,W\in\BC^{(g,g)}\,|\ W=\,{ }^tW,\
I_g-W{\overline W}>0\,\right\}$$ be the generalized unit disk of
degree $g$. The Cayley transform $\Psi:\Dg\lrt\Hg$ defined by
\begin{equation}
\Psi(W)=i\,(I_g+W)(I_g-W)^{-1},\quad W\in\Dg
\end{equation}
is a biholomorphic mapping of $\Dg$ onto $\Hg$ which gives the
bounded realization of $\Hg$ by $\Dg$\,(cf.\,\cite{Si1}). A.
Kor{\'a}nyi and J. Wolf \cite{KW} gave a realization of a bounded
symmetric domain as a Siegel domain of the third kind
investigating a generalized Cayley transform of a bounded
symmetric domain that generalizes the Cayley transform $\Psi$ of
$\BD_g$.

\vskip 0.2cm Let
\begin{equation}
T={1\over {\sqrt{2}} }\,
\begin{pmatrix} \ I_g&\ I_g\\
                   iI_g&-iI_g\end{pmatrix}
\end{equation}
be the $2g\times 2g$ matrix represented by $\Psi.$ Then
\begin{equation}
T^{-1}Sp(g,\BR)\,T=\left\{ \begin{pmatrix} P & Q \\ \OQ & \OP
\end{pmatrix}\,\Big|\ ^tP\OP-\,{}^t\OQ Q=I_g,\ {}^tP\OQ=\,{}^t\OQ
P\,\right\}.
\end{equation}
Indeed, if $M=\begin{pmatrix} A&B\\
C&D\end{pmatrix}\in Sp(g,\BR)$, then
\begin{equation}
T^{-1}MT=\begin{pmatrix} P & Q \\ \OQ & \OP
\end{pmatrix},
\end{equation}
where
\begin{equation}
P= {\frac 12}\,\Big\{ (A+D)+\,i\,(B-C)\Big\}
\end{equation}
and
\begin{equation}
 Q={\frac
12}\,\Big\{ (A-D)-\,i\,(B+C)\Big\}.
\end{equation}

For brevity, we set
\begin{equation*}
G_*= T^{-1}Sp(g,\BR)T.
\end{equation*}
Then $G_*$ is a subgroup of $SU(g,g),$ where
$$SU(g,g)=\left\{\,h\in\BC^{(g,g)}\,\big|\ {}^th I_{g,g}{\overline
h}=I_{g,g}\,\right\},\quad I_{g,g}=\begin{pmatrix} \ I_g&\ 0\\
0&-I_g\end{pmatrix}.$$ In the case $g=1$, we observe that
$$T^{-1}Sp(1,\BR)T=T^{-1}SL_2(\BR)T=SU(1,1).$$
If $g>1,$ then $G_*$ is a {\it proper} subgroup of $SU(g,g).$ In
fact, since ${}^tTJ_gT=-\,i\,J_g$, we get
\begin{equation}G_*=\Big\{\,h\in SU(g,g)\,\big|\
{}^thJ_gh=J_g\,\Big\}=SU(g,g)\cap Sp(g,\BC),
\end{equation}

\noindent where
$$Sp(g,\BC)=\Big\{\alpha\in \BC^{(2g,2g)}\ \big\vert \ ^t\!\alpha\, J_g\,\alpha= J_g\ \Big\}.$$

Let
\begin{equation*}
P^+=\left\{\begin{pmatrix} I_g & Z\\ 0 & I_g
\end{pmatrix}\,\Big|\ Z=\,{}^tZ\in\BC^{(g,g)}\,\right\}
\end{equation*}
be the $P^+$-part of the complexification of $G_*\subset SU(g,g).$
We note that the Harish-Chandra decomposition of an element
$\begin{pmatrix} P & Q\\ {\overline Q} & {\overline P}
\end{pmatrix}$ in $G_*$ is
\begin{equation*}
\begin{pmatrix} P & Q\\ \OQ & \OP
\end{pmatrix}=\begin{pmatrix} I_g & Q\OP^{-1}\\ 0 & I_g
\end{pmatrix} \begin{pmatrix} P-Q\OP^{-1}\OQ & 0\\ 0 & \OP
\end{pmatrix} \begin{pmatrix} I_g & 0\\ \OP^{-1}\OQ & I_g
\end{pmatrix}.
\end{equation*}
For more detail, we refer to \cite[p.\,155]{Kn}. Thus the
$P^+$-component of the following element
$$\begin{pmatrix} P & Q\\ \OQ & \OP
\end{pmatrix}   \cdot\begin{pmatrix} I_g & W\\ 0 & I_g
\end{pmatrix},\quad W\in \BD_g$$ of the complexification of $G_*^J$ is
given by
\begin{equation}
 \begin{pmatrix} I_g & (PW+Q)(\OQ W+\OP)^{-1}
\\ 0 & I_g
\end{pmatrix}.
\end{equation}
\newcommand\PW{ {{\partial}\over{\partial W}} }
\newcommand\PWB{ {{\partial}\over {\partial{\overline W}}} }
\newcommand\OVW{\overline W}

\noindent We note that $Q\OP^{-1}\in\Dg.$ We get the
Harish-Chandra embedding of $\Dg$ into $P^+$\,(cf.
\cite[p.\,155]{Kn} or \cite[pp.\,58-59]{Sa4}). Therefore we see
that $G_*$ acts on $\Dg$ transitively by
\begin{equation}
\begin{pmatrix} P & Q \\ \OQ & \OP
\end{pmatrix}\cdot W=(PW+Q)(\OQ W+\OP)^{-1},\quad \begin{pmatrix} P & Q \\ \OQ & \OP
\end{pmatrix}\in G_*,\ W\in \Dg.
\end{equation}
The isotropy subgroup $K_*$ of $G_*$ at the origin $o$ is given by
$$K_*=\left\{\,\begin{pmatrix} P & 0 \\ 0 & {\overline
P}\end{pmatrix}\,\Big|\ P\in U(g)\ \right\}.$$ Thus $G_*/K_*$ is
biholomorphic to $\Dg$. It is known that the action (1.1) is
compatible with the action (2.15) via the Cayley transform $\Psi$\
(cf.\,(2.7)). In other words, if $M\in Sp(g,\BR)$ and $W\in\BD_g$,
then
\begin{equation}
M\cdot \Psi(W)=\Psi(M_*\cdot W),
\end{equation}

\noindent where $M_*=T^{-1}MT\in G_*.$

\vskip 0.2cm For $W=(w_{ij})\in \Dg,$ we write $dW=(dw_{ij})$ and
$d{\overline W}=(d{\overline{w}}_{ij})$. We put $$\PW=\,\left(\, {
{1+\delta_{ij}}\over 2}\, { {\partial}\over {\partial w_{ij} } }
\,\right) \qquad\text{and}\qquad \PWB=\,\left(\, {
{1+\delta_{ij}}\over 2}\, { {\partial}\over {\partial {\overline
{w}}_{ij} } } \,\right).$$

Using the Cayley transform $\Psi:\Dg\lrt \BH_g$, Siegel showed
(cf.\,\,\cite{Si1}) that
\begin{equation}
ds_*^2=4 \s \Big((I_g-W{\overline W})^{-1}dW\,(I_g-\OVW
W)^{-1}d\OVW\,\Big)\end{equation} is a $G_*$-invariant Riemannian
metric on $\BD_g$ and Maass \cite{M1} showed that its Laplacian is
given by
\begin{equation}
\Delta_*=\,\s \left( (I_g-W\OW)\,{ }^t\!\left( (I_g-W\OVW)
\PWB\right)\PW\right).\end{equation}

\end{section}

%
%
\begin{section}{{\bf Invariant Differential Operators on Siegel Space}}
\setcounter{equation}{0}

\newcommand\POB{ {{\partial}\over {\partial{\overline \Omega}}} }
\newcommand\PZB{ {{\partial}\over {\partial{\overline Z}}} }
\newcommand\PO{ {{\partial}\over{\partial \Omega}} }
\newcommand\PZ{ {{\partial}\over{\partial Z}} }

\vskip 0.2cm For brevity, we write $G=Sp(g,\BR).$ The isotropy
subgroup $K$ at $iI_g$ for the action (1.1) is a maximal compact
subgroup given by
\begin{equation*}
K=\left\{ \begin{pmatrix} A & -B \\ B & A \end{pmatrix} \Big| \
A\,^t\!A+ B\,^t\!B=I_g,\ A\,^t\!B=B\,^t\!A,\ A,B\in
\BR^{(g,g)}\,\right\}.
\end{equation*}

\noindent Let $\fk$ be the Lie algebra of $K$. Then the Lie
algebra $\fg$ of $G$ has a Cartan decomposition $\fg=\fk\oplus
\fp$, where
\begin{equation*}
\fp=\left\{ \begin{pmatrix} X & Y \\ Y & -X \end{pmatrix} \Big| \
X=\,^tX,\ Y=\,^tY,\ X,Y\in \BR^{(g,g)}\,\right\}.
\end{equation*}

The subspace $\fp$ of $\fg$ may be regarded as the tangent space
of $\BH_g$ at $iI_g.$ The adjoint representation of $G$ on $\fg$
induces the action of $K$ on $\fp$ given by
\begin{equation}
k\cdot Z=\,kZ\,^tk,\quad k\in K,\ Z\in \fp.
\end{equation}

Let $T_g$ be the vector space of $g\times g$ symmetric complex
matrices. We let $\psi: \fp\lrt T_g$ be the map defined by
\begin{equation}
\psi\left( \begin{pmatrix} X & Y \\ Y & -X \end{pmatrix}
\right)=\,X\,+\,i\,Y, \quad \begin{pmatrix} X & Y \\ Y & -X
\end{pmatrix}\in \fp.
\end{equation}

\noindent We let $\delta:K\lrt U(g)$ be the isomorphism defined by
\begin{equation}
\delta\left( \begin{pmatrix} A & -B \\ B & A \end{pmatrix}
\right)=\,A\,+\,i\,B, \quad \begin{pmatrix} A & -B \\ B & A
\end{pmatrix}\in K,
\end{equation}

\noindent where $U(g)$ denotes the unitary group of degree $g$. We
identify $\fp$ (resp. $K$) with $T_g$ (resp. $U(g)$) through the
map $\Psi$ (resp. $\delta$). We consider the action of $U(g)$ on
$T_g$ defined by
\begin{equation}
h\cdot Z=\,hZ\,^th,\quad h\in U(g),\ Z\in T_g.
\end{equation}

\noindent Then the adjoint action (3.1) of $K$ on $\fp$ is
compatible with the action (3.4) of $U(g)$ on $T_g$ through the
map $\psi.$ Precisely for any $k\in K$ and $\omega\in \fp$, we get
\begin{equation}
\psi(k\,\omega \,^tk)=\delta(k)\,\psi(\omega)\,^t\delta (k).
\end{equation}

\noindent The action (3.4) induces the action of $U(g)$ on the
polynomial algebra $ \textrm{Pol}(T_g)$ and the symmetric algebra
$S(T_g)$ respectively. We denote by $ \textrm{Pol}(T_g)^{U(g)}$
$\Big( \textrm{resp.}\ S(T_g)^{U(g)}\,\Big)$ the subalgebra of $
\textrm{Pol}(T_g)$ $\Big( \textrm{resp.}\ S(T_g)\,\Big)$
consisting of $U(g)$-invariants. The following inner product $(\
,\ )$ on $T_g$ defined by $$(Z,W)= \, \textrm{tr}
\big(Z\,{\overline W}\,\big),\quad Z,W\in T_g$$

\noindent gives an isomorphism as vector spaces
\begin{equation}
T_g\cong T_g^*,\quad Z\mapsto f_Z,\quad Z\in T_g,
\end{equation}

\noindent where $T_g^*$ denotes the dual space of $T_g$ and $f_Z$
is the linear functional on $T_g$ defined by
$$f_Z(W)=(W,Z),\quad W\in T_g.$$

\noindent It is known that there is a canonical linear bijection
of $S(T_g)^{U(g)}$ onto the algebra ${\mathbb D}(\BH_g)$ of
differential operators on $\BH_g$ invariant under the action (1.1)
of $G$. Identifying $T_g$ with $T_g^*$ by the above isomorphism
(3.6), we get a canonical linear bijection
\begin{equation}
\Phi:\textrm{Pol}(T_g)^{U(g)} \lrt {\mathbb D}(\BH_g)
\end{equation}

\noindent of $ \textrm{Pol}(T_g)^{U(g)}$ onto ${\mathbb
D}(\BH_g)$. The map $\Phi$ is described explicitly as follows.
Similarly the action (3.1) induces the action of $K$ on the
polynomial algebra $ \textrm{Pol}(\fp)$ and $S(\fp)$ respectively.
Through the map $\psi$, the subalgebra $ \textrm{Pol}(\fp)^K$ of $
\textrm{Pol}(\fp)$ consisting of $K$-invariants is isomorphic to $
\textrm{Pol}(T_g)^{U(g)}$. We put $N=g(g+1)$. Let $\left\{
\xi_{\alpha}\,|\ 1\leq \alpha \leq N\, \right\}$ be a basis of
$\fp$. If $P\in \textrm{Pol}(\fp)^K$, then
\begin{equation}
\Big(\Phi (P)f\Big)(gK)=\left[ P\left( {{\partial}\over {\partial
t_{\al}}}\right)f\left(g\,\text{exp}\, \left(\sum_{\al=1}^N
t_{\al}\xi_{\al}\right) K\right)\right]_{(t_{\al})=0},
\end{equation} where $f\in C^{\infty}({\mathbb H}_{g})$. We refer to \cite{He1,He2} for more detail. In
general, it is hard to express $\Phi(P)$ explicitly for a
polynomial $P\in \textrm{Pol}(\fp)^K$.

\vskip 0.3cm According to the work of Harish-Chandra
\cite{HC1,HC2}, the algebra ${\mathbb D}(\BH_g)$ is generated by
$g$ algebraically independent generators and is isomorphic to the
commutative ring $\BC [x_1,\cdots,x_g]$ with $g$ indeterminates.
We note that $g$ is the real rank of $G$. Let $\fg_{\BC}$ be the
complexification of $\fg$. It is known that $\BD(\BH_g)$ is
isomorphic to the center of the universal enveloping algebra of
$\fg_{\BC}$\,(cf.\,\cite{Sh5}).

\vskip 0.3cm Using a classical invariant theory
(cf.\,\cite{Ho,Wey}), we can show that $\textrm{Pol}(T_g)^{U(g)}$
is generated by the following algebraically independent
polynomials
\begin{equation}
q_j (Z)=\,\textrm{tr}\Big( \big(Z {\overline
Z}\big)^j\,\Big),\quad j=1,2,\cdots,g.
\end{equation}

For each $j$ with $1\leq j\leq g,$ the image $\Phi(q_j)$ of $q_j$
is an invariant differential operator on $\BH_g$ of degree $2j$.
The algebra ${\mathbb D}(\BH_g)$ is generated by $g$ algebraically
independent generators $\Phi(q_1),\Phi(q_2),\cdots,\Phi(q_g).$ In
particular,
\begin{equation}
\Phi(q_1)=\,c_1\, \textrm{tr} \left( Y\,\,{
}^t\!\left(Y\POB\right)\PO\right)\quad  \textrm{for\ some
constant}\ c_1.
\end{equation}

\noindent We observe that if we take $Z=X+i\,Y$ with real $X,Y$,
then $q_1(Z)=q_1(X,Y)=\,\textrm{tr}\big( X^2 +Y^2\big)$ and
\begin{equation*}
q_2(Z)=q_2(X,Y)=\, \textrm{tr}\Big(
\big(X^2+Y^2\big)^2+\,2\,X\big(XY-YX\big)Y\,\Big).
\end{equation*}

\vskip 0.3cm We propose the following problem.

\vskip 0.2cm \noindent $ \textbf{Problem.}$ Express the images
$\Phi(q_j)$ explicitly for $j=2,3,\cdots,g.$

\vskip 0.3cm We hope that the images $\Phi(q_j)$ for
$j=2,3,\cdots,g$ are expressed in the form of the $\textit{trace}$
as $\Phi(q_1)$.

\vskip 0.3cm\noindent $ \textbf{Example 3.1.}$ We consider the
case $g=1.$ The algebra $ \textrm{Pol}(T_1)^{U(1)}$ is generated
by the polynomial
\begin{equation*}
q(z)=z\,{\overline z},\quad z\in \BC.
\end{equation*}

Using Formula (3.8), we get

\begin{equation*}
\Phi (q)=\,4\,y^2 \left( { {\partial^2}\over {\partial x^2} }+{
{\partial^2}\over {\partial y^2} }\,\right).
\end{equation*}

\noindent Therefore $\BD (\BH_1)=\BC\big[ \Phi(q)\big].$

\vskip 0.3cm\noindent $ \textbf{Example 3.2.}$ We consider the
case $g=2.$ The algebra $ \textrm{Pol}(T_2)^{U(2)}$ is generated
by the polynomial
\begin{equation*}
q_1(Z)=\,\s \big(Z\,{\overline Z}\,\big),\quad q_2(Z)=\,\s \Big(
\big(Z\,{\overline Z}\big)^2\Big), \quad Z\in T_2.
\end{equation*}

Using Formula (3.8), we may express $\Phi(q_1)$ and $\Phi(q_2)$
explicitly. $\Phi (q_1)$ is expressed by Formula (3.10). The
computation of $\Phi(q_2)$ might be quite tedious. We leave the
detail to the reader. In this case, $\Phi (q_2)$ was essentially
computed in \cite{BC}, Proposition 6. Therefore $\BD
(\BH_2)=\BC\big[ \Phi(q_1), \Phi(q_2)\big].$ The authors of
\cite{BC} computed the center of $U(\fg_{\BC}).$

\end{section}
%
%
\begin{section}{{\bf Siegel's Fundamental Domain}}
\setcounter{equation}{0} We let
$$\CP=\left\{\, Y\in\BR^{(g,g)}\,|\ Y=\,^tY>0\ \right\}$$
be an open cone in $\BR^N$ with $N=g(g+1)/2.$ The general linear
group $GL(g,\BR)$ acts on $\CP$ transitively by
\begin{equation}
g\circ Y:=gY\,^tg,\qquad g\in GL(g,\BR),\ Y\in \CP.\end{equation}
Thus $\CP$ is a symmetric space diffeomorphic to $GL(g,\BR)/O(g).$
\vskip 0.10cm
\newcommand\Mg{{\mathcal M}_g}
\newcommand\Rg{{\mathcal R}_g}
The fundamental domain $\Rg$ for $GL(g,\BZ)\ba \CP$ which was
found by H. Minkowski\,\cite{Min} is defined as a subset of $\CP$
consisting of $Y=(y_{ij})\in \CP$ satisfying the following
conditions (M.1)--(M.2)\ (cf. \cite{Ig} p.\,191 or \cite{M2}
p.\,123): \vskip 0.1cm (M.1)\ \ \ $aY\,^ta\geq y_{kk}$\ \ for
every $a=(a_i)\in\BZ^g$ in which $a_k,\cdots,a_g$ are relatively
prime for $k=1,2,\cdots,g$. \vskip 0.1cm (M.2)\ \ \ \
$y_{k,k+1}\geq 0$ \ for $k=1,\cdots,g-1.$ \vskip 0.1cm We say that
a point of $\Rg$ is {\it Minkowski reduced} or simply {\it M}-{\it
reduced}. $\Rg$ has the following properties (R1)--(R4): \vskip
0.1cm (R1) \ For any $Y\in\CP,$ there exist a matrix $A\in
GL(g,\BZ)$ and $R\in\Rg$ such that $Y=R[A]$\,(cf. \cite{Ig}
p.\,191 or \cite{M2} p.\,139). That is,
$$GL(g,\BZ)\circ \Rg=\CP.$$
\indent (R2)\ \ $\Rg$ is a convex cone through the origin bounded
by a finite number of hyperplanes. $\Rg$ is closed in $\CP$
(cf.\,\cite{M2} p.\,139).

\vskip 0.1cm (R3) If $Y$ and $Y[A]$ lie in $\Rg$ for $A\in
GL(g,\BZ)$ with $A\neq \pm I_g,$ then $Y$ lies on the boundary
$\partial \Rg$ of $\Rg$. Moreover $\Rg\cap (\Rg [A])\neq
\emptyset$ for only finitely many $A\in GL(g,\BZ)$ (cf.\,\cite{M2}
p.\,139). \vskip 0.1cm (R4) If $Y=(y_{ij})$ is an element of
$\Rg$, then
$$y_{11}\leq y_{22}\leq \cdots \leq y_{gg}\quad \text{and}\quad
|y_{ij}|<{\frac 12}y_{ii}\quad \text{for}\ 1\leq i< j\leq g.$$
\indent We refer to \cite{Ig} p.\,192 or \cite{M2} pp.\,123-124.
\vskip 0.1cm\noindent {\bf Remark.} Grenier\,\cite{Gr} found
another fundamental domain for $GL(g,\BZ)\ba \CP.$ \vskip 0.1cm
\newcommand\POB{ {{\partial}\over {\partial{\overline \Omega}}} }
\newcommand\PZB{ {{\partial}\over {\partial{\overline Z}}} }
\newcommand\PX{ {{\partial}\over{\partial X}} }
\newcommand\PY{ {{\partial}\over {\partial Y}} }
\newcommand\PU{ {{\partial}\over{\partial U}} }
\newcommand\PV{ {{\partial}\over{\partial V}} }
\newcommand\PO{ {{\partial}\over{\partial \Omega}} }
\newcommand\PZ{ {{\partial}\over{\partial Z}} }
\vskip 0.2cm For $Y=(y_{ij})\in \CP,$ we put
$$dY=(dy_{ij})\qquad\text{and}\qquad \PY\,=\,\left(\,
{ {1+\delta_{ij}}\over 2}\, { {\partial}\over {\partial y_{ij} } }
\,\right).$$ Then we can see easily that
\begin{equation}
ds^2=\s ( (Y^{-1}dY)^2)\end{equation} is a $GL(g,\BR)$-invariant
Riemannian metric on $\CP$ and its Laplacian is given by
$$\Delta=\s \left( \left( Y\PY\right)^2\right).$$
We also can see that
$$d\mu_g(Y)=(\det Y)^{-{ {g+1}\over2 } }\prod_{i\leq j}dy_{ij}$$
is a $GL(g,\BR)$-invariant volume element on $\CP$. The metric
$ds^2$ on $\CP$ induces the metric $ds_{\mathcal R}^2$ on $\Rg.$
Minkowski \cite{Min} calculated the volume of $\Rg$ for the volume
element $[dY]:=\prod_{i\leq j}dy_{ij}$ explicitly. Later Siegel
computed the volume of $\Rg$ for the volume element $[dY] $ by a
simple analytic method and generalized this case to the case of
any algebraic number field. \vskip 0.1cm Siegel\,\cite{Si1}
determined a fundamental domain $\CCF$ for $\G_g\ba \BH_g.$ We say
that $\Om=X+iY\in \BH_g$ with $X,\,Y$ real is {\it Siegel reduced}
or {\it S}-{\it reduced} if it has the following three properties:
\vskip 0.1cm (S.1)\ \ \ $\det (\text{Im}\,(\g\cdot\Om))\leq \det
(\text{Im}\,(\Om))\qquad\text{for\ all}\ \g\in\G_g$; \vskip 0.1cm
(S.2)\ \ $Y=\text{Im}\,\Om$ is M-reduced, that is, $Y\in \Rg\,;$
\vskip 0.1cm (S.3) \ \ $|x_{ij}|\leq {\frac 12}\quad \text{for}\
1\leq i,j\leq g,\ \text{where}\ X=(x_{ij}).$ \vskip 0.1cm $\CCF$
is defined as the set of all Siegel reduced points in $\BH_g.$
Using the highest point method, Siegel proved the following
(F1)--(F3)\,(cf. \cite{Ig} pp.\,194-197 or \cite{M2} p.\,169):
\vskip 0.1cm (F1)\ \ \ $\G_g\cdot \CCF=\BH_g,$ i.e.,
$\BH_g=\cup_{\g\in\G_g}\g\cdot \CCF.$ \vskip 0.1cm (F2)\ \ $\CCF$
is closed in $\BH_g.$ \vskip 0.1cm (F3)\ \ $\CCF$ is connected and
the boundary of $\CCF$ consists of a finite number of hyperplanes.
\vskip 0.21cm The metric $ds^2$ given by (2.1) induces a metric
$ds_{\mathcal F}^2$ on $\CCF.$ \vskip 0.1cm Siegel\,\cite{Si1}
computed the volume of $\CCF$
\begin{equation}
\text{vol}\,(\CCF)=2\prod_{k=1}^g\pi^{-k}\,\G
(k)\,\zeta(2k),\end{equation} where $\G (s)$ denotes the Gamma
function and $\zeta (s)$ denotes the Riemann zeta function. For
instance,
$$\text{vol}\,({\mathcal F}_1)={{\pi}\over 3},\quad \text{vol}\,({\mathcal F}_2)={{\pi^3}\over {270}},
\quad \text{vol}\,({\mathcal F}_3)={{\pi^6}\over {127575}},\quad
\text{vol}\,({\Cal F}_4)={{\pi^{10}}\over {200930625}}.$$

\vskip 0.3cm For a fixed element $\Om\in \BH_g,$ we set
\begin{equation*}
L_{\Om}:=\BZ^g+\BZ^g\Om,\qquad \BZ^g=\BZ^{(1,g)}.\end{equation*}
It follows from the positivity of $\text{Im}\,\Om$ that $L_{\Om}$
is a lattice in $\BC^g$. We see easily that if $\Om$ is an element
of $\BH_g$, the period matrix $\Om_*:=(I_g,\Om)$ satisfies the
Riemann conditions (RC.1) and (RC.2)\,: \vskip 0.1cm (RC.1) \ \ \
$\Om_*J_g,^t\Om_*=0$. \vskip 0.1cm (RC.2) \ \ \ $-{1 \over
{i}}\Om_*J_g\,^t{\overline{\Om}}_*
>0$.
\vskip 0.2cm \noindent Thus the complex torus
$A_{\Om}:=\BC^g/L_{\Omega}$ is an abelian variety.

\vskip 0.1cm We fix an element $\Om=X+iY$ of $\BH_g$ with
$X=\text{Re}\,\Om$ and $Y=\text{Im}\, \Om.$ For a pair $(A,B)$
with $A,B\in\BZ^g,$ we define the function $E_{\Om;A,B}:\BC^g\lrt
\BC$ by
\begin{equation*}
E_{\Om;A,B}(Z)=e^{2\pi i\left( \s(\,^tAU\,)+\, \s
((B-AX)Y^{-1}\,^tV)\right)},\end{equation*} where $Z=U+iV$ is a
variable in $\BC^g$ with real $U,V$. \vskip 0.1cm\noindent
\begin{lemma} For any $A,B\in \BZ^g,$ the function
$E_{\Om;A,B}$ satisfies the following functional equation
\begin{equation*}
E_{\Om;A,B}(Z+\la \Om+\mu)=E_{\Om;A,B}(Z),\quad
Z\in\BC^g\end{equation*} for all $\la,\mu\in\BZ^g.$ Thus
$E_{\Om;A,B}$ can be regarded as a function on $A_{\Om}.$ \vskip
0.1cm \end{lemma}
\begin{proof} The proof can be found in \cite{YJH8}.
\end{proof}
\newcommand\AO{A_{\Omega}}
\newcommand\Imm{\text{Im}}

 \vskip 0.1cm We
let $L^2(\AO)$ be the space of all functions $f:\AO\lrt\BC$ such
that
$$||f||_{\Om}:=\int_{\AO}|f(Z)|^2dv_{\Om},$$
where $dv_{\Om}$ is the volume element on $\AO$ normalized so that
$\int_{\AO}dv_{\Om}=1.$ The inner product $(\,\,,\,\,)_{\Om}$ on
the Hilbert space $L^2(\AO)$ is given by
\begin{equation*}
(f,g)_{\Om}:=\int_{\AO}f(Z)\,{\overline{g(Z)} }\,dv_{\Om},\quad
f,g\in L^2(\AO).\end{equation*}
\begin{theorem}
The set $\left\{\,E_{\Om;A,B}\,|\ A,B\in\BZ^g\,\right\}$ is a
complete orthonormal basis for $L^2(\AO)$. Moreover we have the
following spectral decomposition of $\Delta_{\Om}$:
$$L^2(\AO)=\bigoplus_{A,B\in \BZ^g}\BC\cdot E_{\Om;A,B}.$$
\end{theorem}

\noindent $\textrm{Proof.}$ The complete proof can be found in
\cite{YJH8}. \hfill $\square$

\end{section}

%
%
\begin{section}{\bf Siegel Modular Forms}
\setcounter{equation}{0}

\vskip 0.5cm \noindent {\bf 5.1. Basic Properties of Siegel
Modular Forms} \vskip 0.3cm Let $\rho$ be a rational
representation of $GL(g,\BC)$ on a finite dimensional complex
vector space $V_{\rho}$.

\vskip 0.3cm\noindent {\bf Definition.} \textrm{A holomorphic
function $f:\BH_g\lrt V_{\rho}$ is called a \textit{Siegel modular
form} with respect to $\rho$ if
\begin{equation}
f(\g\cdot \Omega)=f\big((A\Om+B)(C\Om+D)^{-1}\big)=\rho
(C\Om+D)f(\Om)
\end{equation}
for all $\begin{pmatrix} A & B \\ C & D \end{pmatrix}\in \G_g$ and
all $\Om\in\BH_g.$ Moreover if $g=1$, we require that $f$ is
holomorphic at the cusp $\infty$.}

We denote by $M_{\rho}(\G_g)$ the vector space of all Siegel
modular forms with respect to $\G_g$. If $\rho=\det^k$ for
$k\in\BZ$, a Siegel modular form $f$ with respect to $\rho$
satisfies the condition
\begin{equation}
f(\g\cdot \Omega)=\det (C\Om+D)^k\,f(\Om),
\end{equation}
where $\g$ and $\Om$ are as above. In this case $f$ is called a
(classical) Siegel modular form on $\BH_g$ of weight $k$. We
denote by $M_k(\G_g)$ the space of all Siegel modular forms on
$\BH_g$ of weight $k$.

\vskip 0.3cm\noindent {\bf Remark.} (1) If $\rho=\rho_1\oplus
\rho_2$ is a direct sum of two finite dimensional rational
representations of $GL(g,\BC)$, then it is easy to see that
$M_{\rho}(\G_g)$ is isomorphic to $M_{\rho_1}(\G_g)\oplus
M_{\rho_1}(\G_g).$ Therefore it suffices to study $M_{\rho}(\G_g)$
for an irreducible
representation $\rho$ of $GL(g,\BC)$.\\
\noindent (2) We may equip $V_{\rho}$ with a hermitian inner
product $(\ ,\ )$ satisfying the following condition
\begin{equation}
\big(\rho(x)v_1,v_2\big)=\big( v_1,
\overline{\rho(^tx)}v_2\big),\quad x\in GL(g, \BC),\ v_1,v_2\in
V_{\rho}.
\end{equation}

\vskip 0.3cm For an irreducible finite dimensional representation
$(\rho,V_{\rho})$ of $GL(g, \BC)$, there exist a highest weight
$k(\rho)=(k_1,\cdots,k_g)\in\BZ^g$ with $k_1\geq \cdots \geq k_g$
and a highest weight vector $v_{\rho}(\neq 0)\in V_{\rho}$ such
that
$$\rho\big(
\textrm{diag}(a_1,\cdots,a_g)\big)v_{\rho}=\prod_{i=1}^g
a_i^{k_i}\,v_{\rho},\quad a_1,\cdots,a_g\in \BC^{\times}.$$ Such a
vector $v_{\rho}$ is uniquely determined up to scalars. The number
$k(\rho):=k_g$ is called the $\textit{weight}$ of $\rho$. For
example, if $\rho=\det^k$, its highest weight is $(k,k,\cdots,k)$
and hence its weight is $k$.

\vskip 0.2cm Assume that $(\rho,V_{\rho})$ is an irreducible
finite dimensional rational representation of $GL(g,\BC)$. Then it
is known \cite{Ig,M2} that a Siegel modular form $f$ in
$M_{\rho}(\G_g)$ admits a Fourier expansion
\begin{equation}
f(\Om)=\sum_{T\geq 0} a(T)\,e^{2\pi i\, \s(T\Om)},
\end{equation}
where $T$ runs over the set of all half-integral semi-positive
symmetric matrices of degree $g$. We recall that $T$ is said to be
$ \textit{half-integral}$ if $2T$ is an integral matrix whose
diagonal entries are even.

\begin{theorem}
(1) If $kg$ is odd, then $M_k(\G_g)=0.$ \par\noindent (2) If
$k<0$, then $M_k(\G_g)=0.$
\par\noindent (3) Let $\rho$ be a non-trivial
irreducible finite dimensional representation of $GL(g,\BC)$ with
\par\noindent \ \ \ \ highest weight $(k_1,\cdots,k_g)$. If $M_{\rho}(\G_g)\neq \{ 0\},$
then $k_g\geq 1.$
\par\noindent (4) If $f\in M_{\rho}(\G_g)$, then $f$ is
bounded in any subset ${\Cal H}(c)$ of $\BH_g$ given by the form
\begin{equation*}
{\Cal H}(c):=\left\{ \Om\in \BH_g\,|\ \textrm{Im}\,\Om>
c\,I_g\,\right\} \end{equation*}
\par\noindent \ \ \ \ with any positive real number $c>0$.
\end{theorem}
\vskip 0.5cm \noindent {\bf 5.2. The Siegel Operator } \vskip
0.3cm Let $(\rho,V_{\rho})$ be an irreducible finite dimensional
representation of $GL(g,\BC)$. For any positive integer $r$ with
$0\leq r < g$, we define the operator $\Phi_{\rho,r}$ on
$M_{\rho}(\G_g)$ by
\begin{equation}
\big( \Phi_{\rho,r}f\big)(\Om_1):=\lim_{t\lrt\infty} f\left(
\begin{pmatrix} \Om_1 & 0 \\ 0 & it I_{g-r}
\end{pmatrix}\right),\quad f\in M_{\rho}(\G_g),\ \Om_1\in \BH_r.
\end{equation}

\noindent We see that $\Phi_{\rho,r}$ is well-defined because the
limit of the right hand side of (5.5) exists (cf.\,Theorem 5.1.
(4)). The operator $\Phi_{\rho,r}$ is called the $ \textit{Siegel
operator}$. A Siegel modular form $f\in M_{\rho}(\G_g)$ is said to
be a $ \textit{cusp form}$ if $\Phi_{\rho,g-1}f=0.$ We denote by
$S_{\rho}(\G_g)$ the vector space of all cusp forms on $\BH_g$
with respect to $\rho$. Let $V_{\rho}^{(r)}$ be the subspace of
$V_{\rho}$ spanned by the values
\begin{equation*}
\left\{ \big(\Phi_{\rho,r}f\big)(\Om_1)\,|\ \Om_1\in\BH_r,\ \,f\in
M_{\rho}(\G_g)\,\right\}.
\end{equation*}
According to \cite{W1}, $V_{\rho}^{(r)}$ is invariant under the
action of the subgroup
\begin{equation*}
\left\{ \begin{pmatrix} a & 0 \\ 0 & I_{g-r} \end{pmatrix} \,\Big|
\ a\in GL(r,\BC)\,\right\}.
\end{equation*}
\noindent Then we have an irreducible rational representation
$\rho^{(r)}$ of $GL(r,\BC)$ on $V_{\rho}^{(r)}$ defined by
\begin{equation*}
\rho^{(r)}(a)v:= \rho\left(
\begin{pmatrix} a & 0 \\ 0 & I_{g-r} \end{pmatrix}\right)v,\quad a
\in GL(r,\BC),\ v\in V_{\rho}^{(r)}.
\end{equation*}
\noindent We observe that if $(k_1,\cdots,k_g)$ is the highest
weight of $\rho$, then $(k_1,\cdots,k_r)$ is the highest weight of
$\rho^{(r)}$.

\begin{theorem}
The Siegel operator $\Phi_{\det^k,r}:M_k(\G_g)\lrt M_k(\G_r)$ is
surjective for $k$ even with $k> {{g+r+3}\over 2}$.
\end{theorem}

\vskip 0.2cm\noindent The proof of Theorem 5.2 can be found in
\cite{W2}.

\vskip 0.3cm We define the Petersson inner product $\langle\ ,\
\rangle_P$ on $M_{\rho}(\G_g)$ by
\begin{equation}
\langle f_1,f_2 \rangle_P:=\int_{{\Cal F}_g}\big(
\rho(\textrm{Im}\,\Om)f_1(\Om),f_2(\Om)\big)\,dv_g(\Om),\quad
f_1,f_2\in M_{\rho}(\G_g),
\end{equation}
where ${\mathcal F}_g$ is the Siegel's fundamental domain, $(\ ,\
)$ is the hermitian inner product defined in (5.3) and $dv_g(\Om)$
is the volume element defined by (2.3). We can check that the
integral of (5.6) converges absolutely if one of $f_1$ and $f_2$
is a cusp form. It is easily seen that one has the orthogonal
decomposition
\begin{equation*}
M_{\rho}(\G_g)=S_{\rho}(\G_g)\oplus S_{\rho}(\G_g)^{\perp},
\end{equation*}
\noindent where
\begin{equation*}
S_{\rho}(\G_g)^{\perp}= \big\{ f\in M_{\rho}(\G_g)\,\vert \
\langle f,h \rangle_P=0\ \textrm{for all}\ h\in
S_{\rho}(\G_g)\,\big\}
\end{equation*}
\noindent is the orthogonal complement of $S_{\rho}(\G_g)$ in
$M_{\rho}(\G_g)$.

\vskip 0.5cm \noindent {\bf 5.3. Construction of Siegel Modular
Forms} \vskip 0.3cm
In this subsection, we provide several
well-known methods to construct Siegel modular forms. \vskip 0.2cm
\noindent $\textsc{(A) Klingen's Eisenstein Series}$ \vskip 0.2cm
Let $r$ be an integer with $0\leq r<g.$ We assume that $k$ is a
positive $ \textit{even}$ integer. For $\Om\in \BH_g$, we write
\begin{equation*}
\Om=\begin{pmatrix} \Om_1 & * \\ * & \Om_2 \end{pmatrix},\quad
\Om_1\in \BH_r,\ \Om_2\in \BH_{g-r}.
\end{equation*}
For a fixed {\it cusp} form $f\in S_k(\G_r)$ of weight $k$, H.
Klingen \cite{K1} introduced the Eisenstein series $E_{g,r,k}(f)$
formally defined by
\begin{equation}
E_{g,r,k}(f)(\Om):=\sum_{\g\in P_r\backslash \G_g} f\big( (\g\cdot
\Om)_1\big)\cdot \det (C\Om+D)^{-k}, \ \quad \g= \begin{pmatrix} A
& B\\ C & D
\end{pmatrix}\in \G_g,
\end{equation}
where
\begin{equation*}
P_r=\left\{ \begin{pmatrix} A_1 & 0 & B_1 & * \\ * & U & * & * \\
C_1 & 0 & D_1 & * \\ 0 & 0 & 0 & ^tU^{-1} \end{pmatrix} \in \G_g
\ \Big|\ \begin{pmatrix} A_1  & B_1 \\
C_1  & D_1  \end{pmatrix}\in \G_r,\ U\in GL(g-r,\BZ)\ \right\}
\end{equation*}
\noindent is a parabolic subgroup of $\G_g.$ We note that if
$r=0$, and if $f=1$ is a constant, then
\begin{equation*}
E_{g,0,k}(\Om)=\sum_{C,D} \det (C\Omega+D)^{-k},
\end{equation*}
\noindent where $\begin{pmatrix}A & B\\ C& D \end{pmatrix}$ runs
over the set of all representatives for the cosets
$GL(g,\BZ)\backslash \G_g.$ \vskip 0.13cm Klingen \cite{K1} proved
the following\,:

\begin{theorem}
Let $g\geq 1$ and let $r$ be an integer with $0\leq r < g.$ We
assume that $k$ is a positive $\textit{even}$ integer with $k>
g+r+1.$ Then for any cusp form $f\in S_k(\G_r)$ of weight $k$,the
Eisenstein series $E_{g,r,k}(f)$ converges to a Siegel modular
form on $\BH_g$ of the same weight $k$ and one has the following
property
\begin{equation}
\Phi_{\det^k,\,r}E_{g,r,k}(f)=f.
\end{equation}
\end{theorem}
The proof of the above theorem can be found in \cite{K1,K2,M2}.

\vskip 0.2cm \noindent $\textsc{(B) Theta Series}$ \vskip 0.2cm
Let $(\rho,V_{\rho})$ be a finite dimensional rational
representation of $GL(g,\BC)$. We let $H_{\rho}(r,g)$ be the space
of pluriharmonic polynomials $P:\BC^{(r,g)}\lrt V_{\rho}$ with
respect to $(\rho,V_{\rho})$. That is, $P\in H_{\rho}(r,g)$ if and
only if $P:\BC^{(r,g)}\lrt V_{\rho}$ is a $V_{\rho}$-valued
polynomial on $\BC^{(r,g)}$ satisfying the following conditions
(5.9) and (5.10)\,: if $z=(z_{kj})$ is a coordinate in
$\BC^{(r,g)}$,
\begin{equation}
\sum_{k=1}^r { {\partial^2 P}\over {\partial z_{ki}\partial
z_{kj}} }=0\quad \textrm{for all}\ i,j\ \textrm{with}\ 1\leq
i,j\leq g
\end{equation}
\noindent and
\begin{equation}
P(zh)=\rho({}^th)\,\det (h)^{-{r\over 2}}P(z)\quad \textrm{for
all}\ z\in\BC^{(r,g)}\ \textrm{and}\ h\in GL(g,\BC).
\end{equation}
Now we let $S$ be a positive definite even unimodular matrix of
degree $r$. To a pair $(S,P)$ with $P\in H_{\rho}(r,g)$, we attach
the theta series
\begin{equation}
\Theta_{S,P}(\Om):=\sum_{A\in \BZ^{(r,g)}}P(S^{\frac 12}A)\,e^{\pi
i\,\s (S[A]\Om)}
\end{equation}
which converges for all $\Om\in \BH_g.$ E. Freitag \cite{Fr3}
proved that $\Theta_{S,P}$ is a Siegel modular form on $\BH_g$
with respect to $\rho$, i.e., $\Theta_{S,P}\in M_{\rho}(\G_g).$

\vskip 0.2cm Next we describe a method of constructing Siegel
modular forms using the so-called $ \textit{theta constants}$.
\vskip 0.2cm We consider a theta characteristic
\begin{equation*}
\epsilon=\begin{pmatrix} \epsilon' \\ \e'' \end{pmatrix}\in \{
0,1\}^{2g}\quad \textrm{with}\quad \e',\e''\in \{ 0,1\}^{g}.
\end{equation*}
A theta characteristic $\epsilon=\begin{pmatrix} \epsilon' \\ \e''
\end{pmatrix}$ is said to be $\textit{odd}$ (resp. $
\textit{even}$) if ${}^t\e'\e''$ is odd (resp. even). Now to each
theta characteristic $\epsilon=\begin{pmatrix} \epsilon' \\ \e''
\end{pmatrix}$, we attach the theta series
\begin{equation}
\theta[\e](\Om):=\sum_{m\in\BZ^g} e^{\pi i\,\big\{
\Om\big[m+{\frac 12}\e'\big]\,+\,{}^t\big(m+{\frac
12}\e'\big)\e''\,\big\} },\quad \Om\in \BH_g.
\end{equation}
If $\e$ is odd, we see that $\theta[\e]$ vanishes identically. If
$\e$ is even, $\theta[\e]$ is a Siegel modular form on $\BH_g$ of
weight ${\frac 12}$ with respect to the principal congreuence
subgroup $\G_g(2)$ (cf.\,\cite{Ig,Mf5}). Here
\begin{equation*}
\G_g(2)=\big\{\, \s\in \G_g\,|\ \s\equiv I_{2g}\ (\textrm{mod}\
2)\ \big\}
\end{equation*}
is a congruence subgroup of $\G_g$ of level $2$. These theta
series $\theta[\e]$ are called $ \textit{theta constants}$. It is
easily checked that there are $2^{g-1}(2^g+1)$ even theta
characteristics. These theta constants $\theta[\e]$ can be used to
construct Siegel modular forms with respect to $\G_g$. We provide
several examples. For $g=1$, we have
\begin{equation*}
\left( \theta[\e_{00}]\, \theta[\e_{01}]\,\theta[\e_{11}]
\right)^8 \in S_{12}(\G_1),
\end{equation*}
where
\begin{equation*}
\e_{00}=\begin{pmatrix} 0 \\ 0 \end{pmatrix},\quad
\e_{01}=\begin{pmatrix} 0
\\ 1 \end{pmatrix}\quad \textrm{and}\quad \e_{11}=\begin{pmatrix} 1
\\ 1 \end{pmatrix}.
\end{equation*}
For $g=2$, we get
\begin{equation*}
\chi_{10}:=-2^{-14} \prod_{\e\in {\mathbb E}}\theta[\e]^2 \in
S_{10}(\G_2)
\end{equation*}
and
\begin{equation*}
\left( \prod_{\e\in {\mathbb E}}\theta[\e]\right)\cdot
\sum_{\e_1,\e_2,\e_3}\left(
\theta[\e_1]\,\theta[\e_2]\,\theta[\e_3]\right)^{20} \in
S_{35}(\G_2),
\end{equation*}
where ${\mathbb E}$ denotes the set of all even theta
characteristics and $(\e_1,\e_2,\e_3)$ runs over the set of
triples of theta characteristics such that $\e_1+\e_2+\e_3$ is
odd. For $g=3$, we have
\begin{equation*}
 \prod_{\e\in {\mathbb E}}\theta[\e] \in
S_{18}(\G_3).
\end{equation*}
We refer to \cite{Ig} for more details.

\vskip 0.5cm \noindent {\bf 5.4. Singular Modular Forms} \vskip
0.3cm We know that a Siegel modular form $f\in M_{\rho}(\G_g)$ has
a Fourier expansion
\begin{equation*}
f(\Om)=\sum_{T\geq 0} a(T)\,e^{2\pi i\, \s(T\Om)},
\end{equation*}
where $T$ runs over the set of all half-integral semi-positive
symmetric matrices of degree $g$. A Siegel modular form $f\in
M_{\rho}(\G_g)$ is said to be $ \textit{singular}$ if $a(T)\neq 0$
implies $\det (T)=0$. We observe that the notion of singular
modular forms is opposite to that of cusp forms. Obviously if
$g=1$, singular modular forms are constants.\par We now
characterize singular modular forms in terms of the weight of
$\rho$ and a certain differential operator. For a coordinate
$\Om=X+iY$ in $\BH_g$ with $X$ real and $Y=(y_{ij})\in {\mathcal
P}_g$\,(cf. Section 4), we define the differential
\begin{equation}
M_g:=\det (Y)\cdot \det \left( {{\partial}\over{\partial Y}}
\right)
\end{equation}
which is invariant under the action (4.1) of $GL(g,\BR)$. Here
\begin{equation*}
{ {\partial}\over {\partial Y}}=\left( { {1+\delta_{ij}}\over 2 }
{ {\partial}\over {\partial y_{ij}} } \right).
\end{equation*}
Using the differential operator $M_g$, Maass
\cite[pp.\,202-204]{M2} proved that if a nonzero singular modular
form on $\BH_g$ of weight $k$ exists, then $nk\equiv 0$\,(mod $2$)
and $0<2k\leq g-1.$ The converse was proved by Weissauer
(cf.\,\cite[Satz 4]{W1}).

\begin{theorem} Let $\rho$ be an irreducible rational finite
dimensional representation of $GL(g,\BC)$ with highest weight
$(k_1,\cdots,k_g)$. Then a non-zero Siegel modular form $f\in
M_{\rho}(\G_g)$ is singular if and only if $2k(\rho)=2k_g<g.$
\end{theorem}
The above theorem was proved by Freitag \cite{Fr2}, Weissauer
\cite{W1} et al. By Theorem 5.6, we see that the weight of a
singular modular form is small. For instance, W. Duke and {\"O}.
Imamo$\check{g}$lu \cite{DI2} proved that $S_6(\G_g)=0$ for all
$g$. In a sense we say that there are $ \textit{no}$ cusp forms of
$\textit{small weight}$.

\begin{theorem} Let $f\in M_{\rho}(\G_g)$ be a Siegel modular form
with respect to a rational representation $\rho$ of $GL(g,\BC)$.
Then the following are equivalent\,:
\par (1) $f$ is a singular modular form.\par
(2) $f$ satisfies the differential equation $M_gf=0.$
\end{theorem}
We refer to \cite{M2} and \cite{YJH3} for the proof.

\vskip 0.2cm Let $f\in M_k(\G_g)$ be a nonzero singular modular
form of weight $k$. According to Theorem 5.4, $2k<g.$ We can show
that $k$ is divisible by $4$. Let $S_1,\cdots,S_h$ be a complete
system of representatives of positive definite even unimodular
integral matrices of degree $2k$. Freitag \cite{Fr2,Fr3} proved
that $f(\Om)$ can be written as a linear combination of theta
series $\theta_{S_1},\cdots,\theta_{S_h}$, where
$\theta_{S_{\nu}}\,(1\leq \nu\leq h)$ is defined by
\begin{equation}
\theta_{S_{\nu}}(\Om):=\sum_{A\in \BZ^{(2k,g)}}e^{\pi
i\,\s(S_{\nu} [A]\Om)},\quad 1\leq \nu\leq h.
\end{equation}

\vskip 0.2cm According to Theorem 5.5, we need to investigate some
properties of the weight of $\rho$ in order to understand singular
modular forms. Let $(k_1,\cdots,k_g)$ be the highest weight of
$\rho$. We define the $ \textit{corank}$ of $\rho$ by
\begin{equation*}
\textrm{corank}(\rho):=\Big| \Big\{ j\,|\ 1\leq j\leq g,\ k_j=k_g\
\Big\}\Big|.
\end{equation*}
Let
\begin{equation*}
f(\Om)=\sum_{T\geq 0} a(T)\,e^{2\pi i\, \s(T\Om)}
\end{equation*}
be a Siegel modular form in $M_{\rho}(\G_g)$. The notion of the $
\textit{rank}$ of $f$ and that of the $ \textit{corank}$ of $f$
were introduced by Weissauer \cite{W1} as follows\,:
\begin{equation*}
\textrm{rank}(f):= \textrm{max} \Big\{ \textrm{rank}\,(T)\ |\
a(T)\neq 0 \ \Big\}
\end{equation*}
and
\begin{equation*}
\textrm{corank}(f):= g-\textrm{min} \Big\{ \textrm{rank}\,(T)\ |\
a(T)\neq 0 \ \Big\}.
\end{equation*}
\indent Weissauer \cite{W1} proved the following.
\begin{theorem}
Let $\rho$ be an irreducible rational representation of
$GL(g,\BC)$ with highest weight $(k_1,\cdots,k_g)$ such that $
\textrm{corank}(\rho) < g-k_g.$ Assume that
\begin{equation*}
\Big| \Big\{ j\,|\ 1\leq j\leq g,\ k_j=k_g+1\ \Big\}\Big| <
2\big(g-k_g- \textrm{corank}(\rho)\big).
\end{equation*}
Then $M_{\rho}(\G_g)=0.$
\end{theorem}

\end{section}
\vskip 0.3cm
%
%
\begin{section}{\bf The Hecke Algebra}
\setcounter{equation}{0} \vskip 0.5cm \noindent {\bf 6.1. The
Structure of the Hecke Algebra} \vskip 0.3cm For a positive
integer $g$, we let $\G_g=Sp(g,\BZ)$ and let
\begin{equation*}
\Delta_g:=GSp(g,\BQ)=\big\{\,M\in GL(2g,\BQ)\,|\
{}^tMJ_gM=l(M)J_g,\ l(M)\in \BQ^{\times}\,\big\}
\end{equation*}
be the group of symplectic similitudes of the rational symplectic
vector space $(\BQ^{2g},\langle\ ,\ \rangle)$. We put
\begin{equation*}
\Delta_g^+:=GSp(g,\BQ)^+=\big\{\,M\in \Delta_g\,|\
 l(M)>0\,\big\}.
\end{equation*}

Following the notations in \cite{Fr3}, we let ${\mathscr
H}(\G_g,\Delta_g)$ be the complex vector space of all formal
finite sums of double cosets $\G_g M\G_g$ with $M\in \Delta_g^+$.
A double coset $\G_g M\G_g\,(M\in \Delta_g^+)$ can be written as a
finite disjoint union of right cosets $\G_gM_{\nu}\,(1\leq \nu\leq
h)\,:$
\begin{equation*}
\G_g M\G_g= \cup^h_{\nu=1} \G_gM_{\nu}\quad ( \textrm{disjoint}).
\end{equation*}
Let ${\mathscr L}(\G_g,\Delta_g)$ be the complex vector space
consisting of formal finite sums of right cosets $\G_gM$ with
$M\in \Delta^+$. For each double coset $\G_g M\G_g= \cup^h_{\nu=1}
\G_gM_{\nu}$ we associate an element $j(\G_gM\G_g)$ in ${\mathscr
L}(\G_g,\Delta_g)$ defined by
\begin{equation*}
j(\G_g M\G_g):= \sum^h_{\nu=1} \G_gM_{\nu}.
\end{equation*}
Then $j$ induces a linear map
\begin{equation}
j_*:{\mathscr H}(\G_g,\Delta_g)\lrt {\mathscr L}(\G_g,\Delta_g).
\end{equation}
We observe that $\Delta_g$ acts on ${\mathscr L}(\G_g,\Delta_g)$
as follows:
\begin{equation*}
\big( \sum_{j=1}^h c_j\,\G_gM_j\big)\cdot M=\sum_{j=1}^h c_j\,\G_g
M_jM,\quad M\in \Delta_g.
\end{equation*}
We denote
\begin{equation*}
{\mathscr L}(\G_g,\Delta_g)^{\G_g}:=\big\{\, T\in {\mathscr
L}(\G_g,\Delta_g)\,|\ T\cdot \g=T\ \textrm{for all}\ \g\in
\G_g\,\big\}
\end{equation*}
be the subspace of $\G_g$-invariants in ${\mathscr
L}(\G_g,\Delta_g)$. Then we can show that ${\mathscr
L}(\G_g,\Delta_g)^{\G_g}$ coincides with the image of $j_*$ and
the map
\begin{equation}
j_*:{\mathscr H}(\G_g,\Delta_g)\lrt {\mathscr
L}(\G_g,\Delta_g)^{\G_g}
\end{equation}
is an isomorphism of complex vector spaces
(cf.\,\cite[p.\,228]{Fr3}). From now on we identify ${\mathscr
H}(\G_g,\Delta_g)$ with ${\mathscr L}(\G_g,\Delta_g)^{\G_g}.$

We define the multiplication of the double coset $\G_g M\G_g$ and
$\G_g N$ by
\begin{equation}
(\G_g M\G_g)\cdot(\G_g N)=\sum_{j=1}^h \G_gM_jN,\quad M,N\in
\Delta_g,
\end{equation}
where $\G_g M\G_g= \cup^h_{j=1} \G_gM_{j}\ ( \textrm{disjoint}).$
The definition (6.3) is well defined, i.e., independent of the
choice of $M_j$ and $N$. We extend this multiplication to
${\mathscr H}(\G_g,\Delta_g)$ and ${\mathscr L}(\G_g,\Delta_g)$.
Since
\begin{equation*}
{\mathscr H}(\G_g,\Delta_g)\cdot {\mathscr
H}(\G_g,\Delta_g)\subset {\mathscr H}(\G_g,\Delta_g),
\end{equation*}
${\mathscr H}(\G_g,\Delta_g)$ is an associative algebra with the
identity element $\G_gI_{2g}\G_g=\G_g$. The algebra ${\mathscr
H}(\G_g,\Delta_g)$ is called the $ \textit{Hecke algebra}$ with
respect to $\G_g$ and $\Delta_g$. \vskip 0.2cm We now describe the
structure of the Hecke algebra ${\mathscr H}(\G_g,\Delta_g)$. For
a prime $p$, we let $\BZ[1/p]$ be the ring of all rational numbers
of the form $a\cdot p^{\nu}$ with $a,\nu\in\BZ.$ For a prime $p$,
we denote
\begin{equation*}
\Delta_{g,p}:=\Delta_g \cap GL\big(2g,\BZ[1/p]\big).
\end{equation*}
Then we have a decomposition of ${\mathscr H}(\G_g,\Delta_g)$
\begin{equation*}
{\mathscr H}(\G_g,\Delta_g)=\bigotimes_{p\,:\, \textrm{prime}}
{\mathscr H}(\G_g,\Delta_{g,p})
\end{equation*}
as a tensor product of local Hecke algebras ${\mathscr
H}(\G_g,\Delta_{g,p}).$ We denote by $\check{\mathscr
H}(\G_g,\Delta_g)$\,(resp. $\check{\mathscr H}(\G_g,\Delta_{g,p})$
the subring of ${\mathscr H}(\G_g,\Delta_g)$\,(resp. ${\mathscr
H}(\G_g,\Delta_{g,p})$ by integral matrices.

In order to describe the structure of local Hecke operators
${\mathscr H}(\G_g,\Delta_{g,p})$, we need the following lemmas.

\begin{lemma} Let $M\in \Delta_g^+$ with ${}^tMJ_gM=lJ_g$. Then the
double coset $\G_gM\G_g$ has a unique representative of the form
\begin{equation*}
M_0= \textrm{diag}(a_1,\cdots,a_g,d_1,\cdots,d_g),
\end{equation*}
\noindent where $a_g|d_g,\ a_j>0,\ a_jd_j=l$ for $1\leq j\leq g$
and $a_k| a_{k+1}$ for $1\leq k\leq g-1.$
\end{lemma}

For a positive integer $l$, we let
\begin{equation*}
O_g(l):=\big\{\, M\in GL(2g,\BZ)\,|\ {}^tMJ_gM=lJ_g\ \big\}.
\end{equation*}
Then we see that $O_g(l)$ can be written as a finite disjoint
union of double cosets and hence as a finite union of right
cosets. We define $T(l)$ as the element of ${\mathscr
H}(\G_g,\Delta_g)$ defined by $O_g(l).$

\begin{lemma} (a) Let $l$ be a positive integer. Let
\begin{equation*}
O_g(l)=\cup_{\nu=1}^h\G_gM_{\nu}\quad ( \textrm{disjoint})
\end{equation*}
be a disjoint union of right cosets $\G_gM_{\nu}\,(1\leq \nu\leq
h).$ Then each right coset $\G_gM_{\nu}$ has a representative of
the form
\begin{equation*}
M_{\nu}=\begin{pmatrix} A_{\nu}&B_{\nu}\\
                   0&   D_{\nu}\end{pmatrix}, \quad
                   {}^tA_{\nu}D_{\nu}=lI_g,\quad A_{\nu}\  \textrm{is upper
                   triangular}.
\end{equation*}
(b) Let $p$ be a prime. Then
\begin{equation*}
T(p)=O_g(p)=\G_g\begin{pmatrix} I_g&0\\
                   0&  pI_g\end{pmatrix}\G_g
\end{equation*}
and
\begin{equation*}
T(p^2)=\sum_{i=0}^g T_i(p^2),
\end{equation*}
where
\begin{equation*}
T_k(p^2):=\begin{pmatrix} I_{g-k}&0 & 0 & 0\\
                   0&  pI_k  & 0 & 0 \\ 0 & 0 & p^2I_{g-k} & 0\\
                   0 & 0 & 0 & pI_k\end{pmatrix}\G_g,\quad 0\leq
                   k\leq g.
\end{equation*}
\end{lemma}
\noindent $ \textit{Proof.}$ The proof can be found in
\cite[p.\,225 and p.\,250]{Fr3}. \hfill$\square$

For example, $T_g(p^2)=\G_g(pI_{2g})\G_g$ and
\begin{equation*}
T_0(p^2)=\G_g\begin{pmatrix} I_g&0\\
                   0&  p^2I_g\end{pmatrix}\G_g=T(p)^2.
\end{equation*}

We have the following
\begin{theorem} The local Hecke algebra $\check{\mathscr H}(\G_g,\Delta_{g,p})$
is generated by algebraically independent generators
$T(p),\,T_1(p^2),\cdots,T_g(p^2).$
\end{theorem}
\noindent $ \textit{Proof.}$ The proof can be found in
\cite[p.\,250 and p.\,261]{Fr3}. \hfill$\square$

\vskip 0.2cm On $\D_g$ we have the anti-automorphism $M\mapsto
M^*:=l(M)M^{-1}\,(M\in \D_g)$. Obviously $\G_g^*=\G_g$. By Lemma
6.1, $(\G_gM\G_g)^*=\G_gM^*\G_g=\G_gM\G_g.$ According to
\cite{Sh1}, Proposition 3.8, ${\mathscr{H}}(\G_g,\Delta_{g})$ is
commutative. \vskip 0.2cm Let $X_0,X_1,\cdots,X_g$ be the $g+1$
variables. We define the automorphisms
\begin{equation*}
w_j: \BC\big[X_0^{\pm 1},X_1^{\pm 1},\cdots,X_g^{\pm 1}\big]\lrt
\BC\big[X_0^{\pm 1},X_1^{\pm 1},\cdots,X_g^{\pm 1}\big],\quad
1\leq j\leq g
\end{equation*}
by
\begin{equation*}
w_j(X_0)=X_0X_j^{-1},\ \ \ w_j(X_j)=X_j^{-1},\ \ \ w_j(X_k)=X_k\ \
\textrm{for} \ k\neq 0,j.
\end{equation*}
Let $W_g$ be the finite group generated by $w_1,\cdots,w_g$ and
the permutations of variables $X_1,\cdots,X_g$. Obviously $w_j^2$
is the identity map and $|W_g|=2^g g!$.
\begin{theorem} There exists an isomorphism
\begin{equation*}
Q:{\mathscr H}(\G_g,\Delta_{g,p})\lrt \BC\big[X_0^{\pm 1
},X_1^{\pm 1 },\cdots,X_g^{\pm 1}\big]^{W_g}.
\end{equation*}
In fact, $Q$ is defined by
\begin{equation*}
Q\big( \sum_{j=1}^h\G_gM_j\big)=\sum_{j=1}^h
Q(\G_gM_j)=\sum_{j=1}^h
X_0^{-k_0(j)}\prod_{\nu=1}^g\big(p^{-\nu}X^{\nu}\big)^{k_{\nu}(j)}|\det
A_j|^{g+1},
\end{equation*}
where we choose the representative $M_j$ of $\G_g M_j$ of the form
\begin{equation*}
M_j=\begin{pmatrix} A_j& B_j\\
                   0&  D_j\end{pmatrix},\quad A_j=\begin{pmatrix} p^{k_1(j)}& \ldots & *\\
                   0 & \ddots & \vdots \\
                   0&  0 & p^{k_g(j)}\end{pmatrix}.
\end{equation*}
We note that the integers $k_1(j),\cdots,k_g(j)$ are uniquely
determined.

\end{theorem}
\noindent $ \textit{Proof.}$ The proof can be found in \cite{Fr3}.
\hfill $\square$

\vskip 0.3cm For a prime $p$, we let
\begin{equation*}
{\mathscr H}(\G_g,\Delta_{g,p})_{\BQ}:=\left\{\, \sum c_j\,\G_g
M_j\G_g \in {\mathscr H}(\G_g,\Delta_{g,p})\,|\ c_j\in\BQ\
\right\}
\end{equation*}
be the $\BQ$-algebra contained in ${\mathscr
H}(\G_g,\Delta_{g,p})$. We put
\begin{equation*}
G_p:=GSp(g,\BQ_p)\quad \textrm{and}\quad K_p=GSp(g,\BZ_p).
\end{equation*}
We can identify ${\mathscr H}(\G_g,\Delta_{g,p})_{\BQ}$ with the
$\BQ$-algebra ${\mathscr H}_{g,p}^{\BQ}$ of $\BQ$-valued locally
constant, $K_p$-biinvariant functions on $G_p$ with compact
support. The multiplication on ${\mathscr H}_{g,p}^{\BQ}$ is given
by
\begin{equation*}
(f_1*f_2)(h)=\int_{G_p}f_1(g)\,f_2(g^{-1}h)dg,\quad f_1,f_2\in
{\mathscr H}_{g,p}^{\BQ},
\end{equation*}
where $dg$ is the unique Haar measure on $G_p$ such that the
volume of $K_p$ is $1$. The correspondence is obtained by sending
the double coset $\G_gM\G_g$ to the characteristic function of
$K_pMK_p$. \vskip 0.2cm In order to describe the structure of
${\mathscr H}_{g,p}^{\BQ}$, we need to understand the $p$-adic
Hecke algebras of the diagonal torus ${\mathbb T}$ and the Levi
subgroup ${\mathbb M}$ of the standard parabolic group. Indeed,
${\mathbb T}$ is defined to be the subgroup consisting of diagonal
matrices in $\D_g$ and
\begin{equation*}
{\mathbb M}=\left\{\,\begin{pmatrix} A& 0\\
                   0&  D\end{pmatrix}\in \D_g\ \right\}
\end{equation*}
is the Levi subgroup of the parabolic subgroup
\begin{equation*}
\left\{\,\begin{pmatrix} A& B\\
                   0&  D\end{pmatrix}\in \D_g\ \right\}.
\end{equation*}
Let $Y$ be the co-character group of ${\mathbb T}$, i.e., $Y=
\textrm{Hom}({\mathbb G}_m,{\mathbb T}).$ We define the local
Hecke algebra ${\mathscr H}_p({\mathbb T})$ for ${\mathbb T}$ to
be the $\BQ$-algebra of $\BQ$-valued, ${\mathbb
T}(\BZ_p)$-biinvariant functions on $\BT(\BQ_p)$ with compact
support. Then ${\mathscr H}_p(\BT)\cong \BQ[Y],$ where $\BQ[Y]$ is
the group algebra over $\BQ$ of $Y$. An element $\la\in Y$
corresponds the characteristic function of the double coset
$D_{\la}=K_p\la(p)K_p$. It is known that  ${\mathscr H}_p({\mathbb
T})$ is isomorphic to the ring $\BQ\big[(u_1/v_1)^{\pm
1},\cdots,(u_g/v_g)^{\pm 1},(v_1\cdots v_g)^{\pm 1}\big]$ under
the map
\begin{equation*}
(a_1,\cdots,a_g,c)\mapsto (u_1/v_1)^{a_1}\cdots
(u_g/v_g)^{a_g}(v_1\cdots v_g)^{c}.
\end{equation*}
Similarly we have a $p$-adic Hecke algebra ${\mathscr
H}_p({\mathbb M})$. Let $W_{\D_g}=N(\BT)/\BT$ be the Weyl group
with respect to $(\BT,\D_g)$, where $N(\BT)$ is the normalizer of
$\BT$ in $\D_g$. Then $W_{\D_g}\cong S_g \ltimes (\BZ/2\BZ)^g,$
where the generator of the $i$-th factor $\BZ/2\BZ$ acts on a
matrix of the form $ \textrm{diag}(a_1,\cdots,a_g,d_1,\cdots,d_g)$
by interchanging $a_i$ and $d_i$, and the symmetry group $S_g$
acts by permuting the $a_i${'}$s$ and $d_i${'}$s$. We note that
$W_{\D_g}$ is isomorphic to $W_g$. The Weyl group $W_{\BM}$ with
respect to $(\BT,\BM)$ is isomorphic to $S_g$. We can prove that
the algebra ${\mathscr H}_p({\mathbb T})^{W_{\D_g}}$ of
$W_{\D_g}$-invariants in ${\mathscr H}_p({\mathbb T})$ is
isomorphic to $\BQ \big[ Y_0^{\pm 1},Y_1,\cdots,
Y_g\big]$\,(cf.\,\cite{Fr3}). We let
\begin{equation*}
B=\left\{\,\begin{pmatrix} A& B\\
                   0&  D\end{pmatrix}\in \D_g\,\Big|\ A\
                   \textrm{is upper triangular,}\ D\ \textrm{is lower triangular}\ \right\}
\end{equation*}
be the Borel subgroup of $\D_g$. A set $\Phi^+$ of positive roots
in the root system $\Phi$ determined by $B$. We set $\rho={\frac
12}\sum_{\alpha \in \Phi^+}\alpha.$

\vskip 0.2cm Now we have the map $\alpha_{\BM}:\BM\lrt {\mathbb
G}_m$ defined by
\begin{equation*}
\alpha_{\BM}(M):=l(M)^{-{{g(g+1)}\over 2}} \big( \det
A\big)^{g+1},\quad M=\begin{pmatrix}A & 0\\ 0 & D\end{pmatrix}\in
\BM
\end{equation*}
and the map $\beta_{\BT}:\BT\lrt {\mathbb G}_m$ defined by
\begin{equation*}
\beta_{\BT}(\textrm{diag}(a_1,\cdots,a_g,d_1,\cdots,d_g)):=\prod_{i=1}^g
a_1^{g+1-2i},\quad \textrm{diag}(a_1,\cdots,a_g,d_1,\cdots,d_g)\in
\BT.
\end{equation*}
Let $\theta_{\BT}:=\alpha_{\BM}\,\beta_{\BT}$ be the character of
$\BT.$ The $ \textit{Satake's spherical map}\ S_{p,\BM}:{\mathscr
H}_{g,p}^{\BQ}\lrt {\mathscr H}_p(\BM)$ is defined by
\begin{equation}
S_{p,\BM}(\phi)(m):=|\alpha_\BM(m)|_p
\int_{U(\BQ_p)}\phi(mu)du,\quad \phi\in {\mathscr H}_{g,p}^{\BQ},\
m\in \BM,
\end{equation}
where $|\ \ |_p$ is the $p$-adic norm and $U(\BQ_p)$ denotes the
unipotent radical of $\D_g$. Also another $\textit{Satake's
spherical map}\ S_{\BM,\BT}:{\mathscr H}_p({\BM})\lrt {\mathscr
H}_p(\BT)$ is defined by
\begin{equation}
S_{\BM,\BT}(f)(t):=|\beta_\BT(t)|_p \int_{\BM\cap{\mathbb
N}}f(tn)dn,\quad t\in {\mathscr H}_p(\BT),\ t\in \BT,
\end{equation}
where ${\mathbb N}$ is a nilpotent subgroup of $\D_g$.

\begin{theorem} The Satake's spherical maps $S_{p,\BM}$ and
$S_{\BM,\BT}$ define the isomorphisms of $\BQ$-algebras
\begin{equation}
{\mathscr H}_{g,p}^{\BQ}\cong {\mathscr H}_p(\BT)^{W_{\D_g}} \quad
\textrm{and}\quad {\mathscr H}_p(\BM)\cong {\mathscr
H}_p(\BT)^{W_{\BM}}.
\end{equation}
\end{theorem}
\noindent We define the elements $\phi_k\,(0\leq k\leq g)$ in
${\mathscr H}_p(\BM)$ by
\begin{equation*}
\phi_k:=p^{-{{k(k+1)}\over 2}}\,\BM(\BZ_p)\begin{pmatrix} I_{g-k}
& 0 & 0 \\ 0 & pI_g & 0 \\ 0 & 0 & I_k \end{pmatrix}
\BM(\BZ_p),\quad i=0,1,\cdots,g.
\end{equation*}
Then we have the relation
\begin{equation}
S_{p,\BM}(T(p))=\sum_{k=0}^g\phi_k
\end{equation}
and
\begin{equation}
S_{p,\BM}\big(T_i(p^2)\big)=\sum_{j,k\geq 0,\,i+j\leq k}
m_{k-j}(i)\, p^{- {k-j+1\choose 2}}\phi_j\phi_k,
\end{equation}
where
\begin{equation*}
m_s(i):=\sharp \left\{\,A\in M(s,{\mathbb F}_p)\,|\ {}^tA=A,\quad
\textrm{corank}(A)=i\ \right\}.
\end{equation*}
Moreover, for $k=0,1,\cdots,g$, we have
\begin{equation}
S_{\BM,\BT}(\phi_k)=(v_1\cdots v_g)E_k(u_1/v_1,\cdots,u_g/v_g),
\end{equation}
where $E_k$ denotes the elementary symmetric function of degree
$k$. The proof of (6.7)-(6.9) can be found in
\cite[pp.\,142-145]{A2}.

\vskip 0.5cm \noindent {\bf 6.2. Action of the Hecke Algebra on
Siegel Modular Forms}

\vskip 0.3cm Let $(\rho,V_{\rho})$ be a finite dimensional
irreducible representation of $GL(g,\BC)$ with highest weight
$(k_1,\cdots,k_g)$. For a function $F:\BH_g\lrt V_{\rho}$ and
$M\in \D_g^+,$ we define
\begin{equation*}
(f|_{\rho}M)(\Om)=\rho(C\Omega+D)^{-1}f(M\cdot\Om),\quad
M=\begin{pmatrix}A & B\\ C & D \end{pmatrix}\in \D_g^+.
\end{equation*}
It is easily checked that $f|_{\rho}M_1M_2=\big(
f|_{\rho}M_1\big)|_{\rho}M_2$ for $M_1,M_2\in \D_g^+.$

\vskip 0.2cm We now consider a subset ${\mathscr M}$ of $\D_g$
satisfying the following properties (M1) and (M2)\,:
\par \ \ (M1)\ \ \ ${\mathscr M}=\cup_{j=1}^h \G_gM_j\quad$ (disjoint
union);\par \ \ (M2)\ \ \ ${\mathscr M}\,\G_g\subset {\mathscr
M}.$

\vskip 0.2cm For a Siegel modular form $f\in M_{\rho}(\G_g)$, we
define
\begin{equation}
T({\mathscr M})f:=\sum_{j=1}^h f|_{\rho}M_j.
\end{equation}
This is well defined, i.e., is independent of the choice of
representatives $M_j$ because of the condition (M1). On the other
hand, it follows from the condition (M2) that $T({\mathscr
M})f|_{\rho}\g=T({\mathscr M})f$ for all $\g\in \G_g.$ Thus we get
a linear operator
\begin{equation}
T({\mathscr M}):M_{\rho}(\G_g)\lrt M_{\rho}(\G_g).
\end{equation}
We know that each double coset $\G_gM\G_g$ with $M\in \D_g$
satisfies the condition $(M1)$ and $(M2)$. Thus a linear operator
$T({\mathscr M})$ defined in (6.10 induces naturally the action of
the Hecke algebra ${\mathscr H}(\G_g,\Delta_g)$ on
$M_{\rho}(\G_g)$. More precisely, if ${\mathscr N}=\sum_{j=1}^h
c_j\G_gM_j\G_g\in {\mathscr H}(\G_g,\Delta_g)$, we define
\begin{equation*}
T({\mathscr N})=\sum_{j=1}^h c_jT(\G_gM_j\G_g).
\end{equation*}
Then $T({\mathscr N})$ is an endomorphism of $M_{\rho}(\G_g)$.

\vskip 0.2cm Now we fix a Siegel modular form $F$ in
$M_{\rho}(\G_g)$ which is an eigenform of the Hecke algebra
${\mathscr H}(\G_g,\Delta_g)$. Then we obtain an algebra
homomorphism $\la_F:{\mathscr H}(\G_g,\Delta_g)\lrt \BC$
determined by
\begin{equation*}
T(F)=\la_F(T)F,\quad T\in {\mathscr H}(\G_g,\Delta_g).
\end{equation*}
By Theorem 6.2 or Theorem 6.3, one has
\begin{eqnarray*}
{\mathscr H}(\G_g,\Delta_{g,p})&\cong& {\mathscr
H}_{g,p}^{\BQ}\otimes \BC \cong \BC[Y]^{W_g}\\
&\cong& {\mathscr H}_p({\mathbb T})^{W_g}\otimes \BC\\
&\cong& \BC\big[(u_1/v_1)^{\pm 1},\cdots,(u_g/v_g)^{\pm
1},(v_1\cdots v_g)^{\pm 1}\big]^{W_g}\\
&\cong& \BC[Y_0,Y_0^{-1},Y_1,\cdots,Y_g],
\end{eqnarray*}
where $Y_0,Y_1,\cdots,Y_g$ are algebraically independent.
Therefore one obtains an isomorphism
\begin{equation*}
\textrm{Hom}_{\BC}\big({\mathscr
H}(\G_g,\Delta_{g,p}),\BC\big)\cong
\textrm{Hom}_{\BC}\big({\mathscr H}_{g,p}^{\BQ}\otimes
\BC,\BC\big)\cong (\BC^{\times})^{(g+1)}/W_g.
\end{equation*}
The algebra homomorphism $\la_F\in
\textrm{Hom}_{\BC}\big({\mathscr H}(\G_g,\Delta_{g,p}),\BC\big)$
is determined by the $W_g$-orbit of a certain $(g+1)$-tuple $\big(
\alpha_{F,0},\alpha_{F,1},\cdots,\alpha_{F,g}\big)$ of nonzero
complex numbers, called the $p$-$ \textit{Satake parameters}$ of
$F$. For brevity, we put $\alpha_i=\alpha_{F,i},\ i=0,1,\cdots,g$.
Therefore $\alpha_i$ is the image of $u_i/v_i$ and $\alpha_0$ is
the image of $v_1\cdots v_g$ under the map $\Theta$. Each
generator $w_i\in W_{\D_g}\cong W_g$ acts by
\begin{equation*}
w_j(\alpha_0)=\alpha_0\alpha_j^{-1}\quad
w_j(\alpha_j)=\alpha_j^{-1},\quad w_j(\alpha_k)=0\ \textrm{if}\
k\neq 0,j.
\end{equation*}
These $p$-Satake parameters $\alpha_0,\alpha_1\cdots,\alpha_g$
satisfy the relation
\begin{equation*}
\alpha_0^2\alpha_1\cdots \alpha_g=p^{\sum_{i=1}^gk_i-g(g+1)/2}.
\end{equation*}
Formula (6.12) follows from the fact that
$T_g(p^2)=\G_g(pI_{2g})\G_g$ is mapped to
$$p^{-g(g+1)/2}\,(v_1\cdots
v_g)^2\prod_{i=1}^g (u_i/v_i).$$ We refer to \cite[p.\,258]{Fr3}
for more detail. According to Formula (6.7)-(6.9), the eigenvalues
$\la_F\big(T(p)\big)$ and $\la_F\big(T_i(p^2)\big)$ with $1\leq
i\leq g$ are given respectively by
\begin{equation}
\la_F\big(T(p)\big)=\alpha_0(1+E_1+E_2+\cdots+E_g)
\end{equation}
and
\begin{equation}
\la_F\big(T_i(p^2)\big)=\sum_{j,k\geq 0,\,j+i\leq k}^g
m_{k-j}(i)\,p^{-{ {k-j+1}\choose 2}}\,\alpha_0^2E_jE_k,\quad
i=1,\cdots,g,
\end{equation}
where $E_j$ denotes the elementary symmetric function of degree
$j$ in the variables $\alpha_1,\cdots,\alpha_g$. The point is that
the above eigenvalues $\la_F\big(T(p)\big)$ and
$\la_F\big(T_i(p^2)\big)\ (1\leq i\leq g)$ are described in terms
of the $p$-Satake parameters $\alpha_0,\alpha_1\cdots,\alpha_g$.
\vskip 0.2cm \noindent $ \textbf{Examples.}$ (1) Suppose
$g(\tau)=\sum_{n\geq 1}a(n)\,e^{2\pi i n\tau}$ is a normalized
eigenform in $S_k(\G_1)$. Let $p$ be a prime. Let $\beta$ be a
complex number determined by the relation
\begin{equation*}
(1-\beta X)(1-{\bar \beta}X)=1-a(p)X+p^{k-1}X^2.
\end{equation*}
Then
\begin{equation*}
\beta+{\bar\beta}=a(p)\quad \textrm{and}\quad \beta {\bar
\beta}=p^{k-1}.
\end{equation*}
The $p$-Satake parameters $\alpha_0$ and $\alpha_1$ are given by
\begin{equation*}
(\alpha_0,\alpha_1)=\left( \beta, {{\bar\beta}\over
\beta}\right)\quad or \quad \left( {\bar\beta}, {{\beta}\over
{\bar\beta}}\right).
\end{equation*}
It is easily checked that
$\alpha_0^2\alpha_1=\beta{\bar\beta}=p^{k-1}$ (cf. Formula
(6.12)). \vskip 0.2cm\noindent (b) For a positive integer $k$ with
$k>g+1$, we let
\begin{equation*}
G_k(\Om):=\sum_{M\in \G_{g,0}\backslash \G_g} \det
(C\Om+D)^k,\quad M=\begin{pmatrix} A & B\\ C & D\end{pmatrix}
\end{equation*}
be the Siegel Eisenstein series of weight $k$ in $M_k(\G_g)$,
where
$$\G_{g,0}:=\left\{ \begin{pmatrix} A & B\\ 0 &
D\end{pmatrix}\in \G_g\right\}$$ is a parabolic subgroup of
$\G_g$. It is known that $G_k$ is an eigenform of all the Hecke
operators\,(cf.\,\cite[p.\,268]{Fr3}). Let $S_1,\cdots,S_h$ be a
complete system of representatives of positive definite even
unimodular integral matrices of degree $2k$. If $k>g+1$, the
Eisenstein series $G_k$ can be expressed as the weighted mean of
theta series $\theta_{S_1},\cdots,\theta_{S_h}\,$:
\begin{equation}
G_k(\Om)=\sum_{\nu=1}^h m_{\nu}\,\theta_{S_{\nu}}(\Om),\quad
\Om\in \BH_g,
\end{equation}
where
\begin{equation*}
m_{\nu}={ {A(S_{\nu},S_{\nu})^{-1}}\over
{A(S_1,S_1)^{-1}+\cdots+A(S_h,S_h)^{-1}} },\quad 1\leq \nu\leq h.
\end{equation*}
We recall that the theta series $\theta_{S_{\nu}}$ is defined in
Formula (5.14) and that for two symmetric integral matrices $S$ of
degree $m$ and $T$ of degree $n$, $A(S,T)$ is defined by
\begin{equation*}
A(S,T):=\sharp \big\{\, G\in \BZ^{(m,n)}\,|\
S[G]=\,{}^tGSG=T\,\big\}.
\end{equation*}
Formula (6.14) was obtained by Witt \cite{Wit} as a special case
of the analytic version of Siegel's Hauptsatz.

\end{section}

%
%
\begin{section}{{\bf Jacobi Forms}}
\setcounter{equation}{0} \vskip 0.3cm In this section, we
establish the notations and define the concept of Jacobi forms.
\vspace{0.05in}\\
\indent Let $$Sp(g,\BR)=\lb M\in \BR^{(2g,2g)}\ \vert \ ^t\!MJ_gM=
J_g\ \rb$$ be the symplectic group of degree $g$, where
$$J_g:=\begin{pmatrix} 0 & I_g \\
                   -I_g & 0
                   \end{pmatrix}.$$

For two positive integers $g$ and $h$, we consider the {\it
Heisenberg group}
\begin{equation*}
H_\BR^{(g,h)}:=\big\{ (\la,\mu,\kappa)\,|\ \la,\mu\in
\BR^{(h,g)},\ \kappa\in\BR^{(h,h)},\ \kappa+\mu\,{}^t\la\
\textrm{symmetric}\,\big\}
\end{equation*}
endowed with the following multiplication law
\begin{equation*}
 (\la,\mu,\kappa)\circ (\la',\mu',\kappa):=(\la+\la',\mu+\mu',
 \kappa+\kappa'+\la\,{}^t\mu'-\mu\,{}^t\la').
\end{equation*}
\indent We recall that the Jacobi group
$G^J_{g,h}:=Sp(g,\BR)\ltimes H_{\BR}^{(g,h)}$ is the semidirect
product of the symplectic group $Sp(g,\BR)$ and the Heisenberg
group $H_{\BR}^{(g,h)}$ endowed with the following multiplication
law
$$
(M,(\lambda,\mu,\kappa))\cdot(M',(\lambda',\mu',\kappa')) :=\,
(MM',(\tilde{\lambda}+\lambda',\tilde{\mu}+ \mu',
\kappa+\kappa'+\tilde{\lambda}\,^t\!\mu'
-\tilde{\mu}\,^t\!\lambda'))$$ with $M,M'\in Sp(g,\BR),
(\lambda,\mu,\kappa),\,(\lambda',\mu',\kappa') \in
H_{\BR}^{(g,h)}$ and
$(\tilde{\lambda},\tilde{\mu}):=(\lambda,\mu)M'$. It is easy to
see that $G_{g,h}^J$ acts on the Siegel-Jacobi space
$\BH_{g,h}:=\BH_g\times \BC^{(h,g)}$ transitively by
\begin{equation}
(M,(\lambda,\mu,\kappa))\cdot (\Om,Z):=(M\cdot\Om,(Z+\lambda
\Om+\mu) (C\Om+D)^{-1}),
\end{equation}
where $M=\begin{pmatrix} A&B\\ C&D\end{pmatrix} \in Sp(g,\BR),\
(\lambda,\mu, \kappa)\in H_{\BR}^{(g,h)}$ and $(\Om,Z)\in
\BH_{g,h}.$
\vspace{0.05in}\\
\indent Let $\rho$ be a rational representation of $GL(g,\BC)$ on
a finite dimensional complex vector space $V_{\rho}.$ Let
${\mathcal M}\in \BR^{(h,h)}$ be a symmetric half-integral
semi-positive definite matrix of degree $h$. Let
$C^{\infty}(\BH_{g,h},V_{\rho})$ be the algebra of all
$C^{\infty}$ functions on $\BH_{g,h}$ with values in $V_{\rho}.$
For $f\in C^{\infty}(\BH_{g,h}, V_{\rho}),$ we define
\begin{eqnarray*}
  & \big(f|_{\rho,{\mathcal{M}}}[(M,(\lambda,\mu,\kappa))]\big)(\Om,Z)\hskip 7cm\\
&\quad\quad:=  e^{-2\pi i\sigma({\mathcal{M}}[Z+\lambda
\Om+\mu](C\Om+D)^{-1}C)} \times e^{2\pi
i\sigma({\mathcal{M}}(\lambda
\Om^t\!\lambda+2\lambda^t\!Z+(\kappa+
\mu^t\!\lambda)))} \hskip 5cm\\
&\times\rho(C\Om+D)^{-1}f(M\cdot\Om,(Z+\lambda
\Om+\mu)(C\Om+D)^{-1}), \hskip 3.8cm
\end{eqnarray*}
where $M=\begin{pmatrix} A&B\\ C&D\end{pmatrix}\in Sp(g,\BR),\
(\lambda,\mu,\kappa)\in H_{\BR}^{(g,h)}$ and $(\Om,Z)\in
\BH_{g,h}.$
\vspace{0.1in}\\
\def\l{\lambda}
\noindent{\bf Definition\ 7.1.}\quad Let $\rho$ and $\M$ be as
above. Let
$$H_{\BZ}^{(g,h)}:=\lb \, (\l,\mu,\kappa)\in H_{\BR}^{(g,h)}\, \vert
\, \l,\mu\in \BZ^{(h,g)},\ \kappa\in \BZ^{(h,h)}\, \rb.$$ Let $\G$
be a discrete subgroup of $\G_g$ of finite index. A {\it Jacobi\
form} of index $\M$ with respect to $\rho$ on $\G$ is a
holomorphic function $f\in C^{\infty}(\BH_{g,h},V_{\rho})$
satisfying the following conditions (A) and (B):
\vspace{0.1in}\\
\indent (A) $f|_{\rho,{\mathcal{M}}}[\tilde{\gamma}] = f$ for all
$\tilde{\gamma}\in
\Gamma^J := \G \ltimes H_{\BZ}^{(g,h)}$.\vspace{0.1in}\\
\indent (B) $f$ has a Fourier expansion of the following form :
$$f(\Om,Z) = \sum\limits_{T\ge0\atop \text {half-integral}}
\sum\limits_{R\in \BZ^{(g,h)}} c(T,R)\cdot e^{{{2\pi i}\over
{\l_{\G}}}\,\sigma(T\Om)}\cdot e^{2\pi i\sigma(RZ)}$$ with some
nonzero integer $\la_{\G}\in \BZ$ and $c(T,R)\ne 0$ only if
$\begin{pmatrix}
{1\over {\la_{\G}}}T & \frac 12R\\
\frac 12\,^t\!R&{\mathcal{M}}\end{pmatrix}
\ge 0$.\\
\indent If $g\geq 2,$ the condition (B) is superfluous by the K{\"
o}cher principle\,(\,cf.\,\cite{Zi} Lemma 1.6). We denote by
$J_{\rho,\M}(\G)$ the vector space of all Jacobi forms of index
$\M$ with respect to $\rho$ on $\G$. Ziegler(cf. \cite{Zi} Theorem
1.8 or \cite{EZ} Theorem 1.1) proves that the vector space
$J_{\rho,\M}(\G)$ is finite dimensional. For more results on
Jacobi forms with $g>1$ and $h>1$, we refer to \cite{Ru},
\cite{YJH0}-\cite{YJH4} and \cite{Zi}.
\vspace{0.1in}\\
\noindent{\bf Definition 7.2.}\quad A Jacobi form $f\in
J_{\rho,\M}(\G)$ is said to
be a {\it cusp} (or {\it cuspidal}) form if $\begin{pmatrix} {1\over {\la_{\G}}}T & \frac 12R\\
\frac 12\,^t\!R & {\mathcal{M}}\end{pmatrix} > 0$ for any $T,\ R$
with $c(T,R)\ne 0.$ A Jacobi form $f\in J_{\rho,\M}(\Gamma)$ is
said to be {\it singular} if it admits a Fourier expansion such
that a Fourier coefficient $c(T,R)$ vanishes unless
$\text{det}\,\begin{pmatrix} {1\over {\la_{\G}}}T &\frac 12R\\
\frac 12\,^t\!R& {\mathcal{M}} \end{pmatrix}=0.$
\vspace{0.1in}\\
\noindent{\bf Example 7.3.}\quad Let $S\in \BZ^{(2k,2k)}$ be a
symmetric, positive definite, unimodular even integral matrix and
$c\in \BZ^{(2k,h)}.$ We define the theta series
\begin{equation}
\vartheta_{S,c}^{(g)}(\Om,Z):=\sum_{\l\in \BZ^{(2k,g)}} e^{\pi
i\{\s(S\la \Om\,^t\!\l)+2\s(\,^t\!cS\l\,^t\!Z)\} },\ \ \Om\in
\BH_g, \ Z\in \BC^{(h,g)}.
\end{equation}
We put $\M:=\frac 12 ^t\!cSc.$ We assume that $2k<g+rank\,(\M).$
Then it is easy to see that $\vartheta_{S,c}^{(g)}$ is a singular
Jacobi form in $J_{k,\M}(\G_g)$(cf. \cite{Zi} p.212).
\vspace{0.1in}\\
\noindent{\bf Remark.} Singular Jacobi forms are characterized by
a special differential operator or the weight of the associated
rational representation of the general linear group
$GL(g,\BC)$\,(cf.\,\cite{YJH3}).

\vspace{0.1in} Now we will make  brief historical remarks on
Jacobi forms. In 1985, the names Jacobi group and Jacobi forms got
kind of standard by the classic book \cite{EZ} by Eichler and
Zagier to remind of Jacobi's ``Fundamenta nova theoriae functionum
ellipticorum'', which appeared in 1829 (cf.\,\cite{J}). Before
\cite{EZ} these objects appeared more or less explicitly and under
different names in the work of many authors. In 1966
Pyatetski-Shapiro \cite{PS} discussed the Fourier-Jacobi expansion
of Siegel modular forms and the field of modular abelian
functions. He gave the dimension of this field in the higher
degree. About the same time Satake \cite{Sa3}-\cite{Sa4}
introduced the notion of ``groups of Harish-Chandra type'' which
are non reductive but still behave well enough so that he could
determine their canonical automorphic factors and kernel
functions. Shimura \cite{Sh3}-\cite{Sh4} gave a new foundation of
the theory of complex multiplication of abelian functions using
Jacobi theta functions. Kuznetsov \cite{Kuz} constructed functions
which are almost Jacobi forms from ordinary elliptic modular
functions. Starting 1981, Berndt \cite{Be1}-\cite{Be3} published
some papers which studied the field of arithmetic Jacobi
functions, ending up with a proof of Shimura reciprocity law for
the field of these functions with arbitrary level. Furthermore he
investigated the discrete series for the Jacobi group $G^J_{g,h}$
and developed the spectral theory for $L^2(\G^J\backslash
G^J_{g,h}$) in the case $g=h=1$\,(cf.\,\cite{Be4}-\cite{Be6}). The
connection of Jacobi forms to modular forms was given by Maass,
Andrianov, Kohnen, Shimura, Eichler and Zagier. This connection is
pictured as follows. For $k$ even, we have the following
isomorphisms
$$M_k^*(\G_2)\,\cong\,J_{k,1}(\G_1)\,\cong\,M_{k-\frac12}^+(\G_0^{(1)}(4))\,
\cong\,M_{2k-2}(\G_1).$$ Here $M_k^*(\G_2)$ denotes Maass's
Spezialschar or Maass space and $M_{k-\frac12}^+(\G_0^{(1)}(4))$
denotes the Kohnen plus space. These spaces shall be described in
some more detail in the next section. For a precise detail, we
refer to \cite{M3}-\cite{M5},\,\cite{A1},\,\cite{EZ} and
\cite{Ko1}. In 1982 Tai \cite{Ta} gave asymptotic dimension
formulae for certain spaces of Jacobi forms for arbitrary $g$ and
$h=1$ and used these ones to show that the moduli ${\mathcal A}_g$
of principally polarized abelian varieties of dimension $g$ is
{\it of\ general\ type} for $g\geq 9.$ Feingold and Frenkel
\cite{FF} essentially discussed Jacobi forms in the context of
Kac-Moody Lie algebras generalizing the Maass correspondence to
higher level. Gritsenko \cite{Gri} studied Fourier-Jacobi
expansions and a non-commutative Hecke ring in connection with the
Jacobi group. After 1985 the theory of Jacobi forms for $g=h=1$
had been studied more or less systematically by the Zagier school.
A large part of the theory of Jacobi forms of higher degree was
investigated by Kramer \cite{Kr1}-\cite{Kr2}, \cite{Ru}, Yang
\cite{YJH0}-\cite{YJH4}and Ziegler \cite{Zi}. There were several
attempts to establish $L$-functions in the context of the Jacobi
group by Murase \cite{Mu1}-\cite{Mu2}and Sugano \cite{MS} using
the so-called ``Whittaker-Shintani functions''. Kramer
\cite{Kr1}-\cite{Kr2} developed an arithmetic theory of Jacobi
forms of higher degree. Runge \cite{Ru} discussed some part of the
geometry of Jacobi forms for arbitrary $g$ and $h=1.$ For a good
survey on some motivation and background for the study of Jacobi
forms, we refer to \cite{Be7}. The theory of Jacobi forms has been
extensively studied by many people until now and has many
applications in other areas like geometry and physics.

\end{section}

%
%
\begin{section}{{\bf Lifting of Elliptic Cusp forms to Siegel Modular Forms}}
\setcounter{equation}{0} \vskip 0.2cm In this section, we presents
some results about the liftings of elliptic cusp forms to Siegel
modular forms. And we discuss the Duke-Imamo${\check g}$lu-Ikeda
lift.
\par In order to discuss these lifts, we need two kinds of
$L$-function or zeta functions associated to Siegel Hecke
eigenforms. These zeta functions are defined by using the Satake
parameters of their associated Siegel Hecke eigenforms. \par Let
$F\in M_{\rho}(\G_g)$ be a nonzero Hecke eigenform on $\BH_g$ of
type $\rho$, where $\rho$ is a finite dimensional irreducible
representation of $GL(g,\BC)$ with highest weight
$(k_1,\cdots,k_g)$. Let
$\alpha_{p,0},\alpha_{p,1},\cdots,\alpha_{p,g}$ be the $p$-Satake
parameters of $F$ at a prime $p$. Using these Satake parameters,
we define the $ \textit{local spinor zeta function}\ Z_{F,p}(s)$
of $F$ at $p$ by
\begin{equation*}
Z_{F,p}(t)=(1-\alpha_{p,0}t)\prod_{r=1}^g\prod_{1\leq i_1< \cdots<
i_r\leq g}(1-\alpha_{p, 0}\alpha_{p,i_1}\cdots\alpha_{p,i_r}t).
\end{equation*}
Now we define the $\textit{ spinor zeta function}\ Z_F(s)$ by
\begin{equation}
Z_F(s):=\prod_{p\,:\, \textrm{prime}} Z_{F,p}(p^{-s})^{-1},\quad
\textrm{Re}\,s\gg 0.
\end{equation}
For example, if $g=1$, the spinor zeta function $Z_f(s)$ of a
Hecke eigenform $f$ is nothing but the Hecke $L$-function $L(f,s)$
of $f$.

\vskip 0.2cm Secondly one has the so-called $ \textit{standard
zeta function}\ D_F(s)$ of a Hecke eigenform $F$ in
$S_{\rho}(\G_g)$ defined by
\begin{equation}
D_F(s):=\prod_{p\,:\, \textrm{prime}} D_{F,p}(p^{-s})^{-1},
\end{equation}
where
\begin{equation*}
D_{F,p}(t)=(1-t)\prod_{i=1}^g
(1-\alpha_{p,i}t)(1-\alpha_{p,i}^{-1}t).
\end{equation*}
For instance, if $g=1$, the standard zeta function $D_f(s)$ of a
Hecke eigenform $f(\tau)=\sum_{n=1}^{\infty} a(n)e^{2\pi i n\tau}$
in $S_k(\G_1)$ has the following
\begin{equation*}
D_f(s-k+1)=\prod_{p\,:\, \textrm{prime}}\big(
1+p^{-s+k-1}\big)^{-1}\cdot\sum_{n=1}^{\infty}a(n^2)n^{-s}.
\end{equation*}
For the present time being, we recall the Kohnen plus space and
the Maass space. Let ${\mathcal M}$ be a positive definite,
half-integral symmetric matrix of degree $h$. For a fixed element
$\Om\in\BH_g$, we denote $\Theta_{{\mathcal M},\Om}^{(g)}$ the
vector space consisting of all the functions
$\theta:\BC^{(h,g)}\lrt \BC$ satisfying the condtition\,:
\begin{equation}
\theta(Z+\la\Om+\mu)=e^{-2\pi i\,\s({\mathcal
M}[\la]\Om+2\,{}^tZ{\mathcal M}\la)},\quad Z\in\BC^{(h,g)}
\end{equation}
for all $\la,\mu\in\BZ^{(h,g)}.$ For brevity, we put
$L:=\BZ^{(h,g)}$ and $L_{\mathcal M}:=L/(2{\mathcal M})L.$ For
each $\g\in{L}_{\mathcal M}$, we define the theta series
$\theta_{\g}(\Om,Z)$ by
\begin{equation*}
\theta_{\g}(\Om,Z)=\sum_{\la\in L} e^{2\pi i\,\s({\mathcal
M}[\la+(2{\mathcal M})^{-1}\g]\Om+2\,{}^t\!Z {\mathcal
M}(\la+(2{\mathcal M})^{-1}\g))},
\end{equation*}
where $(\Om,Z)\in\BH_g\times\BC^{(h,g)}.$ Then
$\left\{\,\theta_{\g}(\Om,Z)\,|\ \g\in L_{\mathcal M}\,\right\}$
forms a basis for $\Theta_{{\mathcal M},\Om}^{(g)}$. For any
Jacobi form $\phi(\Om,Z)\in J_{k,{\mathcal M}}(\G_g),$ the
function $\phi(\Om,\cdot)$ with fixed $\Om$ is an element of
$\Theta_{{\mathcal M},\Om}^{(g)}$ and $\phi(\Om,Z)$ can be
written as a linear combination of theta series
$\theta_{\g}(\Om,Z)\,(\,\g\in L_{\mathcal M})\,:$
\begin{equation}
\phi(\Om,Z)=\sum_{\g\in L_{\mathcal M}}\phi_{\g}(\Om)
\theta_{\g}(\Om,Z).
\end{equation}
We observe that $\phi=(\phi_{\g}(\Om))_{\g\in L_{\mathcal M}}$ is
a vector valued automorphic form with respect to a theta
multiplier system.

\vskip 0.2cm We now consider the case\,: $h=1,\ {\mathcal
M}=I_h=1,\ L=\BZ^{(1,g)}\cong \BZ^g.$ We define the theta series
$\theta^{(g)}(\Om)$ by
\begin{equation}
\theta^{(g)}(\Om)=\sum_{\la\in L}e^{2\pi
i\,\s(\la\Om\,{}^t\!\la)}=\theta_0(\Om,0),\quad \Om\in\BH_g.
\end{equation}
Let
\begin{equation*}
\G_0^{(g)}(4)=\left\{\,\begin{pmatrix} A & B\\ C & D
\end{pmatrix}\in \G_g\,\Big|\ C\equiv 0\,\,(
\textrm{mod}\,4)\,\right\}
\end{equation*}
be the congruence subgroup of $\G_g.$ We define the automorphic
factor $j:\G_0^{(g)}(4)\times \BH_g\lrt \BC^{\times}$ by
\begin{equation*}
j(\g,\Om)={ {\theta^{(g)}(\g\cdot\Om)}\over {\theta^{(g)}(\Om)}
},\quad \g\in \G_0^{(g)}(4),\ \Om\in \BH_g.
\end{equation*}
Thus one obtains the relation
\begin{equation*}
j(\g,\Om)^2=\varepsilon (\g)\,\det (C\Om+D),\quad
\varepsilon(\g)^2=1,
\end{equation*}
for any $\g=\begin{pmatrix} A & B\\ C & D
\end{pmatrix}\in\G_0^{(g)}(4).$

\vskip 0.2cm Kohnen \cite{Ko1} introduced the so-called $
\textit{Kohnen plus space}$ $M_{k-{\frac
12}}^+\big(\G_0^{(g)}(4)\big)$ consisting of holomorphic functions
satisfying the following conditions $(K1)$ and $(K2)$\,:\par (K1)
\ \ \ $f(\g\cdot\Om)=j(\g,\Om)^{2k-1}f(\Om)\quad $ for all $\g\in
\G_0^{(g)}(4)\,;$\par (K2) $f$ has the Fourier expansion
\begin{equation*}
f(\Om)=\sum_{T\geq 0}a(T)\,e^{2\pi i\,\s(T\Om)},
\end{equation*}
where $T$ runs over the set of semi-positive half-integral
symmetric matrices of degree $g$ such that $a(T)=0$ unless
$T\equiv -\mu\,{}^t\!\mu$\,mod\,$4S^*_g(\BZ)$ for some $\mu\in
\BZ^{(g,1)}$. Here we put
\begin{equation*}
S_g^*(\BZ)=\left\{\,T\in\BR^{(g,g)}\,|\ T=\,{}^tT,\ \s(TS)\in\BZ\
\ \textrm{for all}\ S=\,{}^tS\in\BZ^{(g,g)}\,\right\}.
\end{equation*}
\indent For a Jacobi form $\phi\in J_{k,1}(\G_g),$ according to
Formula (8.4), one has
\begin{equation}
\phi(\Om,Z)=\sum_{\g\in L/2L}f_{\g}(\Om)\,
\theta_{\g}(\Om,Z),\quad \Om\in \BH_g,\ Z\in \BC^{(h,g)} .
\end{equation}
Now we put
\begin{equation*}
f_{\phi}(\Om):=\sum_{\g\in L/2L}f_{\g}(4\Om),\quad \Om\in \BH_g.
\end{equation*}
Then $f_{\phi}\in M_{k-{\frac 12}}^+\big(\G_0^{(g)}(4)\big).$

\begin{theorem} $ \textbf{(Kohnen-Zagier (g=1),\ Ibukiyama ($g >1$))}$ Suppose
$k$ is an even positive integer. Then there exists the isomorphism
given by
\begin{equation*}
J_{k,1}(\G_g)\cong M_{k-{\frac
12}}^+\big(\G_0^{(g)}(4)\big),\quad\ \phi\mapsto f_{\phi}.
\end{equation*}
Moreover the isomorphism is compatible with the action of Hecke
operators.
\end{theorem}
For a positive integer $k\in\BZ^+$, H. Maass \cite{M3,M4,M5}
introduced the so-called $ \textit{Maass space}$
$M_k^*(\G_2)$ consisting of all Siegel modular forms
$F(\Om)=\sum_{g\geq 0}a_F(T)\,e^{2\pi i\,\s(T\Om)}$ on $\BH_2$ of
weight $k$ satisfying the following condition
\begin{equation}
a_F(T)=\sum_{d|(n,r,m),\,d>0}d^{k-1}\,a_F \begin{pmatrix}
{{dm}\over {d^2}} & {r \over {2d}} \\ {r \over {2d}} & 1
\end{pmatrix}
\end{equation}
for all $T=\begin{pmatrix} n & {\frac r2} \\
{\frac r2} & m
\end{pmatrix}\geq 0$ with $n,r,m\in\BZ.$ For $F\in M_k(\G_2)$, we
let
\begin{equation*}
F(\Om)=\sum_{m\geq 0}\phi_m(\tau,z)\,e^{2\pi im\tau'},\quad \ \Om=\begin{pmatrix} \tau & z \\
z & \tau'
\end{pmatrix}\in\BH_2
\end{equation*}
be the Fourier-Jacobi expansion of $F$. Then for any nonnegative
integer $m$, we obtain the linear map
\begin{equation*}
\rho_m:M_k(\G_2)\lrt J_{k,m}(\G_1),\ \quad\ F\mapsto \phi_m.
\end{equation*}
We observe that $\rho_0$ is nothing but the Siegel
$\Phi$-operator. Maass \cite{M3,M4,M5} showed that for $k$ even,
there exists a nontrivial map $V:J_{k,1}(\G_1)\lrt M_k(\G_2)$ such
that $\rho_1\circ V$ is the identity. More precisely, we let
$\phi\in J_{k,1}(\G_1)$ be a Jacobi form with Fourier coefficients
$c(n,r)\,(n,r\in \BZ,\ r^2\leq 4n)$ and define for any nonnegative
integer $m\geq 0$
\begin{equation}
\big(V_m\phi\big)(\tau,z)=\sum_{n,r\in\BZ,\,r^2\leq 4mn}\left(
\sum_{d|(n,r,m)}d^{k-1}\,c\left( {{mn}\over {d^2}},{\frac
rd}\right)\right)\,e^{2\pi i(n\tau+rz)}.
\end{equation}
It is easy to see that $V_1\phi=\phi$ and $V_m\phi\in
J_{k,m}(\G_1).$ We define
\begin{equation}
(V\phi)(\Om)=\sum_{m\geq 0}(V_m\phi)(\tau,z)\,e^{2\pi
im\tau'},\quad \Om=\begin{pmatrix} \tau & z \\
z & \tau'
\end{pmatrix}\in \BH_2.
\end{equation}
We denote by $T_n\,(\,n\in\BZ^+)$ the usual Hecke operators on
$M_k(\G_2)$ resp. $S_k(\G_2)$. For instance, if $p$ is a prime,
$T_p=T(p)$ and $T_{p^2}=T_1(p^2)$. We denote by
$T_{J,n}\,(m\in\BZ^+)$ the Hecke operators on $J_{k,m}(\G_1)$
resp. $J_{k,m}^{\textrm{cusp}}(\G_1)$\,(cf.\,\cite{EZ}).

\begin{theorem}$ \textbf{(Maass \cite{M2,M3,M4},\ Eichler-Zagier \cite{EZ}, Theorem 6.3})$
Suppose $k$ is an even positive integer. Then the map $\phi\mapsto
V\phi$ gives an isomorphism of $J_{k,m}(\G_1)$ onto $M_k^*(\G_2)$
which sends cusp Jacobi forms to cusp forms and is compatible with
the action of Hecke operators. If $p$ is a prime, one has
\begin{equation*}
T_p\circ V=V\circ \big(T_{J,p}+p^{k-2}(p+1)\big)
\end{equation*}
and
\begin{equation*}
T_{p^2}\circ V=V\circ
\big(T_{J,p}^2+p^{k-2}(p+1)T_{J,p}+p^{2k-2}\big).
\end{equation*}
\end{theorem}
In Summary, we have the following isomorphisms
\begin{eqnarray}
M_k^*(\G_2) &\cong J_{k,m}(\G_1) &\cong  M_{k-{\frac
12}}^+\big(\G_0^{(1)}(4)\big) \cong M_{2k-2}(\G_1),\\
 V_{\phi} &\longleftarrow\ \ \phi  \ &\lrt f_{\phi} \nonumber
\end{eqnarray}
where the last isomorphism is the Shimura correspondence. All the
above isomorphisms are compatible with the action of Hecke
operators.

\vskip 0.2cm In 1978, providing some evidences, Kurokawa and Saito
conjectured that there is a one-to-one correspondence between
Hecke eigenforms in $S_{2k-2}(\G_1)$ and Hecke eigenforms in
$M_k(\G_2)$ satisfying natural identity between their spinor zeta
functions. This was solved mainly by Maass and then completely
solved by Andrianov \cite{A1} and Zagier \cite{Z}.
\begin{theorem} Suppose $k$ is an even positive integer and let $F\in M_k^*(\G_2)$
be a nonzero Hecke eigenform. Then there exists a unique
normalized Hecke eigenform $f$ in $M_{2k-2}(\G_1)$ such that
\begin{equation}
Z_F(s)=\zeta(s-k+)\,\zeta(s-k+2)L(f,s),
\end{equation}
where $L(f,s)$ is the Hecke $L$-function attached to $f$.
\end{theorem}
$F$ is called the $ \textit{Saito-Kurokawa lift}$ of $f$. Theorem
8.3 implies that $Z_F(s)$ has a pole at $s=k$ if $F$ is an
eigenform in $M_k^*(\G_2)$. If $F\in S_k(\G_2)$ is a Hecke
eigenform, it was proved by Andrianov \cite{A2} that $Z_F(s)$ has
an analytic continuation to the whole complex plane which is
holomorphic  everywhere if $k$ is odd and is holomorphic except
for a possible simple pole at $s=k$ if $k$ is even. Moreover the
global function
\begin{equation*}
Z_F^*(s):=(2\pi)^{-s}\G (s)\G(s-k+2)Z_F(s)
\end{equation*}
is $(-1)^k$-invariant under $s\mapsto 2k-2-s.$ It was proved that
Evdokimov and Oda that $Z_F(s)$ is holomorphic if and only if $F$
is contained in the orthogonal complement of $M_k^*(\G_2)$ in
$M_k(\G_2)$. We remark that $M_k(\G_2)=\BC \,G_k\oplus
S_k^*(\G_2)$, where $G_k$ is the Siegel Eisenstein series of
degree $2$ (cf.\,(6.14)) and $S_k^*(\G_2)=S_k(\G_2)\cap
M_k^*(\G_2).$

\vskip 0.2cm Around 1996, Duke and Imamo${\check g}$lu \cite{DI1}
conjectured a generalization of Theorem 7.3. More precisely, they
formulated the conjecture that if $f$ is a normalized Hecke
eigenform in $S_{2k}(\G_1)\,(k\in\BZ^+)$ and $n$ is a positive
integer with $n\equiv k$\,(mod 2), then there exists a Hecke
eigenform $F$ in $S_{k+n}(\G_{2n})$ such that the standard zeta
function $D_F(s)$ of $F$ equals
\begin{equation}
\zeta(s) \sum_{j=1}^{2n}L(f,s+k+n-j),
\end{equation}
where $L(f,s)$ is the Hecke $L$-function of $f$. Later some
evidence for this conjecture was given by Breulmann and Kuss
\cite{BK}. In 1999, Ikeda \cite{I} proved that the conjecture of
Duke and Imamo${\check g}$lu is true. Such a Hecke eigenform $F$
in $S_{k+n}(\G_{2n})$ is called the $\textit{Duke-Imamo${\check
g}$lu-Ikeda}$ {\it lift} of a normalized Hecke eigenform $f$ in
$S_{2k}(\G_1)$. According to the Shimura isomorphism $M_{k+{\frac
12}}^+\big(\G_0^{(1)}(4)\big) \cong M_{2k}(\G_1)$ in Formula
(8.10), one has the so-called $\textit{Ikeda's lift map}$
\begin{equation}
I_{k,n}:S_{k+{\frac 12}}^+\big(\G_0^{(1)}(4)\big)\lrt
S_{k+n}(\G_{2n})
\end{equation}
defined by
\begin{equation*}
f(\tau)=\sum_{(-1)^km\equiv 0,1( \textrm{mod}\,4)}c(m)e^{2\pi
im\tau}\longmapsto F(\Om)=\sum_{T>0} A(T)\,e^{2\pi i\,\s(T\Om)},
\end{equation*}
where
\begin{eqnarray*} A(T)&=& c(|D_{T,0}|)\,f_T^{k-{\frac
12}}\prod_{p|D_T}{\tilde F}_p(T;\alpha_p)\\
&=&\sum_{a|f_T}a^{k-1}\phi(a;T)\,c(|D_T|/a^2).
\end{eqnarray*}
We refer to \cite{I} and \cite{Ko2} respectively for a precise
definition of $A(T)$. For brevity, we set $S_{k+{\frac
12}}^+=S_{k+{\frac 12}}^+\big(\G_0^{(1)}(4)\big)$. Kohnen and
Kojima \cite{KK} characterized the image of the Ikeda's lift map
$I_{k,n}$. If $F(\Om)=\sum_{T>0} A(T)\,e^{2\pi i\,\s(T\Om)}$ is an
element in the image of $I_{k,n}$, then $A(T)=A({\tilde T})$ if
$T$ and ${\tilde T}$ are positive definite half-integral matrices
of degree $2n$ with $D_T$ and $D_{\tilde T}$ and such that for all
positive divisors $a$ of $f_T=f_{\tilde T}$, one has
$\phi(a;T)=\phi(a;{\tilde T}).$ Here $D_T$ denotes the
discriminant of $T$ defined by $D_T:=(-1)^n\det (2T).$ They called
the image of $I_{k,n}$ in $S_{k+n}(\G_{2n})$ the $ \textit{Maass
space}$. If $n=1,\ M^*_k(\G_2)$ coincides with the image of
$I_{k,1}$. Thus this generalizes the original Maass space.
Breulmann and Kuss \cite{BK} dealt with the special case of the
lift map $I_{6,2}:S_{12}(\G_1)(\cong S_{13/2}^+)\lrt S_8(\G_4)$.
In the article \cite{BFW}, starting with the Leech lattice
$\Lambda$, the authors constructed a nonzero Siegel cusp form of
degree $12$ and weight $12$ which is the image of a cusp form
$\Delta\in S_{12}(\G_1)$ under the Ikeda lift map $I_{6,6}$. Here
$\Delta$ is the cusp form in $S_{12}(\G_1)$ defined by
\begin{equation*}
\Delta(\tau)=(2\pi)^{12}\, q\,
\prod_{n=1}^{\infty}(1-q^n)^{24},\quad \tau\in\BH_1,\ q=e^{2\pi
i\tau}.
\end{equation*}
It is known that there exist $24$ Niemeier lattices of rank $24$,
say, $L_1,\cdots,L_{24}$. The theta series
\begin{equation*}
\theta_{L_i}(\Om)=\sum_{G\in \BZ^{(24,12)}}e^{2\pi i\,
\s(L_i[G]\Om)},\quad \Om\in \BH_{12},\quad i=1,\cdots,24
\end{equation*}
generate a subspace $V_*$ of $M_{12}(\G_{12})$. These
$\theta_{L_i}\,(1\leq i\leq 24)$ are linearly independent. It can
be seen that the intersection $V_*\cap S_{12}(\G_{12})$ is one
dimensional. This nontrivial cusp form in $V_*\cap
S_{12}(\G_{12})$ up to constant is just the Siegel modular form
constructed by them. Under the assumption $n+r\equiv k\, (
\textrm{mod}\,2)$ with $k,n,r\in\BZ^+,$ using the lift map
$I_{k,n+r}:S_{k+{\frac 12}}^+\lrt S_{k+n+r}(\G_{2n+2r}),$ recently
Ikeda \cite{I2} constructed the following map
\begin{equation}
J_{k,n,r}:S_{k+{\frac 12}}^+\times S_{k+n+r}(\G_r)\lrt
S_{k+n+r}(\G_{2n+r})
\end{equation}
defined by
\begin{equation*}
J_{k,n,r}(h,G)(\Om):=\int_{\G_r\ba \BH_r}I_{k,n+r}(h)\left(
\begin{pmatrix} \Om & 0 \\ 0 & \tau \end{pmatrix}\right)
\overline{G^c(\tau)}\, \big( \det
\textrm{Im}\,\tau\big)^{k+n-1}d\tau,
\end{equation*}
where $ h\in S_{k+{\frac 12}}^+,\,G\in  S_{k+n+r}(\G_r),\,
\Om\in\BH_{2n+r},\,\tau\in \BH_r,\
G^c(\tau)=\overline{G(-\overline{\tau})}$ and $( \det
\textrm{Im}\,\tau)^{-(r+1)}d\tau$
is an invariant volume element\,(cf.\,\S2\,(2.3)). He proved that
the standard zeta function $D_{J_{k,n,r}(h,G)}(s)$ of
$J_{k,n,r}(h,G)$ is equal to
\begin{equation*}
D_{J_{k,n,r}(h,G)}(s)=D_G(s)\,\prod_{j=1}^n L(f,s+k+n-j),
\end{equation*}
where $f$ is the Hecke eigenform in $S_{2k}(\G_1)$ corresponding
to $h\in S_{k+{\frac 12}}^+$ under the Shimura correspondence.
\vskip 0.2cm \noindent $ \textbf{Question\,:}$ Can you describe a
geometric interpretation of the Duke-Imamo${\check g}$lu-Ikeda
lift or the map $J_{k,n,r}$ ?

\end{section}

%
%
\begin{section}{\bf Holomorphic Differential Forms on Siegel Space}
\setcounter{equation}{0} \vskip 0.5cm In this section, we describe
the relationship between Siegel modular forms and holomorphic
differential forms on the Siegel space. We also discuss the Hodge
bundle. First of all we need to review the theory of toroidal
compactifications of the Siegel space. \vskip 0.13cm Let $\BD_g$
be the generalized unit disk of degree $g$ (cf.\,Section 2).
$\BH_g$ is realized as a bounded symmetric domain via the Cayley
transform (cf.\,(2.7)). Let ${\overline \BD}_g$ be the topological
closure in $T_g$, where $T_g$ denotes the vector space of all
$g\times g$ complex symmetric matrices. Then ${\overline \BD}_g$
is the disjoint union of all boundary components of $\BD_g$. Let
\begin{equation*}
F_r:=\left\{ \begin{pmatrix} W_1 & 0 \\ 0 & I_{g-r}
\end{pmatrix}\in {\overline \BD}_g\,\Big| \ W_1 \in
\BD_r\,\right\},\quad r=0,1,\cdots,g
\end{equation*}
be the standard rational boundary components of $\BD_g$. Then any
boundary component $F$ of $\BD_g$ is of the form $F=\s\cdot F_r$
for some $\s\in Sp(g,\BR)$ and some integer $r$ with $0\leq r\leq
g.$ In addition, if $F$ is a rational boundary component of
$\BD_g$, then it is of the form $F=\g\cdot F_r$ for some $\g\in
\G_g$ and some integer $r$ with $0\leq r\leq g.$ We note that
$F_0=\{ I_g\}$ and $F_g=\BD_g.$ We set
\begin{equation}
\BD^*_g:=\bigcup_{0\leq r\leq g}\G_g\cdot F_r.
\end{equation}
Then $\BD^*_g$ is the union of all rational boundary components of
$\BD_g$ and is called the $ \textit{rational closure}$ of $\BD_g$.
Then we obtain the so-called $ \textit{Satake compactification}\
A_g^*:=\G_g\backslash \BD_g^*$ of $A_g:=\G_g\backslash \BD_g.$ Let
$F$ be a rational boundary component of $\BD_g$. We let
$P(F),\,W(F)$ and $U(F)$ be the parabolic subgroup associated with
$F$,\ the unipotent radical of $P(F)$ and the center of $W(F)$
respectively. We set $V(F):=W(F)/U(F).$ Since $P(\s\cdot F)=\s
P(F)\s^{-1}$ for $\s\in Sp(g,\BR),$ it is enough to investigate
the structures of these groups for the standard rational boundary
components $F_r$ with $0\leq r\leq g.$ \vskip 0.1cm Now we take
$F=F_r$ for some integer $r$ with $0\leq r\leq g.$ We define
$D(F):=U(F)_{\BC}\cdot \BD_g\subset {\hat \BD}_g$. Here ${\hat
\BD}_g:=B\backslash Sp(g,\BR)_{\BC}$ is the compact dual of
$\BD_g$ with $B$ a parabolic subgroup of $Sp(g,\BR)_{\BC}$. We
denote by $G_{\BC}$ the complexification of a real Lie group $G$.
It is obvious that $U(F)_{\BC}\cong T_{g-r}$ and $D(F)\cong
F\times V(F)\times U(F)_{\BC}$ analytically. We observe that
$U(F)$ acts on $D(F)$ as the linear translation on the factor
$U(F)_{\BC}.$ Indeed, the isomorphism $\varphi:D(F)\lrt F\times
V(F)\times U(F)_{\BC}$ is given by
\begin{equation}
\varphi\left( \begin{pmatrix} W_1 & W_2 \\ * & W_3
\end{pmatrix}\right):=(W_1,W_2,W_3),
\end{equation}
where $W_1\in\BD_r,\ W_2\in\BC^{(r,g-r)}$ and $W_3\in T_{g-r}.$ We
define the mapping $\Phi_F:D(F)\lrt U(F)$ by
\begin{equation}
\Phi_F\left( \begin{pmatrix} W_1 & W_2 \\ * & W_3
\end{pmatrix}\right)= \textrm{Im}\,W_3-{}^t( \textrm{Im}\,W_2)( \textrm{Im}\,W_1)^{-1}( \textrm{Im}\,W_2),
\end{equation}
where $(W_1,W_2,W_3)\in D(F)$. Then $\BD_g\cong \BH_g$ is
characterized by the condition $\Phi_F(W)>0$ for all $W\in \BD_g$.
This is the realization of a Siegel domain of the third kind. We
let $C(F)$ be the cone of real positive symmetric matrices of
degree $g-r$ in $U(F)\cong T_{g-r}(\BR),$ where $T_{g-r}(\BR)$
denotes the vector space of all symmetric real matrices of degree
$g-r$. Clearly one has $\BD_g=\Phi^{-1}(C(F)).$ We define
\begin{equation*}
G_h(F):= \textrm{Aut}(F)\quad ( \textrm{modulo\ finite\ group})
\end{equation*}
and
\begin{equation*}
G_l(F):= \textrm{Aut}\big(U(F),C(F)\big).
\end{equation*}
Then it is easy to see that
\begin{equation*}
P(F)=\big( G_h(F)\times G_l(F)\big) \ltimes W(F)\ \quad (
\textrm{the\ semi-direct\ product)}.
\end{equation*}
We obtain the natural projections $p_h:P(F)\lrt G_h(F)$ and
$p_l:P(F)\lrt G_l(F).$

\vskip 0.3cm\noindent $ \textbf{Step I\,:}$ Partial
compactification for a rational boundary component. \vskip 0.2cm
Now we let $\G$ be an arithmetic subgroup of $Sp(g,\BR)$. We let
\begin{eqnarray*}
\G(F):&=& \G\cap P(F),\\
{\overline \G}(F):&=& p_l(\G(F))\subset G_l(F),\\
U_{\G}(F):&=&\G \cap U(F),\quad \textrm{a lattice in} \ U(F),\\
W_{\G}(F):&=& \G\cap W(F).
\end{eqnarray*}
We note that ${\overline \G}(F)$ is an arithmetic subgroup of
$G_l(F)$. \vskip 0.2cm Let $\Sigma_F=\big\{ \s_{\al}^F\big\}$ be a
${\overline \G}(F)$-admissible polyhedral decomposition of $C(F)$.
We set $D(F)':=D(F)/U(F)_{\BC}.$ Since $D(F)'\cong F\times V(F)$,
the projection $\pi_F:D(F)\lrt D(F)'$ is a principal
$U(F)_{\BC}$-bundle over $D(F)'.$ The map
\begin{equation}
\pi_{F,\G}:U_{\G}(F)\ba D(F)\cong F\times V(F)\times
\big(U_{\G}(F)\ba U(F)_{\BC}\big)\lrt D(F)'
\end{equation}
is a principal $T(F)$-bundle with the structure group
$T(F):=U_{\G}(F)\ba U(F)_{\BC}\cong \big(\BC^*\big)^q,$ where
$q={\frac 12}(g-r)(g-r+1).$ Let $X_{\Sigma_F}$ be a normal torus
embedding of $T(F)$. We note that $X_{\Sigma_F}$ is determined by
$\Sigma_F$. Then we obtain a fibre bundle
\begin{equation}
{\mathfrak X}(\Sigma_F):=\big( U_{\G}(F)\ba
D(F)\big)\times_{T(F)}X_{\Sigma_F}
\end{equation}
over $D(F)'$ with fibre $X_{\Sigma_F}$. We denote by $
\textbf{X}(\Sigma_F)$ the interior of the closure of $U_{\G}(F)\ba
\BD_g$ in ${\mathfrak X}(\Sigma_F)$ (because $\BD_g\subset D(F)$).
$\textbf{X}(\Sigma_F)$ has a fibrewise $T(F)$-orbit decomposition
$\coprod_{\mu}O(\mu)$ such that
\begin{eqnarray*}
 &&\textrm{(a)\ each}\ O(\mu)\ \textrm{is an algebraic torus bundle over}\ D(F)' ,\\
&& \textrm{(b)}\ \s_{\mu}^F\prec \s_{\nu}^F\quad \textrm{iff}\ {\overline{O(\mu)}}\supseteq O(\nu),   \\
&&\textrm{(c)}\ \textrm{dim}\,\s_{\mu}^F+ \dim O(\mu)= \dim D(F),\\
&&\textrm{(d)}\ \textrm{for}\ \s_{\mu}^F=0,\ O(\mu)=U_{\G}(F)\ba
D(F).
\end{eqnarray*}
We define
\begin{equation*}
O(F):=\bigcup_{\s_{\al}^F\cap \,C(F)\neq \emptyset}O(\al)\subset
\textbf{X}(\Sigma_F)
\end{equation*}
and
\begin{equation*}
{\bar O}(F):=\G(F)/U_{\G}(F)\ba O(F).
\end{equation*}
We note that $O(F_g)=\BD_g$ and ${\bar O}(F_g)=\G\ba \BD_g.$ We
set
\begin{equation*}
\textbf{Y}(\Sigma_F):=\G(F)/U_{\G}(F)\ba \textbf{X}(\Sigma_F).
\end{equation*}
We note that $\G(F)/U_{\G}(F)$ acts on $\textbf{Y}(\Sigma_F)$
properly discontinuously. Then we can show that
$\textbf{Y}(\Sigma_F)$ has a canonical quotient structure of a
normal analytic space and ${\bar O}(F)$ is a closed analytic set
in $\textbf{Y}(\Sigma_F)$. \vskip 0.3cm\noindent $ \textbf{Step
II\,:}$ Gluing. \vskip 0.2cm Let $\Sigma:=\big\{ \Sigma_F\,|\ F\
\textrm{is a rational boundary component of}\ \BD_g\,\big\}$ be a
$\G$-admissible family of polyhedral decompositions. We put
\begin{equation*}
{\widetilde {\G\ba
\BD_g}}:=\bigcup_{F:\,\textrm{rational}}\textbf{X}(\Sigma_F).
\end{equation*}
We define the equivalence relation $\backsim $ on
${\widetilde{\G\ba \BD_g}}$ as follows:
\begin{equation*}
X_1\,\backsim\,X_2,\ \ X_1\in {\bold X}(\Sigma_{F_1}),\ \ X_2\in
{\bold X}(\Sigma_{F_2})
\end{equation*}
iff there exist a rational boundary component $F$, an element
$\gamma\in \G$ such that $F_1\prec F,\ \gamma\,F_2\prec F$ and
there exists an element $X\in {\bold X}(\Sigma_F)$ such that
$\pi_{F,F_1}(X)=X_1,\ \pi_{F,F_2}(X)= \gamma X_2,$ where
$$\pi_{F,F_1}:{\bold X}(\Sigma_F)\lrt {\bold X}(\Sigma_{F_1}),\ \ \
\pi_{F,F_2}:{\bold X}(\Sigma_F)\lrt {\bold X}(\Sigma_{\gamma
F_2}).$$ The space $\overline{ \G\ba \BD_g}:= {\widetilde{ \G\ba
\BD_g}}/\backsim$ is called the {\it toroidal\ compactification}
of $\G\ba \BD_g$ associated with $\Sigma$. It is known that
$\overline{\left( \G\ba \BD_g\right)}$ is a Hausdorff analytic
variety containing $\G\ba \BD_g$ as an open dense subset.

For a neat arithmetic subgroup $\G$, e.g., $\G=\G_g(n)$ with
$n\geq 3$, we can obtain a smooth projective toroidal
compactification of $\G\ba \BD_g.$ The theory of toroidal
compactifications of bounded symmetric domains was developed by
Mumford's school\,(cf.\,\cite{AMRT} and \cite{N}). We set
\begin{equation*}
{\mathcal A}_g:=\G_g\ba \BH_g\quad \textrm{and}\quad {\mathcal
A}^*_g:=\G_g\ba \BH_g^*=\bigcup_{0\leq i\leq g} \G_i\ba \BH_i\quad
(\textrm{disjoint union}).
\end{equation*}
I. Satake \cite{Sa1} showed that ${\mathcal A}^*_g$ is a normal
analytic space and W. Baily \cite{B1} proved that ${\mathcal
A}^*_g$ is a projective variety. Let ${\tilde {\mathcal A}}_g$ be
a toroidal compactification of ${\mathcal A}_g$. Then the boundary
${\tilde {\mathcal A}}_g-{\mathcal A}_g$ is a divisor with normal
crossings and one has a universal semi-abelian variety over
${\tilde {\mathcal A}}_g$ in the orbifold. We refer to \cite{HS}
for the geometry of ${\mathcal A}_g$.

\vskip 0.2cm Let $\theta$ be the second symmetric power of the
standard representation of $GL(g,\BC)$. For brevity we set
$N={\frac 12}g(g+1).$ For an integer $p$ with $0\leq p\leq N$, we
denote by $\theta^{[p]}$ the $p$-th exterior power of $\theta$.
For any integer $q$ with $0\leq q\leq N$, we let
$\Om^q(\BH_g)^{\G_g}$ be the vector space of all $\G_g$-invariant
holomorphic $q$-forms on $\BH_g$. Then we obtain an isomorphism
\begin{equation*}
\Om^q(\BH_g)^{\G_g}\lrt M_{\theta^{[q]}}(\G_g).
\end{equation*}

\begin{theorem} $\textbf{(Weissauer \cite{W1})}$ For an integer
$\al$ with $0\leq \al\leq g$, we let $\rho_{\al}$ be the
irreducible representation of $GL(g,\BC)$ with the highest weight
\begin{equation*}
(g+1,\cdots,g+1,g-\al,\cdots,g-\al)
\end{equation*}
such that $\text{corank}\,(\rho_{\al})=\al$ for $1\leq \al\leq g.$
If $\al=-1,$ we let $\rho_{\al}=(g+1,\cdots,g+1).$ Then
\begin{equation*}
\Omega^q(\BH_g)^{\G_g}=\begin{cases} M_{\rho_{\al}}(\G_g)\ \ &
\textrm{if}\ \,
q={{g(g+1)}\over 2}-{{\al(\al+1)}\over 2} \\
0 &\text{otherwise}.\end{cases}
\end{equation*}
\end{theorem}

\noindent{\bf Remark.} If $2\al>g,$ then any $f\in
M_{\rho_{\al}}(\G_g)$ is singular (cf.\,Theorem 5.4). Thus if
$q<{{g(3g+2)}\over 8},$ then any $\G_g$-invariant holomorphic
$q$-form on $\BH_g$ can be expressed in terms of vector valued
theta series with harmonic coefficients. It can be shown with a
suitable modification that the just mentioned statement holds for
a sufficiently small congruence subgroup of $\G_g.$
\vspace{0.1in}\\
\indent Thus the natural question is to ask how to determine the
$\G_g$-invariant holomorphic $p$-forms on $\BH_g$ for the
nonsingular range $\dfrac{g(3g+2)}{8}\leq p \leq
\dfrac{g(g+1)}{2}.$ Weissauer \cite{W2} answered the above
question for $g=2.$ For $g>2,$ the above question is still open.
It is well known that the vector space of vector valued modular
forms of type $\rho$ is finite dimensional. The computation or the
estimate of the dimension of $\Omega^p(\BH_g)^{\G_g}$ is
interesting because its dimension is finite even though the
quotient space ${\mathcal A}_g$ is noncompact.
\vspace{0.1in}\\
{\bf Example 1.} Let
\begin{equation}
\varphi=\sum_{i\leq j}f_{ij}(\Om)\,d\omega_{ij}
\end{equation}
be a $\G_g$-invariant holomorphic $1$-form on $\BH_g$. We put
\begin{equation*}
f(\Om)=\big( f_{ij}(\Om)\big)\ \textrm{with} \ f_{ij}=f_{ji}\quad
\textrm{and}\quad d\Om=(d\om_{ij}).
\end{equation*}
Then $f$ is a matrix valued function on $\BH_g$ satisfying the
condition
\begin{equation*}
f(\g\cdot\Om)=(C\Om+D)f(\Om)\,^t(C\Om+D)\quad \textrm{for all}\
\g=\begin{pmatrix} A & B\\ C & D \end{pmatrix}\in \G_g \
\textrm{and}\ \Om\in\BH_g.
\end{equation*}
This implies that $f$ is a Siegel modular form in
$M_{\theta}(\G_g)$, where $\theta$ is the irreducible
representation of $GL(g,\BC)$ on $T_g=\textrm{Symm}^2(\BC^g)$
defined by
\begin{equation*}
\theta(h)v=h\,v\,{}^th,\quad h\in GL(g,\BC),\ v\in T_g.
\end{equation*}
We observe that $(9.6)$ can be expressed as $\varphi=\s(f\,d\Om).$
\vspace{0.1in}\\
{\bf Example 2.} Let
\begin{equation*}
\om_0=d\om_{11}\wedge d\om_{12}\wedge \cdots \wedge d\om_{gg}
\end{equation*}
be a holomorphic $N$-form on $\BH_g$. If $\om=f(\Om)\,\om_0$ is
$\G_g$-invariant, it is easily seen that
\begin{equation*}
f(\g\cdot\Om)=\det(C\Om+D)^{g+1}f(\Om)\quad \textrm{for all}\
\g=\begin{pmatrix} A & B\\ C & D \end{pmatrix}\in \G_g \
\textrm{and}\ \Om\in\BH_g.
\end{equation*}
Thus $f\in M_{g+1}(\G_g).$ It was shown by Freitag \cite{Fr3} that
$\om$ can be extended to a holomorphic $N$-form on
${\tilde{\mathcal A}}_g$ if and only if $f$ is a cusp form in
$S_{g+1}(\G_g).$ Indeed, the mapping
\begin{equation*}
S_{g+1}(\G_g)\lrt \Om^N\big( {\tilde{\mathcal
A}}_g\big)=H^0\big({\tilde{\mathcal A}}_g,\Om^N\big),\qquad
f\mapsto f\,\om_0
\end{equation*}
is an isomorphism. Let $\om_k=F(\Om)\,\om_0^{\otimes k}$ be a
$\G_g$-invariant holomorphic form on $\BH_g$ of degree $kN$. Then
$F\in M_{k(g+1)}(\G_g).$
\vspace{0.1in}\\
{\bf Example 3.} We set
\begin{equation*}
\eta_{ab}=\epsilon_{ab}\bigwedge_{\substack{1\leq\mu\leq\nu\leq g\\
(\mu,\nu)\neq (a,b)}}d\om_{\mu\nu},\quad 1\leq a\leq b\leq g,
\end{equation*}
where the signs $\epsilon_{ab}$ are determined by the relations
$\epsilon_{ab}\,\eta_{ab}\wedge d\om_{ab}=\om_0.$ We assume that
\begin{equation*}
\eta_*=\sum_{1\leq a\leq b\leq g}F_{ab}\,\eta_{ab}
\end{equation*}
is a $\G_g$-invariant holomorphic $(N-1)$-form on $\BH_g$. Then
the matrix valued function $F=\big(\epsilon_{ab}\,F_{ab}\big)$
with $\epsilon_{ab}=\epsilon_{ba}$ and $F_{ab}=F_{ba}$ is an
element of $M_{\tau}(\G_g)$, where $\tau$ is the irreducible
representation of $GL(g,\BC)$ on $T_g$ defined by
\begin{equation*}
\tau(h)v=( \det h)^{g+1}\,^th^{-1}vh^{-1},\quad h\in GL(g,\BC),\
v\in T_g.
\end{equation*}
\vskip 0.2cm
\indent We will mention the results due to Weissauer \cite{W2}. We
let $\G$ be a congruence subgroup of $\G_2.$ According to Theorem
9.1, $\G$-invariant holomorphic forms in $\Omega^2(\BH_2)^{\G}$
are corresponded to modular forms of type (3,1). We note that
these invariant holomorphic $2$-forms are contained in the {\it
nonsingular\ range}. And if these modular forms are not cusp
forms, they are mapped under the Siegel $\Phi$-operator to cusp
forms of weight $3$ with respect to some congruence
subgroup\,(\,dependent on $\G$\,) of the elliptic modular group.
Since there are finitely many cusps, it is easy to deal with these
modular forms in the adelic version. Observing these facts, he
showed that any $2$-holomorphic form on $\G\ba \BH_2$ can be
expressed in terms of theta series with harmonic coefficients
associated to binary positive definite quadratic forms. Moreover
he showed that $H^2(\G\ba \BH_2,\BC)$ has a pure Hodge structure
and that the Tate conjecture holds for a suitable compactification
of $\G\ba \BH_2.$  If $g\geq 3,$ for a congruence subgroup $\G$ of
$\G_g$ it is difficult to compute the cohomology groups
$H^{\ast}(\G\ba \BH_g,\BC)$ because $\G\ba \BH_g$ is noncompact
and highly singular. Therefore in order to study their structure,
it is natural to ask if they have pure Hodge structures or mixed
Hodge structures.

\vskip 0.2cm We now discuss the Hodge bundle on the Siegel modular
variety ${\mathcal A}_g$. For simplicity we take $\G=\G_g(n)$ with
$n\geq 3$ instead of $\G_g$. We recall that $\G_g(n)$ is a
congruence subgroup of $\G_g$ consisting of matrices $M\in \G_g$
such that $M\equiv I_{2g}\,( \textrm{mod}\ n).$ Let
\begin{equation*}
{\mathfrak X}_g(n):=\G_g(n)\ltimes \BZ^{2g}\ba \BH_g\times \BC^g
\end{equation*}
be a family of abelian varieties of dimension $g$ over ${\mathcal
A}_g(n):=\G_g(n)\ba \BH_g.$ We recall that $\G_g(n)\ltimes
\BZ^{2g}$ acts on $\BH_g\times \BC^g$ freely by
\begin{equation*}
\big(\g,(\la,\mu)\big)\cdot(\Om,Z)=\big(\g\cdot\Om,(Z+\la
\Om+\mu)(C\Om+D)^{-1}\big),
\end{equation*}
where $\g=\begin{pmatrix} A & B\\ C & D \end{pmatrix}\in \G_g(n),\
\la,\mu\in\BZ^g,\ \Om\in \BH_g$ and $Z\in\BC^g.$ If we insist on
using $\G_g$, we need to work with orbifolds or stacks to have a
universal family
\begin{equation*}
{\mathfrak X}_g:={\mathfrak X}_g(n)/Sp(g,\BZ/n\BZ)
\end{equation*}
available. We observe that $\G_g(n)$ acts on $\BH_g$ freely.
Therefore we obtain a vector bundle ${\mathbb E}={\mathbb E}_g$
over ${\mathcal A}_g(n)$ of rank $g$
\begin{equation*}
{\mathbb E}={\mathbb E}_g:=\G_g(n)\ba \big(\BH_g \times
\BC^g\big).
\end{equation*}
This bundle ${\mathbb E}$ is called the $ \textit{Hodge bundle}$
over ${\mathcal A}_g(n)$. The finite group $Sp(g,\BZ/n\BZ)$ acts
on ${\mathbb E}$ and a $Sp(g,\BZ/n\BZ)$-invariant section of
$(\det {\mathbb E})^{\otimes k}$ with a positive integer $k$ comes
from a Siegel modular form of weight $k$ in $M_k(\G_g)$. The
canonical line bundle $\kappa_g(n)$ of ${\mathcal A}_g(n)$ is
isomorphic to $ (\det {\mathbb E})^{\otimes (g+1)}$. A holomorphic
section of $\kappa_g(n)$ corresponds to a Siegel modular form in
$M_{g+1}(\G_g(n))$\,(cf. Example 2). We note that the sheaf
$\Om_{{\mathcal A}_g(n)}^1$ of holomorphic $1$-forms on ${\mathcal
A}_g(n)$ is isomorphic to $ \textrm{Symm}^2({\mathbb E}).$ This
sheaf can be extended over a toroidal compactification ${\tilde
{\mathcal A}}_g$ of ${\mathcal A}_g$ to an isomorphism
\begin{equation*}
\Om_{{\tilde {\mathcal A}}_g}^1(\log D)\cong
\textrm{Symm}^2({\mathbb E}),
\end{equation*}
where the boundary $D={\tilde {\mathcal A}}_g-{\mathcal A}_g$ is
the divisor with normal crossings. Similarly to each finite
dimensional representation $(\rho,V_{\rho})$ of $GL(g,\BC)$, we
may associate the vector bundle
\begin{equation*}
{\mathbb E}_{\rho}:=\G_g(n)\ba \big(\BH_g \times V_{\rho}\big)
\end{equation*}
by identifying $(\Om,v)$ with $(\g\cdot \Om,\rho(C\Om+D)v)$, where
$\Om\in\BH_g,\ v\in V_{\rho}$ and $\g=\begin{pmatrix} A & B\\ C &
D \end{pmatrix}\in \G_g(n).$ Obviously ${\mathbb E}_{\rho}$ is a
holomorphic vector bundle over ${\mathcal A}_g(n)$ of rank $ \dim
V_{\rho}.$

\end{section}

%
%
\begin{section}{{\bf Subvarieties of the Siegel Modular
Variety}} \setcounter{equation}{0}
\renewcommand{\theequation}{10.\arabic{equation}}
\newcommand\BP{\mathbb P}
\vskip 0.3cm  Here we assume that the ground field is the complex
number field $\BC.$
\vspace{0.1in}\\
\noindent{\bf Definition 9.1.}\quad A nonsingular variety $X$ is
said to be {\it rational} if $X$ is birational to a projective
space $\BP^n(\BC)$ for some integer $n$. A nonsingular variety $X$
is said to be {\it stably\ rational} if $X\times \BP^k(\BC)$ is
birational to $\BP^N(\BC)$ for certain nonnegative integers $k$
and $N$. A nonsingular variety $X$ is called {\it unirational} if
there exists a dominant rational map $\varphi:\BP^n(\BC)\lrt X$
for a certain positive integer $n$, equivalently if the function
field $\BC(X)$ of $X$ can be embedded in a purely transcendental
extension $\BC(z_1,\cdots,z_n)$ of $\BC.$
\vspace{0.1in}\\
\noindent{\bf Remarks 9.2.}\quad (1) It is easy to see that the
rationality implies the stably rationality and that the stably
rationality implies the unirationality.
\vspace{0.05in}\\
\noindent (2) If $X$ is a Riemann surface or a complex surface,
then the notions of rationality, stably rationality and
unirationality are equivalent one another.
\vspace{0.05in}\\
\noindent (3) Griffiths and Clemens \cite{C-G} showed that most of
cubic threefolds in $\BP^4(\BC)$ are unirational but {\it not}
rational.
\vspace{0.1in}\\
The following natural questions arise :
\vspace{0.1in}\\
{\sc Question 1.}\quad Is a stably rational variety {\it
rational}\,? Indeed, the question was raised by Bogomolov.
\vspace{0.1in}\\
{\sc Question 2.}\quad Is a general hypersurface $X\subset
\BP^{n+1}(\BC)$ of degree $d\leq n+1$ {\it unirational}\,?
\vspace{0.1in}\\
\noindent{\bf Definition 9.3.}\quad Let $X$ be a nonsingular
variety of dimension $n$ and let $K_X$ be the canonical divisor of
$X$. For each positive integer $m\in \BZ^+$, we define the
$m$-{\it genus} $P_m(X)$ of $X$ by
$$P_m(X):=\text{dim}_{\BC}\,H^0(X,{\mathcal{O}}(mK_X)).$$
The number $p_g(X):=P_1(X)$ is called the {\it geometric\ genus}
of $X$. We let
$$N(X):=\left\{\,m\in \BZ^+\,\vert\, P_m(X)\geq 1\,\right\}.$$
For the present, we assume that $N(X)$ is nonempty. For each $m\in
N(X),$ we let $\left\{ \phi_0,\cdots,\phi_{N_m}\right\}$ be a
basis of the vector space $H^0(X,{\mathcal{O}}(mK_X)).$ Then we
have the mapping $\Phi_{mK_X}\,: X\lrt \BP^{N_m}(\BC)$ by
$$\Phi_{mK_X}(z):=(\phi_0(z):\cdots:\phi_{N_m}(z)),\ \ \ z\in X.$$
We define the {\it Kodaira\ dimension} $\kappa(X)$ of $X$ by
$$\kappa(X):=\text{max}\,\left\{\,\text{dim}_{\BC}\,\Phi_{mK_X}(X)\,\vert\,\,
m\in N(X)\,\right\}.$$ If $N(X)$ is empty, we put
$\kappa(X):=-\infty.$ Obviously $\kappa(X)\leq
\text{dim}_{\BC}\,X.$ A nonsingular variety $X$ is said to be {\it
of\ general\ type} if $\kappa(X)=\text{dim}_{\BC}X.$ A singular
variety $Y$ in general is said to be rational, stably rational,
unirational or of general type if any nonsingular model $X$ of $Y$
is rational, stably rational, unirational or of general type
respectively. We define
$$P_m(Y):=P_m(X)\ \ \ \ \text{and}\ \ \ \ \kappa(Y):=\kappa(X).$$
A variety $Y$ of dimension $n$ is said to be {\it of\ logarithmic\
general\ type} if there exists a smooth compactification ${\tilde
Y}$ of $Y$ such that $D:={\tilde Y}-Y$ is a divisor with normal
crossings only and the transcendence degree of the logarithmic
canonical ring
$$\oplus_{m=0}^{\infty}\,H^0({\tilde Y},\,m(K_{\tilde Y}+[D]))$$
is $n+1$, i.e., the {\it logarithmic\ Kodaira\ dimension} of $Y$
is $n$. We observe that the notion of being of logarithmic general
type is weaker than that of being of general type.
\newcommand\CA{\mathcal A}
\vspace{0.1in}\\
\indent Let $\CA_g:=\G_g\ba \BH_g$ be the Siegel modular variety
of degree $g$, that is, the moduli space of principally polarized
abelian varieties of dimension $g$. It has been proved that
$\CA_g$ is of general type for $g\geq 6.$ At first Freitag
\cite{Fr1} proved this fact when $g$ is a multiple of $24$. Tai
\cite{Ta} proved this fact for $g\geq 9$ and Mumford \cite{Mf4} proved
this fact for $g\geq 7.$
Recently Grushevsky
and Lehavi \cite{G-L} announced that they proved that the Siegel
modular variety ${\mathcal A}_6$ of genus $6$ is of general type
after constructing a series of new effective geometric divisors on
${\mathcal A}_g.$ Before 2005 it had been known that ${\mathcal
A}_g$ is of general type for $g\geq 7$. On the other hand, $\CA_g$
is known to be unirational for $g\leq 5\,:$ Donagi \cite{Do} for
$g=5,$ Clemens \cite{Cl} for $g=4$ and classical for $g\leq 3.$
For $g=3,$ using the moduli theory of curves, Riemann \cite{Ri},
Weber \cite{We} and Frobenius \cite{Fro} showed that
$\CA_3(2):=\G_3(2)\ba \BH_3$ is a rational variety and moreover
gave $6$ generators of the modular function field $K(\G_3(2))$
written explicitly in terms of derivatives of odd theta functions
at the origin. So $\CA_3$ is a unirational variety with a Galois
covering of a rational variety of degree
$[\G_3:\G_3(2)]=1,451,520.$ Here $\G_3(2)$ denotes the principal
congruence subgroup of $\G_3$ of level $2.$ Furthermore it was
shown that $\CA_3$ is stably rational(cf. \cite{KS}, \cite{Bo}).
For a positive integer $k$, we let $\G_g(k)$ be the principal
congruence subgroup of $\G_g$ of level $k$. Let $\CA_g(k)$ be the
moduli space of abelian varieties of dimension $g$ with $k$-level
structure. It is classically known that $\CA_g(k)$ is of
logarithmic general type for $k\geq 3$\,(cf. \cite{Mf3}). Wang
\cite{Wa2} proved that $\CA_2(k)$ is of general type for $k\geq
4.$ On the other hand, van der Geer \cite{G1} showed that
$\CA_2(3)$ is rational. The remaining unsolved problems are
summarized as follows\,:
\vspace{0.1in}\\
{\bf Problem 1.}\quad Is $\CA_3$ rational\,?
\vspace{0.1in}\\
{\bf Problem 2.}\quad Are $\CA_4,\ \CA_5$ stably rational or
rational\,?
\vspace{0.1in}\\
{\bf Problem 3.}\quad What type of varieties are $\CA_g(k)$ for
$g\geq 3$ and $k\geq 2$\,?
\vspace{0.1in}\\
\indent We already mentioned that $\CA_g$ is of general type if
$g\geq 6.$ It is natural to ask if the subvarieties of
$\CA_g\,(g\geq 6)$ are of general type, in particular the
subvarieties of $\CA_g$ of codimension one. Freitag \cite{Fr4}
showed that there exists a certain bound $g_0$ such that for
$g\geq g_0,$ each irreducible subvariety of $\CA_g$ of codimension
one is of general type. Weissauer \cite{W3} proved that every
irreducible divisor of $\CA_g$ is of general type for $g\geq 10.$
Moreover he proved that every subvariety of codimension $\leq
g-13$ in $\CA_g$ is of general type for $g\geq 13.$ We observe
that the smallest known codimension for which there exist
subvarieties of $\CA_g$ for large $g$ which are not of general
type is $g-1.\ \CA_1\times \CA_{g-1}$ is a subvariety of $\CA_g$
of codimension $g-1$ which is not of general type.
\vspace{0.1in}\\
\noindent{\bf Remark.}\quad Let $\M_g$ be the coarse moduli space
of curves of genus $g$ over $\BC.$ Then $\M_g$ is an analytic
subvariety of $\CA_g$ of dimension $3g-3.$ It is known that $\M_g$
is unirational for $g\leq 10.$ So the Kodaira dimension
$\kappa(\M_g)$ of $\M_g$ is $-\infty$ for $g\leq 10.$ Harris and
Mumford \cite{H-M} proved that $\M_g$ is of general type  for odd
$g$ with $g\geq 25$ and $\kappa(\M_{23})\geq 0.$

\end{section}

%
%
\begin{section}{{\bf Proportionality Theorem }} \setcounter{equation}{0}
\newcommand\CA{\mathcal A}
\vskip 0.3cm In this section we describe the proportionality
theorem for the Siegel modular variety following the work of
Mumford \cite{Mf3}. Historically F. Hirzebruch \cite{H} first
described a beautiful proportionality theorem for the case of a $
\textit{compact}$ locally symmetric variety in 1956. We shall
state his proportionality theorem roughly. Let $D$ be a bounded
symmetric domain and let $\G$ be a discrete torsion-free
co-compact group of automorphisms of $D$. We assume that the
quotient space $X_\G:=\G\ba D$ is a $\textit{compact}$ locally
symmetric variety. We denote by $\check{D}$ the $ \textit{compact
dual}$ of $D$. Hirzebruch \cite{H} proved that the Chern numbers
of $X_\G$ are proportional to the Chern numbers of $\check{D}$,
the constant of proportionality being the volume of $X_\G$ in a
natural metric. Mumford \cite{Mf3} generalized Hirzebruch's
proportionality theorem to the case of a noncompact arithmetic
variety. \vskip 0.2cm Before we describe the proportionality
theorem for the Siegel modular variety, first of all we review the
compact dual of the Siegel upper half plane $\BH_g$. We note that
$\BH_g$ is biholomorphic to the generalized unit disk $\BD_g$ of
degree $g$ through the Cayley transform (2.7). We suppose that
$\Lambda=(\BZ^{2g},\langle\ ,\ \rangle)$ is a symplectic lattice with a
symplectic form $\langle\ ,\ \rangle.$ We extend scalars of the lattice
$\Lambda$ to $\BC$. Let
\begin{equation*}
{\mathfrak Y}_g:=\left\{\,L\subset \BC^{2g}\,|\ \dim_\BC L=g,\ \
\langle x,y \rangle=0\quad \textrm{for all}\ x,y\in L\,\right\}
\end{equation*}
be the complex Lagrangian Grassmannian variety parameterizing
totally isotropic subspaces of complex dimension $g$. For the
present time being, for brevity, we put $G=Sp(g,\BR)$ and
$K=U(g).$ The complexification $G_\BC=Sp(g,\BC)$ of $G$ acts on
${\mathfrak Y}_g$ transitively. If $H$ is the isotropy subgroup of
$G_\BC$ fixing the first summand $\BC^g$, we can identify
${\mathfrak Y}_g$ with the compact homogeneous space $G_\BC/H.$ We
let
\begin{equation*}
{\mathfrak Y}_g^+:=\big\{\,L\in {\mathfrak Y}_g\,|\ -i \langle x,{\bar
x}\rangle
>0\quad \textrm{for all}\ x(\neq 0)\in L\,\big\}
\end{equation*}
be an open subset of ${\mathfrak Y}_g$. We see that $G$ acts on
${\mathfrak Y}_g^+$ transitively. It can be shown that ${\mathfrak
Y}_g^+$ is biholomorphic to $G/K\cong \BH_g.$ A basis of a lattice
$L\in {\mathfrak Y}_g^+$ is given by a unique $2g\times g$ matrix
${}^t(-I_g\,\,\Om)$ with $\Om\in\BH_g$. Therefore we can identify
$L$ with $\Om$ in $\BH_g$. In this way, we embed $\BH_g$ into
${\mathfrak Y}_g$ as an open subset of ${\mathfrak Y}_g$. The
complex projective variety ${\mathfrak Y}_g$ is called the $
\textit{compact dual}$ of $\BH_g.$

\vskip 0.2cm Let $\G$ be an arithmetic subgroup of $\G_g$. Let
$E_0$ be a $G$-equivariant holomorphic vector bundle over
$\BH_g=G/K$ of rank $n$. Then $E_0$ is defined by the
representation $\tau:K\lrt GL(n,\BC).$ That is, $E_0\cong
G\times_K \BC^n$ is a homogeneous vector bundle over $G/K$. We
naturally obtain a holomorphic vector bundle $E$ over
$\CA_{g,\G}:=\G\ba G/K.$ $E$ is often called an
$\textit{automorphic}$ or $ \textit{arithmetic}$ vector bundle
over $\CA_{g,\G}$. Since $K$ is compact, $E_0$ carries a
$G$-equivariant Hermitian metric $h_0$ which induces a Hermitian
metric $h$ on $E$. According to Main Theorem in \cite{Mf3}, $E$
admits a $ \textit{unique}$ extension ${\tilde E}$ to a smooth
toroidal compactification ${\tilde \CA}_{g,\G}$ of $\CA_{g,\G}$
such that $h$ is a singular Hermitian metric $ \textit{good}$ on
${\tilde \CA}_{g,\G}$. For the precise definition of a
$\textit{good metric}$ on $\CA_{g,\G}$ we refer to \cite[p.\,242]{Mf3}.
According to Hirzebruch-Mumford's Proportionality
Theorem\,(cf.\,\cite[p.\,262]{Mf3}), there is a natural metric on
$G/K=\BH_g$ such that the Chern numbers satisfy the following
relation
\begin{equation}
c^{\al}\big({\tilde E}\big)=(-1)^{{\frac 12}g(g+1)}\,
\textmd{vol}\left( \G\ba \BH_g\right)\,c^{\al}\big( {\check
E}_0\big)
\end{equation}
for all $\al=(\al_1,\cdots,\al_n)$ with nonegative integers
$\al_i\,(1\leq i\leq n)$ and $\sum_{i=1}^n\al_i={\frac 12}g(g+1),$
where ${\check E}_0$ is the $G_{\BC}$-equivariant holomorphic
vector bundle on the compact dual ${\mathfrak Y}_g$ of $\BH_g$
defined by a certain representation of the stabilizer $
\textrm{Stab}_{G_\BC}(e)$ of a point $e$ in ${\mathfrak Y}_g$.
Here $\textmd{vol}\left( \G\ba \BH_g\right)$ is the volume of
$\G\ba\BH_g$ that can be computed\,(cf.\,\cite{Si1}).

\vskip 0.2cm \noindent $ \textbf{Remark 11.1.}$ Goresky and Pardon
\cite{Go} investigated Chern numbers of an automorphic vector
bundle over the Baily-Borel compactification ${\overline X}$ of a
Shimura variety $X$. It is known that ${\overline X}$ is usually a
highly singular complex projective variety. They also described
the close relationship between the topology of $X$ and the
characteristic classes of the unique extension ${\widetilde{TX}}$
of the tangent bundle $ {TX}$ of $X$ to a smooth toroidal
compactification ${\tilde X}$ of $X$.

\end{section}

%
%
\begin{section}{{\bf Motives and Siegel Modular Forms}}
\setcounter{equation}{0} \vskip 0.3cm Assuming the existence of
the hypothetical motive $M(f)$ attached to a Siegel modular form
$f$ of degree $g$, H. Yoshida \cite{Yo2} proved an interesting
fact that $M(f)$ has at most $g+1$ period invariants. I shall
describe his results in some detail following his paper. \vskip
0.1cm Let $E$ be an algebraic number field with finite degree
$l=[E:\BQ].$ Let $J_E$ be the set of all isomorphisms of $E$ into
$\BC.$ We put $R=E\otimes_\BQ \BC.$ Let $M$ be a motive over $\BQ$
with coefficients in $E$. Roughly speaking motives arise as direct
summands of the cohomology of a smooth projective algebraic
variety defined over $\BQ.$ Naively they may be defined by a
collection of realizations satisfying certain axioms. A motive $M$
has at least three realizations : the Betti realization, the de
Rham realization and the $\la$-adic realization. \par First we let
$H_B(M)$ be the Betti realization of $M$. Then $H_B(M)$ is a free
module over $E$ of rank $d:=d(M).$ We put
$H_B(M)_\BC:=H_B(M)\otimes_\BQ \BC$. We have the involution
$F_{\infty}$ acting on $H_B(M)_\BC\ E$-linearly. Therefore we
obtain the the eigenspace decomposition
\begin{equation}
H_B(M)_\BC=H_B^+(M)\oplus H_B^-(M),
\end{equation}
where $H_B^+(M)$ (resp. $H_B^-(M)$) denotes the $(+1)$-eigenspace
(resp. the $(-1)$-eigenspace) of $H_B(M)$. We let $d^+$\,(resp.
$d^-$) be the dimension $H_B^+(M)$\,(resp. $H_B^-(M)$).
Furthermore $H_B(M)_\BC$ has the Hodge decomposition into
$\BC$-vector spaces\,:
\begin{equation}
H_B(M)_\BC=\bigoplus_{p,q\in\BZ} H^{p,q}(M),
\end{equation}
where $H^{p,q}(M)$ is a free $R$-module. A motive $M$ is said to
be $ \textit{of pure weight}\ w:=w(M)$ if $H^{p,q}(M)=\{0\}$
whenever $p+q\neq w.$ From now on we shall assume that $M$ is of
pure weight.\par Secondly we let $H_{\textrm{DR}}(M)$ be the de
Rham realization of $M$ that is a free module over $E$ of rank
$d$. Let
\begin{equation}
H_{\textrm{DR}}(M)=F^{i_1}\supsetneqq F^{i_2}\supsetneqq
\cdots\supsetneqq F^{i_m}\supsetneqq F^{i_{m+1}}=\big\{ 0\big\}
\end{equation}
be a decreasing Hodge filtration so that there are no different
filtrations between successive members. The choice of members
$i_{\nu}$ may not be unique for $F^{i_\nu}$. For the sake of
simplicity, we assume that $i_{\nu}$ is chosen for $1\leq \nu\leq
m$ so that it is the maxium number. We put
\begin{equation*}
s_\nu= \textrm{rank}\ H^{i_\nu ,w-i_\nu}(M),\qquad 1\leq \nu\leq
m,
\end{equation*}
where rank means the rank as a free $R$-module. Let
\begin{equation*}
I:H_B(M)_\BC\lrt
H_{\textrm{DR}}(M)_\BC=H_{\textrm{DR}}(M)\otimes_E \BC
\end{equation*}
be the comparison isomorphism which satisfies the conditions
\begin{equation}
I\left( \bigoplus_{p'\geq p} H^{p',q}(M)\right)=F^p\otimes_\BQ
\BC.
\end{equation}
According to $(12.4)$, we get
\begin{equation*}
s_\nu= \dim_E\,F^{i_\nu}-\dim_E\,F^{i_{\nu+1}},\quad
\dim_E\,F^{i_\nu}=s_\nu+s_{\nu+1}+\cdots+s_m,\quad 1\leq \nu\leq
m.
\end{equation*}
We choose a basis $\big\{ w_1,\cdots,w_d\big\}$ of
$H_{\textrm{DR}}(M)$ over $E$ so that $\big\{
w_{s_1+s_2+\cdots+s_{\nu-1}+1},\cdots,w_d\big\}$ is a basis of
$F^{i_\nu}$ for $1\leq \nu\leq m.$ We observe that
\begin{equation}
d=s_1+s_2+\cdots+s_m\quad \textrm{all}\ s_\nu>0\ \textrm{with}\
1\leq \nu\leq m.
\end{equation}
We are in a position to describe the fundamental periods of $M$
that Yoshida introduced. Let $\big\{
v_1^+,v_2^+,\cdots,v_{d^+}^+\big\}$\,$\big($resp.\,$\big\{
v_1^-,v_2^-,\cdots,v_{d^-}^-\big\}\,\big)$ be a basis of
$H_B^+(M)$ (resp. $H_B^-(M)$) over $E$. Writing
\begin{equation}
I(v_j^\pm)=\sum_{i=1}^d x_{ij}^\pm w_i,\quad x_{ij}^\pm\in R,\quad
1\leq j\leq d^\pm,
\end{equation}
we obtain a matrix $X^+=\big(x_{ij}^+\big)\in R^{(d,d^+)}$ and a
matrix $X^-=\big(x_{ij}^-\big)\in R^{(d,d^-)}$. We recall that
$R^{(m,n)}$ denotes the set of all $m\times n$ matrices with
entries in $R$. Let $P_M$ be the lower parabolic subgroup of
$GL(d)$ which corresponds to the partition (12.5). Let $P_M(E)$ be
the group of $E$-rational points of $P_M$. Then the coset of
$X^+$\,(resp.\,$X^-$) in
\begin{equation*}
P_M(E)\ba R^{(d,d^+)}/GL(d^+,E)\quad \big( \textrm{resp.}\
P_M(E)\ba R^{(d,d^-)}/GL(d^-,E)\big)
\end{equation*}
is independent of the choice of a basis. We set $X_M=(X^+,X^-)\in
R^{(d,d)}.$ Then it is easily seen that the coset of $X_M$ in
\begin{equation*}
P_M(E)\ba R^{(d,d)}/\big(GL(d^+,E)\times GL(d^-,E)\big)
\end{equation*}
is independent of the choice of a basis, i.e., well defined. A
$d\times d$ matrix $X_M=(X^+,X^-)$ is called a $ \textit{period
matrix}$ of $M$.\par For an $m$-tuple $(a_1,\cdots,a_m)\in \BZ^m$
of integers, we define a character $\la_1$ of $P_M$ by
\begin{equation*}
\lambda _1 \left(
\begin{pmatrix} P_{1} & 0 & \hdots & 0 \\
\ast & P_{2} & \hdots & 0 \\
\ast & \ast & \ddots & \vdots \\
\ast & \ast & \ast & P_{m}
\end{pmatrix} \right)
=\prod_{j=1}^m \det (P_{j})^{a_j},\quad P_j\in GL(s_j),\quad 1\leq
j\leq m.
\end{equation*}
For a pair $(k^+,k^-)$ of integers, we define a character $\la_2$
of $GL(d^+)\times GL(d^-)$ by
\begin{equation*}
\lambda _2 \left ( \begin{pmatrix} A & 0 \\ 0 & B \end{pmatrix}
\right) =(\det A)^{k^+} (\det B)^{k^-}, \qquad A \in GL(d^+), \ B
\in GL(d^-).
\end{equation*}
A polynomial $f$ on $R^{(d,d)}$ rational over $\BQ$ is said to be
$\textit{of the type}\ \big\{ (a_1,\cdots,a_m);(k^+,k^-)\big\}$ or
$ \textit{of the type}\ (\la_1,\la_2)$ if $f$ satisfies the
following condition
\begin{equation}
f(pxq )=\lambda _1(p)\lambda _2(q )f(x) \qquad \text {for all \ $p
\in P_M$, $q \in GL(d^+) \times GL(d^-)$}.
\end{equation}
We now assume that $f$ is a nonzero polynomial on $R^{(d,d)}$ of
the type $\big\{ (a_1,\cdots,a_m);(k^+,k^-)\big\}$. Let
$X_M=(X^+,X^-)$ be a period matrix of a motive $M$ as before. Then
it is clear that $f(X_M)$ is uniquely determined up to
multiplication by elements in $E^{\times}$. We call $f(X_M)$ a $
\textit{period invariant}$ of $M$ of the type $\big\{
(a_1,\cdots,a_m);(k^+,k^-)\big\}$. Hereafter we understand the
equality between period invariants mod $E^\times$.\par We now
consider the following special polynomials of the type
$(\la_1,\la_2):$\par \noindent $ \textbf{I.}$ Let $f(x)=\det (x)$
for $x\in R^{(d,d)}.$ \par It is easily seen that $f(x)$ is of the
type $\big\{ (1,1,\cdots,1);(1,1)\big\}$. Then $f(X_M)$ is nothing
but Deligne's period $\delta(M)$.
\par \noindent $ \textbf{II.}$ Let $f^+(x)$ be the determinant of
the upper left $d^+\times d^+$-submatrix of $x\in R^{(d,d)}$. It
is easily checked that $f^+(x)$ is of the type
\begin{equation*}
\big\{ (\overbrace{1,1, \ldots ,1}^{p^+}, 0,\ldots ,0); (1,0)
\big\},
\end{equation*} where $p^+$ is a positive integer such that
$s_1+s_2+\cdots+s_{p^+}=d^+.$ We note that $f^+(X_M)$ is Deligne's
period $c^+(M)$.
\par \noindent $ \textbf{III.}$
Let $f^-(x)$ be the determinant of the upper right $d^- \times
d^-$-submatrix of $x$. Then $f^-(x)$ is of the type
\begin{equation*}
\big\{ (\overbrace{1,1, \ldots ,1}^{p^-}, 0,\ldots ,0); (0,1)
\big\}
\end{equation*} and $f^-(X_M)$ is Deligne's period $c^-(M)$. Here
$p^-$ is a positive integer such that
$s_1+s_2+\cdots+s_{p^-}=d^-.$
\par
Either one of the above conditions is equivalent to that
$F^{\mp}(M)$, hence also $c^{\pm}(M)$ can be defined (cf.
\cite{De}, \S1, \cite{Yo1}, \S2). We have
$F^{\mp}(M)=F^{i_{p^{\pm}+1}}(M)$; $F^{\pm}(M)$ can be defined if
$M$ has a critical value. Let $\Cal P=\Cal P(M)$ denote the set of
integers $p$ such that $s_1+s_2+ \cdots +s_p<\min (d^+, d^-)$.
Yoshida (cf.\,\cite{Yo2}, Theorem 3) showed that for every $p \in
\mathcal P$, there exists a non-zero polynomial $f_p$ of the type
\begin{equation*}
\big\{(\overbrace {2, \ldots , 2}^p, \overbrace {1, \ldots ,
1}^{m-2p}, \overbrace {0, \ldots , 0}^p); (1,1) \big\}
\end{equation*}
and that every polynomial satisfying (12.7) can be written
uniquely as a monomial of $\det (x)$, $f^+(x)$, $f^-(x)$,
$f_p(x)$, $p \in \Cal P$. We put $c_p(M)=f_p(X_M)$. We call
$\delta (M)$, $c^{\pm}(M)$, $c_p(M)$, $p \in \mathcal P$ the {\it
fundamental periods} of $M$. Therefore any period invariant of $M$
can be written as a monomial of the fundamental periods. Moreover
Yoshida showed that if a motive $M$ is constructed from motives
$M_1,\cdots,M_t$ of pure weight by standard algebraic operations
then the fundamental periods of $M$ can be written as monomials of
the fundamental periods of $M_1,\cdots,M_t.$ He proved that a
motive $M$ has at most $ \textrm{min}(d^+,d^-)+2$ fundamental
periods including Deligne's periods $\delta(M)$ and $c^\pm(M)$.
\par
Thirdly we let $H_\la (M)$ be the $\la$-adic realization of $M$.
We note that $H_\la (M)$ is a free module over $E_{\la}$ of rank
$d$. We have a continuous $\la$-adic representation of the
absolute Galois group $G_\BQ= \textrm{Gal}\big({\overline
\BQ}/\BQ\big)$ on $H_\la (M)$ for each prime $\la$. Also there is
an isomorphism $I_\la:H_B(M)\otimes_E E_\la\lrt H_\la (M)$ which
transforms the involution $F_{\infty}$ into the complex
conjugation.


\vskip 0.1cm We recall that an integer $s=n$ is said to be $
\textit{critical}$ for a motive $M$ if both the infinite Euler
factors $L_{\infty}(M,s)$ and $L_{\infty}({\check M},s)$ are
holomorphic at $s=n$. Here $L(M,s)$ denotes the complex
$L$-function attached to $M$ and ${\check M}$ denotes the dual
motive of $M$. Such values $L(M,n)$ are called $ \textit{critical
values }$ of $L(M,s).$ Deligne proposed the following. \vskip
0.1cm\noindent $ \textbf{Conjecture\,(Deligne \cite{De})}.$ {\it
Let $M$ be a motive of pure weight and $L(M,s)$ the $L$-function
of $M$. Then for critical values $L(M,n)$, one has}
\begin{equation*}
{ {L(M,n)}\over { (2\pi i)^{d^\pm}\,c^\pm(M) } }\in E,\quad
d^\pm:=d^\pm (M),\ \pm 1=(-1)^n.
\end{equation*}
Indeed Deligne showed that $c^\pm(M)\in R^\times$ and Yoshida
showed that other period invariants are elements of $R^\times$.
\vskip 0.1cm \noindent $ \textbf{Remark 12.1.}$ The Hodge
decomposition (12.2) determines the gamma factors of the
conjectural functional equation of $L(M,s)$. Conversely the gamma
factor of the functional equation of $L(M,s)$ determines the Hodge
decomposition if $M$ is of pure weight.

\vskip 0.2cm Let $f\in S_k(\G_g)$ be a nonzero Hecke eigenform on
$\BH_g$. Let $L_{ \textrm{st}}(s,f)$ and $L_{ \textrm{sp}}(s,f)$
be the standard zeta function and the spinor zeta function of $f$
respectively. For the sake of simplicity we use the notations $L_{
\textrm{st}}(s,f)$ and $L_{ \textrm{sp}}(s,f)$ instead of $D_f(s)$
and $Z_f(s)$ (cf.\ \S8) in this section. We put $w=kg-{\frac
12}g(g+1)$. We have a normalized Petersson inner product $\langle\
,\ \rangle$ on $S_k(\G_g)$ given by
\begin{equation*}
\langle F,F\rangle= \textrm{vol}\big(\G_g\ba
\BH_g\big)^{-1}\,\int_{\G_g\ba \BH_g} |f(\Om)|^2\,\big(\det
Y\big)^{k-g-1}[dX][dY],\quad F\in S_k(\G_g),
\end{equation*}
where $\Om=X+iY\in \BH_g$ with real $X=(x_{\mu\nu}),\
Y=(y_{\mu\nu})$, $[dX]=\bigwedge_{\mu\leq\nu}dx_{\mu\nu}$ and
$[dY]=\bigwedge_{\mu\leq\nu}dy_{\mu\nu}$.

\vskip 0.1cm We assume the following (A1)-(A6)\,:\par\noindent
 {\bf (A1)} The Fourier coefficients of $f$ are contained in a totally
real algebraic number field $E$.\par\noindent {\bf (A2)} There
exist motives $M_{ \textrm{st}}(f)$ and $M_{ \textrm{sp}}(f)$ over
$\BQ$ with coefficients in $E$ satisfying the conditions
\begin{equation*}
L\big(M_{ \textrm{st}}(f),s)= \big( L_{
\textrm{st}}(s,f^\s)\big)_{\s\in J_E}\quad \textrm{and}\quad
L\big(M_{ \textrm{sp}}(f),s)= \big( L_{
\textrm{sp}}(s,f^\s)\big)_{\s\in J_E}.
\end{equation*}
\par\noindent {\bf (A3)} Both $M_{ \textrm{st}}(f)$ and $M_{
\textrm{sp}}(f)$ are of pure weight.

\par\noindent {\bf (A4)} We assume
\begin{equation*}
\bigwedge^{2g+1}M_{ \textrm{st}}(f)\cong T(0),
\end{equation*}
\begin{eqnarray*}
H_B(M_{\textrm{st}}(f)) \otimes _{\BQ} \BC
=& H^{0,0}(M_{\text {st}}(f)) \hskip 5cm\\
\bigoplus _{i=1}^g & \Big(H^{-k+i,k-i}(M_{\textrm{st}}(f)) \oplus
H^{k-i,-k+i}(M_{\textrm{st}}(f))\Big).
\end{eqnarray*}
We also assume that the involution $F_\infty$ acts on
$H^{0,0}(M_{\text {st}}(f))$ by $(-1)^g.$

\par\noindent {\bf (A5)} We assume
\begin{equation*}
\bigwedge ^{2^g} M_{\textrm{sp}}(f) \cong T(2^{g-1}w),
\end{equation*}
\begin{equation*}
H_B(M_{\textrm{sp}}(f)) \otimes _{\BQ} \BC=\bigoplus _{p, q}
H^{p,q}(M_{\textrm{sp}}(f)),
\end{equation*}

\begin{eqnarray*}
& p=(k-i_1)+(k-i_2)+ \cdots + (k-i_r), \quad
q=(k-j_1)+(k-j_2)+ \cdots + (k-j_s), \\
& r+s=g, \qquad 1 \leq i_1< \cdots <i_r \leq g, \qquad
1 \leq j_1<\cdots <j_s \leq g, \\
& \{ i_1, \cdots , i_r \} \cup \{ j_1, \cdots , j_s \} =\{ 1, 2,
\ldots , g \},
\end{eqnarray*}
including the cases $r=0$ or $s=0.$

\par\noindent {\bf (A6)} If
$w=kg-{\frac 12}g(g+1)$ is even, then the eigenvalues $+1$ and
$-1$ of $F_{\infty}$ on $H^{p,p}(M_{\text {sp}}(f))$ occur with
the equal multiplicities.

\vskip 0.2cm Let $J_E=\{ \sigma _1, \sigma _2, \ldots , \sigma _l
\}$, $l=[E:\BQ]$ and write $x \in R \cong \BC^{J_E}$ as
$x=(x^{(1)}, x^{(2)}, \cdots , x^{(l)})$, $x^{(i)} \in \BC$ so
that $x^{(i)}=x^{\sigma _i}$ for $x \in E$. Yoshida showed that
when $k>2g$, assuming Deligne's conjecture, one has
\begin{equation*}
c^{\pm}(M_{\text {st}}(f))=\pi ^{kg} \big(\langle f^{\sigma},
f^{\sigma} \rangle \big)_{\sigma \in J_E}.
\end{equation*}
He proved the following interesting
result\,(cf.\,Yoshida\,\cite{Yo2},\,Theorem 14).
\begin{theorem}
Let the notation be the same as above. We assume that two motives
over $\BQ$ having the same $L$-function are isomorphic (Tate's
conjecture). Then there exist $p_1$, $p_2$, $\cdots$, $p_r \in
\BC^{\times}$, $1 \leq r \leq g+1$ such that for any fundamental
period $c \in R^{\times}$ of $M_{\text {st}}(f)$ or $M_{\text
{sp}}(f)$, we have
\begin{equation*}
c^{(1)}=\alpha\,\pi ^A\,p_1^{a_1}p_2^{a_2} \cdots p_r^{a_r}
\end{equation*}
with $\alpha \in \overline {\BQ}^{\times}$ and non-negative
integers $A$, $a_i$, $1 \leq i \leq r$.
\end{theorem}
\vskip 0.1cm \noindent $ \textbf{Remark 12.2.}$ It is widely
believed that the zeta function of the Siegel modular variety
${\mathcal A}_g:=\G_g\ba \BH_g$ can be expressed using the spinor
zeta functions of (not necessarily holomorphic) Siegel modular
forms:
\begin{equation*}
\zeta \big(s, {\mathcal A}_g\big) \fallingdotseq\prod _f L_{\text
{sp}}(s, f).
\end{equation*}
Yoshida proposed the following conjecture. \vskip 0.1cm\noindent $
\textbf{Conjecture\,(Yoshida \cite{Yo2})}.$ {\it If one of two
motives $M_{\text {st}}(f)$ and $M_{\text {sp}}(f)$ is not of pure
weight, then the associated automorphic representation to $f$ is
not tempered. Furthermore $f$ can be obtained as a lifting from
lower degree forms.}

\end{section}

%
%
\begin{section}{{\bf Remark on Cohomology of a Shimura Variety}}
\setcounter{equation}{0} \vskip 0.2cm First we recall the
definition of a Shimura variety. Let $(G,X)$ be a Shimura datum as
in \cite[p.\,322]{Mil} or \cite[pp.\,321-322]{Hi} so that $G$ is a
connected reductive group over $\BQ$ and $X$ is a finite disjoint
union of symmetric Hermitian domains, homogeneous under $G(\BR)$.
The points of $X$ correspond to homomorphism
\begin{equation*}
h_x:{\mathbb S}_\BR= \textrm{Res}_{\BC/\BR}({\mathbb G}_m)\lrt
G_\BR=G\times_\BQ \BR
\end{equation*}
satisfying the following axioms (D.1) and (D.2). For convenience,
we list the axioms for a datum
$(G,X)$\,(cf.\,\cite[p.\,322]{Mil})\,:
\par (D.1) for each $x\in X$, the Hodge structure on ${\mathfrak
g}$ defined by $h_x$ is of type $\big\{ (-1,1),(0,0),$ \newline
$(1,-1)\big\}$;
\par (D.2) for each $x,\ \textrm{ad}\,h_x(i)$ is a Cartan
involution on $G_\BD^{ \textrm{ad}}$;
\par (D.3) $G^{ \textrm{ad}}$ has no factor defined over $\BQ$
whose real points form a compact group\,;
\par (D.4) the identity component of $G(G)^0$ of the center
$Z(G)$ of $G$ splits over a CM-field. \vskip 0.1cm\noindent Here
${\mathfrak g}$ denotes the Lie algebra of $G$ and $G^{
\textrm{ad}}$ denotes the adjoint group of $G$. Axiom (D.4) is not
in Deligne's list of axioms\,[in {\em Variet{\'e}s de Shimura\,:
interpretation modulaire, et techniques de construction de
mod{\'e}les canoniques,} Proc. Symp. Pure Math., A.M.S. {\bf 33},
Part 2 (1979), 247--290]. For an open compact subgroup $K$ of $G(
{\mathbb A}_f)$, we consider
\begin{equation*}
 Sh_K(G,X):=G(\BQ)\ba X\times G( {\mathbb A}_f)/K,
\end{equation*}
where
\begin{equation*}
q(x,a)k=(qx,qak),\quad q\in G(\BQ),\ x\in X,\ a\in G( {\mathbb
A}_f),\ k\in K.
\end{equation*}
Endowed with the quotient topology this is a Hausdorff space with
finitely many connected components, each of which is isomorphic to
$\G\ba X^+$ for any connected component $X^+$ of $X$ and some
arithmetic subgroup $\G\subset G(\BQ)$. The space $Sh_K(G,X)$ is a
quasi-projective complex algebraic variety. $Sh_K(G,X)$ has a
canonical model over the reflex field $E(G,X)$. This is a normal
quasi-projective scheme over $E(G,X)$ together with an isomorphism
between the complex space of its $\BC$-valued points and
$Sh_K(G,X)$. Let
\begin{equation*}
 Sh(G,X):=\lim_{\longleftarrow\atop K}Sh_K(G,X).
\end{equation*}
Then this is a scheme over $\BC$ whose complex points are
\begin{equation*}
 Sh(G,X):=G(\BQ)\ba X\times G( {\mathbb A}_f)/Z(\BQ)^-,
\end{equation*}
where $Z(\BQ)^-$ is the closure of $Z(\BQ)$ in $Z({\mathbb A}_f)$.
One has a natural continuous action of $G( {\mathbb A}_f)$ on
$Sh(G,X)$ given by
\begin{equation}
[x,a]h=[x,ah],\quad x\in X,\ a,h\in G( {\mathbb A}_f).
\end{equation}
The scheme $Sh(G,X)$ together with the above action (13.1) is
called the $ \textit{Shimura variety}$ defined by $(G,X)$. By a $
\textit{model}$ of $Sh(G,X)$ over a subfield $E$ of $\BC$, we mean
a scheme $S$ defined over $\BQ$ together with a $E$-rational
action of $G( {\mathbb A}_f)$ such that there is a $G( {\mathbb
A}_f)$-equivariant isomorphism (over $\BC$)
\begin{equation*}
 Sh(G,X) \cong S\otimes_E \BC.
\end{equation*}
\vskip 0.1cm \noindent {\bf Example 13.1.} Let $G=GSp(2)$ be the
group of symplectic similitudes of degree $2$. We fix a morphism
\begin{equation*}
h_0:{\mathbb S}_\BR\lrt G_\BR=G\times_\BQ \BR
\end{equation*}
by requiring that
\begin{equation*}
\BC^\times={\mathbb S}_\BR\ni x+iy\mapsto \begin{pmatrix} xI_2 &
yI_2\\ -y I_2 & xI_2 \end{pmatrix}\in GSp(2,\BR).
\end{equation*}
The $G(\BR)$-conjugacy class of the homomorphism $h_0$ is
analytically isomorphic to the union $X^\pm:=X^+\cup X^-$ of the
Siegel upper and lower half planes of degree two. The pair
$(GSp_2,X^\pm)$ defines a Shimura variety $Sh(GSp_2,X^\pm)$ which
is the Siegel modular variety of degree two. Its reflex filed is
$\BQ$. \vskip 0.1cm \noindent {\bf Example 13.2.} Let $F$ be a
totally real number field and let $G=GL(2,F)$ so that
$G(\BR)=\prod_{ \textrm{Hom}(F,\BR)}GL(2,\BR).$ Let $X$ be the set
of $G(\BR)$-conjugates of $h_0:{\mathbb S}_\BR\lrt
G_\BR=G\times_\BQ \BR$ given by
\begin{equation*}
h_0(a+ib)=\left( \begin{pmatrix} a & -b \\ b & \
a\end{pmatrix},\begin{pmatrix} a & -b \\ b & \
a\end{pmatrix},\cdots,\begin{pmatrix} a & -b \\ b & \
a\end{pmatrix}\right),\quad a+ib\in\BC \,(a,b\in \BR).
\end{equation*}
Then $X$ is a point of $[E:\BQ]$ copies of $\BC-\BR$, and $(G,X)$
satisfies the axioms (D.1)--(D.2). The Shimura variety is nothing
but the so-called {\it Hilbert modular variety.}

\vskip 0.2cm It is well known that every Shimura variety $Sh(G,X)$
has a unique canonical model $Sh(G,X)_E$ over the reflex field
$E(G,X)$\,(cf.\,\cite[Theorem 5.5]{Mil}). $Sh(G,X)$ may be viewed
as a parameter space for a family of
motives\,(cf.\,\cite[II.\,\S3]{Mil}). The determination of the
zeta function of a Shimura variety and its expression in terms of
$L$-functions of automorphic representations is viewed by
Langlands as a higher dimensional version of Artin reciprocity. It
is known that there exists a smooth compactification of a Shimura
variety $Sh(G,X)$\,(cf.\,\cite[p.\,395]{Mil}). There is a notion
of a $ \textit{mixed Shimura variety}$ which we will not consider
here. Goresky and Pardon \cite{Go} showed that the minimal
compactification of a Shimura variety resembles a smooth
projective variety in that one can define its Chern classes in
cohomology with complex coefficients and that one can also define
the Chern classes of an automorphic vector bundle as cohomology
classes on the minimal compactification with complex coefficients.
Zucker\,\cite{Zu2} treated several topological
compactifications\,(e.g., the reductive Borel-Serre
compactification, the minimal compactification) of a Shimura
variety as algebraic varieties and constructed mixed Hodge
structures on their cohomology groups. In \cite{Zu1} Zucker proved
that the $L^p$-cohomology of a certain Shimura variety is
canonically isomorphic to the ordinary cohomology of its reductive
Borel-Serre compactification. Cojectures of Beilinson and Deligne
predict the existence of extensions of mixed motives which should
be become visible in the cohomology of open varieties over number
fields. Harris and Zucker\,\cite{HZ1} proved that mixed Hodge-De
Rham structures arise naturally in the boundary cohomology of
automorphic vector bundles on Shimura varieties. The theory of
automorphic vector bundles on Shimura varieties has been studied
by Harris et al (cf.\,\cite{Har1,Har2},\,\cite{Go}). For the
theory of the cohomology of arithmetic varieties we also refer to
\cite{S1,S2,S3} and \cite{LS}. R. Charney and R. Lee
\cite{C-R} showed that the stable cohomology of the Satake
cohomology ${\bar {\mathcal A}}_g$ of ${\mathcal A}_g$ contains a polynomial algebra
which coincides with the stable cohomology of the compact dual
${\mathfrak Y}_g$ of $\BH_g$. We note that the intersection
cohomology $IH^*\big({\bar {\mathcal A}}_g,\BC)$ contains a copy of
$H^*({\mathfrak Y}_g,\BC).$
For the theory of the cohomology of
the Siegel modular variety (in particular of degree two) we refer
to \cite{C-R},\,\cite{HW1,HW2},\,\cite{L},\,\cite{LW1} and \cite{Sch1,Sch2}.
\par
Labesse and Schwermer \cite{Lab} used two kinds of lifting of an
irreducible automorphic representation $\pi$ of $GL(2,{\mathbb
A}_F)$ with $F$ a global field to obtain the nonvanishing of
certain cusp cohomology classes which correspond to nondiscrete
representations. One lift is the lifting of $\pi$ to a
representation of $GL(2,{\mathbb A}_K)$ introduced by Langlands
\cite{Lan}, where $K$ is a cyclic extension of $F$ of prime degree
or a cubic extension of $F$. The other is the so called
Gelbart-Jacquet lifting of $\pi$ to a representation of
$GL(3,{\mathbb A}_F)$ via the adjoint representation
(cf.\,\cite{GJ}).

\end{section}

%
%


%
%


%
%

\end{document}